\documentclass[a4paper,10pt]{article}
\usepackage[utf8x]{inputenc}
\usepackage{tracefnt,amsmath,tabu,array}
\usepackage{amssymb,graphicx,setspace,amsfonts,amsbsy}
\usepackage{pifont,latexsym,ifthen,amsthm,rotating,calc,textcase,booktabs,cancel,slashed}

\newtheorem{theorem}{Theorem}[section]
\newtheorem{lemma}[theorem]{Lemma}
\newtheorem{corollary}[theorem]{Corollary}
\newtheorem{proposition}[theorem]{Proposition}
\newtheorem{definition}[theorem]{Definition}

\newtheorem{remark}[theorem]{Remark}
\newcommand{\filledbox}{\leavevmode
  \hbox to.77778em{%
  \hfil\vbox to.675em{\hrule width.6em height.6em}\hfil}}

\newcommand{\Rm}{{\mathbb R}}

\begin{document}
%\doublespacing
\tabulinesep=1.0mm
% Enter full title and short title for running headers
\title{Non-radiative solutions and long-time dynamics of 5D focusing energy-critical wave equation in the radial case}

\author{Ruipeng Shen\\
Centre for Applied Mathematics\\
Tianjin University\\
Tianjin, China
}

\maketitle

\begin{abstract}
 In this article we discuss the long-time dynamics of the radial solutions to the focusing energy-critical wave equation in 5-dimensional space. We give some details about the asymptotic behaviour, topological structure and time evolution of the non-radiative solutions to this equation. As an application we prove a quantitative version of soliton resolution theorem for solutions defined for all time $t>0$, which immediately verifies the soliton resolution conjecture in the radial case, without a priori boundedness assumption on the energy norm of solution as time tends to infinity. The main tool of this work is the radiation theory of wave equations and the major observation of this work is a correspondence between the radiation and the soliton collision behaviour of solutions. 
\end{abstract}

\section{Introduction} 

\subsection{Background} 

In this work we consider the long-time behaviour of the radial solutions to the focusing, energy critical wave equation in 5-dimensional space 
\[
 \left\{\begin{array}{ll} \partial_t^2 u - \Delta u = |u|^{4/3} u, & (x,t) \in \Rm^5 \times \Rm;  \\ (u,u_t)|_{t=0} = (u_0,u_1)\in \dot{H}^1\times L^2. & \end{array} \right. \qquad \hbox{(CP1)}
\]
For convenience we use the notation $F(u) = |u|^{4/3} u$ in this work. The energy is conserved for all $t$ in the maximal lifespan $(-T_-,T_+)$:
\[
 E = \int_{\Rm^5} \left(\frac{1}{2}|\nabla u(x,t)|^2 + \frac{1}{2}|u_t(x,t)|^2 - \frac{3}{10}|u(x,t)|^{10/3}\right) {\rm d} x.
\]
This equation is invariant under the natural dilation. More precisely, if $u$ is a solution to (CP1), then 
\[ 
 u_\lambda = \frac{1}{\lambda^{3/2}} u\left(\frac{x}{\lambda}, \frac{t}{\lambda}\right), \qquad \lambda\in \Rm^+
\]
is also a solution to (CP1). This equation is called energy critical since the initial data of $u$ and $u_\lambda$ share the same $\dot{H}^1\times L^2$ norm and energy. 

Unlike the $d$-dimensional defocusing energy-critical wave equation $\partial_t^2 u - \Delta u = - |u|^{\frac{4}{d-2}} u$, in which all finite-energy solutions are defined for all $t\in \Rm$ and scatter in both two time directions(see, \cite{mg1, enscatter1, enscatter2, ss1, ss2}, for instance), the long time behaviour of solutions in the focusing case are quite complicated and subtle. We give a few examples: 

\paragraph{Finite time blow-up} If the solution blows up at time $T_+ \in \Rm^+$, then we must have 
\[
 \|u\|_{L^{7/3} L^{14/3} ([0,T_+)\times \Rm^3)} = + \infty. 
\]
We may further divide finite time blow-up solutions into two types:
\begin{itemize} 
 \item Type I blow-up solutions satisfy
  \[
   \limsup_{t\rightarrow T_+} \|(u,u_t)\|_{\dot{H}^1\times L^2} = + \infty. 
  \]
 An explicit example can be given by
  \[
   u(x,t) = \left(\frac{15}{4}\right)^{3/4} (T_+ -t)^{-3/2}.
  \]
 A smooth cut-off technique and the finite speed of propagation then gives a type I blow-up solution with initial data in the energy space. It has been proved in Donninger \cite{stablet1} that the type I blow-up of this example is stable under a small perturbation in the energy space. 
 \item Type II blow-up solution satisfies 
 \[
  \limsup_{t\rightarrow T_+} \|(u,u_t)\|_{\dot{H}^1\times L^2} < +\infty. 
 \]
 This kind of solutions was first constructed in 3-dimensional space by Krieger-Schlag-Tataru \cite{slowblowup1}, Krieger-Schlag \cite{slowblowup2} and Donninger-Huang-Krieger-Schlag \cite{moreexamples}. Then similar type II blow-up solutions were given by Hillairet-Rapha\"{e}l \cite{4dtypeII} in 4D and Jendrej \cite{5dtypeII} in 5D. %The behaviour of these solutions as $t\rightarrow T_+$ will be introduced in the soliton resolution part below. The instability and the stable manifolds of the specific examples given in the first two papers above have also been discussed by Krieger \cite{stablem2}, Krieger-Nahas \cite{instable} and Burzio-Krieger \cite{stablemanit2}. 
\end{itemize}

\paragraph{Global solutions} Global solutions are defined for all $t\in \Rm^+$. One typical example is the scattering solution, whose asymptotic behaviour  resembles that of a free wave. More precisely, $u$ scatters in the positive time direction if and only if there exists a linear free wave $v^+$, such that 
\[
 \lim_{t\rightarrow +\infty} \|\vec{u}(t) - \vec{v}^+(t)\|_{\dot{H}^1\times L^2} = 0. 
\]
Here the notations $\vec{u} = (u,u_t)$ and $\vec{v}^+ = (v^+, v_t^+)$ will be frequently used in this work. In particular, a solution with small initial data in the energy space must be a scattering solution. Another important example of global solution is the ground state
\[
 W(x) = \left(1 + \frac{|x|^2}{15}\right)^{-3/2}. 
\]
This comes with a smallest energy among all non-trivial stationary solutions of (CP1), i.e. solutions to the elliptic equation $-\Delta u = F(u)$. In fact, all radial finite-energy stationary solutions are exactly given by (see Section 4)
\[
 \left\{0\right\}\cup\left\{\pm W_\lambda: \lambda \in \Rm^+\right\}. 
\]
Here $W_\lambda(x)$ is the rescaled version of $W$ defined by 
\[
 W_\lambda (x) = \frac{1}{\lambda^{3/2}} W\left(\frac{x}{\lambda}\right). 
\]
%More examples of global nonscattering solutions are given in Donninger-Krieger \cite{nonscaglobal1}. 

\paragraph{Soliton resolution} 
Soliton resolution conjecture predicts that a global solution (or type II blow-up solution) to (CP1) decomposes to a sum of decoupled solitary waves, a radiation term (a linear free wave) and a small error term, as the time tends to infinity (or the blow-up time $T_+$). Radially symmetric assumption guarantees that all possible solitary waves are simply rescaled ground states $\pm W_\lambda$ and that the decoupled condition is equivalent to scale separation. More precisely we have 
\begin{equation} \label{J bubble resolution}
 \vec{u}(t) = \sum_{j=1}^J \zeta_j (W_{\lambda_j(t)},0) + \vec{v}_L(t) + o(1). 
\end{equation}
Here $\zeta_j \in \{+1,-1\}$, $v_L$ is a free wave, the scale functions $\lambda_j(t)$ satisfy 
\begin{align*}
 &\lambda_J(t) \ll \lambda_{J-1}(t) \ll \cdots \ll \lambda_1(t) \ll t, & & t\rightarrow +\infty; \\
 &\lambda_J(t) \ll \lambda_{J-1}(t) \ll \cdots \ll \lambda_1(t) \ll T_+-t, & &t\rightarrow T_+. 
\end{align*}
The radial case of soliton resolution conjecture has been verified in the past decade. The 3-dimensional case was first proved by Duyckaerts-Kenig-Merle \cite{se} via a combination of profile decomposition and channel of energy method. Duyckaerts-Kenig-Merle \cite{oddhigh}, Duyckaerts-Kenig-Martel-Merle \cite{soliton4d} and Collot-Duyckaerts-Kenig-Merle \cite{soliton6d} then proved the odd higher dimensional, 4-dimensional and 6-dimensional cases, respectively, by following roughly the same idea. Recently another proof for 4 or higher dimensional case was given by Jendrej-Lawrie \cite{anothersoliton}. Their proof is a combination of the sequential soliton resolution result and a ``no-return'' theory by the modulation analysis. The soliton resolution conjecture without radially symmetric assumption, however, is still an open problem, although a weaker version of it, i.e. the soliton resolution along a sequence of time, has been prove by Duyckaerts-Jia-Kenig \cite{djknonradial}. 

\paragraph{Multi-soliton solutions} If the soliton resolution of a solution comes with $N$ solitary waves, then we usually call it an $N$-soliton solution or $N$-bubble solution. Multi-soliton solutions have been constructed in many papers, for instance, Martel-Merle \cite{msoliton5d, msoliton5ds} and Yuan \cite{msoliton5dy}. Martel and Merle's work \cite{msoliton5ds} also shows that their 2-soliton solutions scatter in the other time direction, which indicates the inelasticity nature of soliton collision. Please note that all these results are in the non-radial case. 

\subsection{Main idea and result}

This work discusses the long-time behaviour of radial global solutions to (CP1). More precisely we prove a quantitative version of soliton resolution theorem with immediate states for radial global solution to (CP1).  Before we discuss the main idea and result of this work, we first make a brief review of the conception of radiation part and non-radiative solutions. 

\paragraph{Radiation part} Usually the first step to prove the soliton resolution conjecture is to separate the radiation part from a nonlinear solution. If $u$ is a solution to (CP1) defined for all $t\geq 0$, then there exists a linear free wave $v_L$ such that (see Lemma \ref{radiation part})
\[
 \lim_{t\rightarrow +\infty} \int_{|x|>t - A} |\nabla_{t,x} (u-v_L)|^2 {\rm d} x = 0, \qquad \forall A \in \Rm. 
\]
We call $v_L$ the radiation part of $u$ in the positive time direction. The theory of radiation fields(see Section \ref{sec: radiation fields} and \ref{sec: nonlinear profiles}) implies that there exists a function $G_+ \in L^2(\Rm)$, called the radiation profile, such that the following limit holds for any fixed $A \in \Rm$
\begin{align} \label{radiation profile u G}
 \lim_{t\rightarrow +\infty} \int_{t-A}^\infty \left(|r^2  u_r(r,t)+G_+(r-t)|^2 + |r^2 u_t(r,t) - G_+(r-t)|^2 \right) {\rm d} r = 0. 
\end{align}
Next we apply the finite speed of wave propagation and small data theory to deduce that the domain of solution $u$ may extend to the region $\{(x,t): |x|>|t|+R\}$ for a large radius $R$. Here we need to use the conception of exterior solution given in Subsection \ref{sec: exterior solutions}. As a result, given any large time $t_1>R$, there exists a linear free wave $v_{t_1,L}$ asymptotically equivalent to $u$ in the exterior region $\{(x,t): |x|>|t-t_1|\}$, i.e. (see Lemma \ref{vL lemma})
\[
 \lim_{t\rightarrow \pm \infty} \int_{|x|>|t-t_1|} \left|\nabla_{t,x} (u-v_{t_1,L})(x,t)\right|^2 {\rm d} x = 0. 
\]
We call this free wave $v_{t_1,L}$ the radiation part of $u$ in the exterior region $\{(x,t): |x|>|t-t_1|\}$. If the restriction of $v_{t_1,L}$ in this exterior region comes with a small Strichartz norm, then we say that the radiation of $u$ in this exterior region is weak in term of Strichartz norm.  

\paragraph{Non-radiative solutions} A solution $u$ to (CP1) defined in the region $\{(x,t): |x|>|t|+R\}$ is called ($R$-weakly) non-radiative if and only if 
\[
 \lim_{t\rightarrow +\infty} \int_{|x|>|t|+R} |\nabla_{t,x} u(x,t)|^2 {\rm d} x = 0. 
\] 
In other words, non-radiative solutions are exactly those whose radiation part is zero in a suitable exterior region $\{(x,t): |x|>|t|+R\}$. Non-radiative solutions are one of the most important topics in the channel of energy theory, which plays an essential role in the discussion of soliton resolution conjecture. Please refer to Duyckaerts-Kenig-Merle \cite{tkm1} and Kenig-Lawrie-Schlag \cite{channel5d} for more details on the theory of channel of energy. The asymptotic behaviour of non-radiative solutions to energy-critical wave equations has been discussed in Duyckaerts-Kenig-Merle \cite{oddtool} and Collot-Duyckaerts-Kenig-Merle \cite{classNR}. 

\paragraph{Goal and main idea} The author's previous work \cite{dynamics3d} discussed the behaviour of global radial solutions to the focusing, energy-critical wave equation in the 3-dimensional case. Briefly speaking, it was proved that the approximated soliton resolution holds for all time $t>0$ except for several ``relatively short'' time interval. In this work we show that roughly the same result holds for 5-dimensional case as well.  The general idea remains the same as in the 3-dimensional case: if $u$ is a solution with a weak radiation (in term of Strichartz norm) in the main light cone , then we may decompose it into the radiation part, a non-radiative part and a small error term, at least outside a suitable light cone $|x|=|t|+R$. However, the higher dimensional case is much different from the 3-dimensional case. In 3-dimensional case all radial non-radiative solutions are ground states, which are stationary solutions, thus the non-radiative part must be a sum of decoupled ground states. This is exactly the situation of soliton resolution. In the higher-dimensional case, however, there exist non-stationary radial non-radiative solutions, which are much different from the ground states. Thus at least theoretically it is possible to find a solution with weak radiation but its initial data are far away from any soliton resolution situation. 

In order to overcome this difficulty, we need to consider the time evolution of the decomposition mentioned above. First of all, we describe the manner in which a non-stationary non-radiative solutions evolves. In fact, all non-stationary non-radiative solutions in 5-dimensional case evolve in the same way, up to a sign, a dilation and a time translation. A typical example $U$ of these non-radiative solutions is odd in time, approaches the ground state $W(x)$ in one time direction, and approaches the ground state $-W(x)$ in the other direction, if we ignore a time-dependent dilation. In particular, $\|U_t(0)\|_{L^2(\Rm^5)} = +\infty$. Next we show that if the radiation profile $G_+(s)$ of a global solution $u$ to (CP1) is small in a long time interval $[-\ell T, -\ell^{-1} T]$, where $\ell > 1$ is a large constant and $T$ is a large time, then the non-radiative part comes with a much smaller scale than $T$ and evolves in the same manner as the pure non-radiative solutions around the time $T$, possibly besides a few decoupled ground states with significantly larger scales. As a result, there exists a time $t_1$ near $T$, such that we may decompose $\vec{u}(t_1)$ into the radiation part, a rescaled version of $\vec{U}(0)$, a small error term and possibly a few decoupled ground states with larger scales, unless the non-radiative part in the decomposition at time $T$ is very close to a sum of decoupled ground states, which is the case we expect. Finally the fact $\|U_t(0)\|_{L^2(\Rm^5)} = +\infty$ implies that if the decomposition at $t_1$ comes with a term of $\vec{U}(0)$, then the $\dot{H}^1\times L^2$ norm has to blow up near the time $t_1$, which gives a contradiction thus excludes this bad situation. 

\paragraph{Main results} We may divide the main results of this work into two parts: The first part is a series of properties about the radial non-radiative solutions to (CP1), such as their asymptotic behaviour, time evolution, classification and topological structure. Most previous works about non-radiative solutions discuss the asymptotic behaviours, i.e. the properties of non-radiative solutions near the infinity; this work gives the first quantitative property near the origin for a non-radiative solution to (CP1) other than the ground states and the zero solution, both of which can be defined explicitly. To be more precise, we prove that the odd non-radiative solution $U$ mentioned above satisfies the inequality
\[
 \frac{0.05}{|x|^{5/2}} < U_t (x,0) < \frac{1.86}{|x|^{5/2}}, \qquad |x| \ll 1. 
\]
These results about the non-radiative solutions are given in Section 4 and Section 5.

The second part of our main result is the dynamics of global solutions to (CP1). The main theorem given below is actually a quantitative version of soliton resolution theorem. 

\begin{theorem} \label{main thm}
 Given any positive constants $\kappa, \varepsilon\ll 1$ and $E_0 > E(W,0)$, there exists a small constant $\delta = \delta(\kappa,\varepsilon, E_0)> 0$ and two large constants $\ell = \ell(\kappa, \varepsilon, E_0)$, $L = L(\kappa,\varepsilon,E_0)>0$ such that if $u$ is a radial solution to (CP1) satisfying
 \begin{itemize} 
  \item $u$ is defined for all time $t \geq 0$; 
  \item $R$ is a sufficiently large radius such that $\|\vec{u}(0)\|_{\dot{H}^1\times L^2(\{x: |x|>R\})} < \delta$; 
  \item the energy $E$ of $u$ satisfies $E(W,0)\leq E < E_0$;
 \end{itemize}
 then there exists a time sequence $\ell R \leq a_1 < b_1 < a_2 < b_2 < \cdots < a_m < b_m  = +\infty$ such that  
  \begin{itemize}
  \item[(a)] (Soliton resolution in stable periods) For any time interval $[a_k, b_k]$, where $k \in \{1, \cdots, m\}$, there exists a nonnegative integer $J_k$, a linear free wave $v_{k,L}$, a sequence $\zeta_{k,1}, \zeta_{k,2}, \cdots, \zeta_{k,J_k} \in \{\pm 1\}$ and a sequence of scale functions $\lambda_{k,1}(t) > \lambda_{k,2}(t) > \cdots > \lambda_{k,J_k}(t)$ satisfying
 \begin{align*}
  \max\left\{\frac{\lambda_{k,1}(t)}{t}, \frac{\lambda_{k,2}(t)}{\lambda_{k,1}(t)}, \cdots, \frac{\lambda_{k,J_k}(t)}{\lambda_{k,J_k-1}(t)}\right\} \leq \kappa, \qquad t\in [a_k, b_k];
 \end{align*}
 such that
  \begin{align*}
   \left\|\vec{u}(t) - \sum_{j=1}^{J_k} \zeta_{k,j} \left(W_{\lambda_{k,j}(t)},0\right) - \vec{v}_{k,L}(t)\right\|_{\dot{H}^1\times L^2} \leq \varepsilon, \qquad t\in [a_k, b_k].  
  \end{align*}
  Here $t=b_m =+\infty$ is excluded if $k=m$. In addition, the linear free wave $v_{k,L}$ satisfies 
  \begin{align*}
    &\|v_{k,L}\|_{L^{7/3} L^{14/3}([a_k,+\infty)\times \Rm^5)} \leq \varepsilon; & & \|\nabla_{t,x} v_{k,L}(\cdot,t)\|_{L^2(\{x: |x|<t-\ell^{-1}a_k\})} \leq \varepsilon, \quad \forall t\geq \ell^{-1} a_k. 
  \end{align*}  
   We call these time periods ``stable periods''. 
  \item[(b)] (Radiation concentration in collision periods) For each $k \in \{1,2,\cdots, m-1\}$, the bubble numbers $J_{k}$ and $J_{k+1}$ satisfy $J_k >J_{k+1}$.  In addition, the nonlinear radiation profile $G_+$ of $u$ and times $b_k, a_{k+1}$ satisfy ($\sigma_4 = 8\pi^2/3$ is the area of $\mathbb{S}^4$)
  \begin{align*}
   &\|G_+\|_{L^2(-\ell t, -\ell^{-1} t)} \geq \delta, \quad \forall t\in [b_k, a_{k+1}];& &\\
   &\left|\sigma_4 \|G_+\|_{L^2([-a_{k+1},-b_k])}^2 - (J_k-J_{k+1}) E(W,0)\right| \leq \varepsilon^2; & & \frac{a_{k+1}}{b_k} \leq L. 
  \end{align*}
  We call these time periods ``collision periods''. In contrast, for each stable period, we have 
 % \[
%    \left|E - J_k E(W,0) - 4\pi \|G_{+}\|_{L^2([-t,+\infty))}^2\right| \leq \varepsilon^2, \qquad t\in [a_k,b_k].
 % \]
 % Thus
  \begin{align*}
   \sigma_4 \|G_+\|_{L^2((-b_{k}, -a_k])}^2 \leq \varepsilon^2. 
  \end{align*}
 \item[(c)] (Length of preparation period) In addition, we may give an upper bound for the initial time of the first stable period $a_1 \leq LR$.
 \end{itemize}
\end{theorem}

\begin{remark}
 From the proof of the main theorem, we see that the radiation profile of $v_{k, L}$ in the positive time direction can be given by 
 \[
  G_{k,+} (s) = \left\{\begin{array}{ll} G_+ (s), & s > -b_k; \\ 0, & s < - b_k. \end{array}\right.
 \]
 In particular, the last free wave $v_{m, L}$ is exactly the scattering part $v_L$ of $u$. It is not difficult to see that the soliton resolution conjecture of global solutions is a direct consequence of Theorem \ref{main thm}. Please note that our proof does not depend on the type II assumption 
 \begin{equation}
  \limsup_{t\rightarrow +\infty} \|\vec{u}(t)\|_{\dot{H}^1\times L^2} < +\infty, \label{technical type II}
 \end{equation}
 which is typically assumed in the soliton resolution theorem given by previous works(see \cite{oddhigh, anothersoliton}). A direct corollary of our main theorem is that any radial solution to (CP1) defined for all $t\geq 0$ must satisfy \eqref{technical type II}. 
 \end{remark}
 
% \begin{remark} 
% We may give an upper bound of the integer $m$: 
% \[
 % m \leq J_0 + 1 \leq \left\lfloor \frac{E+\varepsilon^2}{E(W,0)} \right\rfloor + 1.
 %\]
 %This follows a combination of the strictly decreasing property of $J_k$ and an upper bound of $J_0$ given in the proof of the main theorem. Please see Section 5.
%\end{remark}

\begin{figure}[h]
 \centering
 \includegraphics[scale=1.1]{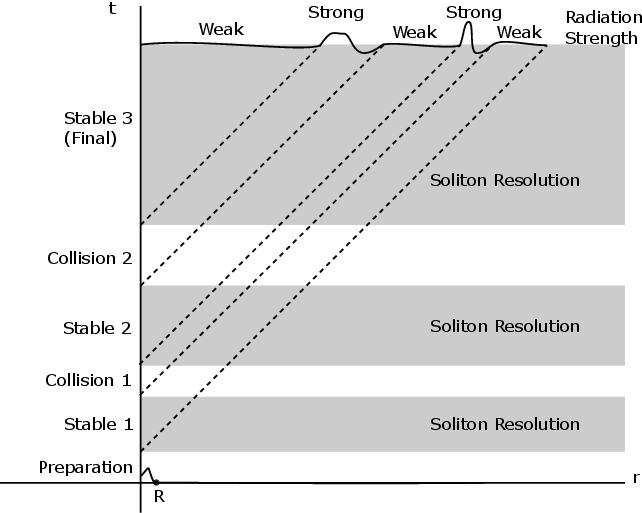}
 \caption{The relationship of stable/collision periods and radiation strength} \label{figure collision}
\end{figure}

\begin{remark}
 According to Theorem \ref{main thm}, we may split the time interval $[0,+\infty)$ into a ``preparation period'' $[0,a_1]$, several ``stable periods'' $[a_k,b_k]$, and several ``collision periods'' $[b_k, a_{k+1}]$ between consecutive stable periods. In each stable period the approximated soliton resolution holds and the bubbles evolve steadily. In each collision period, in contrast, strong interaction happens and at least one bubble is eliminated. It is natural to view the radiation waves travelling in the channel $t-t_2 < |x| < t - t_1$, whose strength can be measured by $\sigma_4 \|G_+\|_{L^2([-t_2,-t_1])}^2$, as the emission of the system during the time period $[t_1,t_2]$. As a result, Theorem \ref{main thm} shows that after the preparation period, almost all radiation comes from the collision periods, whose length is bounded if we apply the logarithm transformation $t' = \ln t$. In addition, the energy of radiation in each collision period is roughly equal to the energy of bubbles eliminated in the collision. This gives a way to understand the long-time dynamics of solutions from their radiation part.  In summary, Long time dynamics in the 5D case is similar to the 3D case: BUBBLE COLLISION GENERATES RADIATION. Please see figure \ref{figure collision} for an illustration of stable/collision periods and their corresponding radiation strength. 
\end{remark}

\begin{remark}
 Given a global solution, Theorem \ref{main thm} may help us predict the upper bound of time at which the solution first reaches an approximated soliton resolution state, simply from the energy $E$ and scale $R$ of the initial data. Please note that this first soliton resolution state is not necessarily the final soliton resolution state described in the soliton resolution conjecture. In fact, it is impossible to predict the time when the solution reaches its final state only from the assumption on $u$ in Theorem \ref{main thm}, as shown by an example constructed in Duyckaerts-Merle \cite{threshold}. According to Duyskaerts-Merle's work, there exists a solution $v$ to (CP1), such that
 \begin{itemize}
  \item $\vec{v}(t)$ converges to $(W,0)$ in $\dot{H}^1\times L^2$ as $t\rightarrow -\infty$;
  \item $v$ scatters in the positive time direction. 
 \end{itemize}
 By choosing the initial data $(u_0,u_1)= \vec{v}(t_1)$ for a large negative number $t_1$, we may make the initial data sufficiently closed to $(W,0)$ in $\dot{H}^1\times L^2$, and make the time when the solution reaches the final(scattering) state arbitrarily large. This also provides an example of global solutions with at least two stable periods and one collision period.  
\end{remark}

%\begin{remark}
% The proof of Theorem \ref{main thm} utilizes neither the profile decomposition nor a sequential version of the soliton resolution. It is much different from the previously known proof of the soliton resolution conjecture. The major tool of the proof is the radiation theory. The radiation theory discusses not only the energy in the exterior region, which is the main topic of the channel of energy method, but also the radiation profile defined above. 
%\end{remark}

\paragraph{Structure of this work} This work is organized as follows: We introduce notations and give some preliminary results in Section 2. Then in Section 3 we present two key lemmata, which help us compare a solution to (CP1) and an approximated solution to (CP1) with the same radiation part. These two lemmata are major tools of this work. Sections 4 and 5 are devoted to the discussion of non-radiative solutions to (CP1). In Section 4 we discuss the asymptotic behaviour, maximal domain, topological structure and time evolution of non-radiative solutions to (CP1). In section 5 we investigate the behaviour of odd non-radiative solutions near the origin, which plays an essential role in the subsequent sections. In Section 6 we show that the minimum value of the energy norm $\|\vec{u}(t)\|_{\dot{H}^1\times L^2}$ in a long time interval is dominated by a linear function of $E$, and then extract the radiation part of a global solution to (CP1) as an application. Finally we combine all the ingredients given above and prove the main theorem in the last two sections. 

\section{Preliminary results} 

\subsection{Exterior solutions} \label{sec: exterior solutions}

For convenience of our discussion, it is helpful to introduce solutions to (CP1) defined only in an exterior region. Before we discuss the basic conception of exterior solutions, we introduce a few notations. Given $R\geq 0$, we call the following region 
\[
 \Omega_R = \left\{(x,t)\in \Rm^5 \times \Rm: |x|>|t|+R\right\}
\]
an exterior region and use the notation $\chi_R$ for the characteristic function of $\Omega_R$. In some situations we also allow $R<0$ and define $\chi_R$ in the same manner. Given a time interval $J$, we use the following notation for the classic Strichartz norm 
\[
 \|u\|_{Y(J)} = \|u\|_{L^{7/3} L^{14/3}(J \times \Rm^3)} = \left(\int_J \left(\int_{\Rm^3} |u(x,t)|^{14/3} {\rm d} x\right)^{1/2} {\rm d} t\right)^{3/7}. 
\]

\paragraph{Exterior solutions} Let $u$ be a function defined in the exterior region 
\[
 \Omega =  \left\{(x,t)\in \Rm^5\times (-T_1,T_2): |x|>|t|+R\right\}.
\]
Here $T_1,T_2$ are either positive real numbers or $\infty$. We call $u$ an exterior solution to (CP1) in the region $\Omega$ with initial data $(u_0,u_1)\in \dot{H}^1 \times L^2(\Rm^5)$, if and only if $\|\chi_R u\|_{Y(J)} < +\infty$ for any bounded closed time interval $J \subset (-T_1,T_2)$ and the following identity holds:
 \[
  u = \mathbf{S}_L (u_0,u_1) + \int_0^t \frac{\sin (t-t')\sqrt{-\Delta}}{\sqrt{-\Delta}} [\chi_R(\cdot,t') F(u(\cdot,t'))] {\rm d} t', \qquad |x|>R+|t|, \; t\in (-T_1,T_2).
 \]
 Here $\mathbf{S}_L$ represents the linear wave propagation operator. We multiply $u$ and $F(u)$ by the characteristic function $\chi_R$ to emphasize that $u$ and $F(u)$ are only defined in the exterior region $\Omega$. More precisely we understand the product in the following way
 \begin{align*}
  &\chi_R u = \left\{\begin{array}{ll} u(x,t), & (x,t)\in \Omega_R; \\ 0, & (x,t)\notin \Omega_R. \end{array} \right.
  &\chi_R F(u) = \left\{\begin{array}{ll} F(u(x,t)), & (x,t) \in \Omega_R; \\ 0, & (x,t)\notin \Omega_R. \end{array} \right.
 \end{align*}
 Although we define the initial data for all $x \in \Rm^5$ in the definition above, finite speed of propagation implies that the values of initial data in the ball $\{x: |x|<R\}$ are irrelevant. For convenience we let $\mathcal{H}(R)$ be the Hilbert space consisting of restrictions of radial $\dot{H}^1\times L^2$ functions on the exterior region $\{x: |x|>R\}$. The norm of $\mathcal{H}(R)$ is defined by 
 \[
  \|(u_0,u_1)\|_{\mathcal{H}(R)}^2 =  \int_{|x|>R} \left(|\nabla u_0(x)|^2 + |u_1(x)|^2 \right) {\rm d} x. 
 \]
 When we talk about a radial exterior solution defined as above, we may specify its initial data by $(u_0,u_1)\in \mathcal{H}(R)$. 
 Similarly we may define an exterior solution $u$ to the wave equation 
 \[
  \left\{\begin{array}{l} \partial_t^2 u - \Delta u = F(x,t), \qquad (x,t)\in \Omega;\\ (u,u_t)|_{t=0} = (u_0,u_1). \end{array}\right.
 \]
 in the same manner, if $F$ is defined in the exterior region $\Omega$ and satisfies $\|\chi_R F\|_{L^1 L^2 (J \times \Rm^5)}< +\infty$ for any bounded closed interval $J \subset (-T_-,T_+)$. 
 
% \paragraph{Restriction and global extension} If $u$ is an exterior solution to (CP1) defined in the exterior region $\Omega_R$ for $t \in I$. Then finite speed of propagation shows that its restriction in the exterior region $\Omega_r$ with $r>R$ is also an exterior solution to (CP1) with the same initial data. On the other hand, the function defined for $(x,t)\in \Rm^3 \times I$ by the formula
 %\[ 
 %  \mathbf{S}_L (u_0,u_1) + \int_0^t \frac{\sin (t-t')\sqrt{-\Delta}}{\sqrt{-\Delta}} [\chi_R(\cdot,t') F(t', \cdot, u)] dt'
 %\]
%coincides with $u$ in the exterior region $\Omega_R$ and solves the modified wave equation below defined in the whole space $\Rm^3$
% \[
%  \partial_t^2 u - \Delta u = \chi_R F(x,t,u), \qquad (x,t) \in \Rm^3 \times I.
% \]
%We call it the global extension of $u$ and still use the notation $u$ to represent it. Please note that the global extension does depend on the choice of initial data in the interior region $\{x: |x|<R\}$. The Strichartz estimates then show that the global extension satisfies ($J \subset I$ is a bounded closed interval)
% \begin{align*}
%  &\|u\|_{Y(J)} < +\infty;& &(u, u_t) \in C(I; \dot{H}^1\times L^2(\Rm^3)).&
% \end{align*}

\paragraph{Local theory} The local well-posedness of initial value problem in the exterior region immediately follows from a combination of the Strichartz estimates (see \cite{strichartz} for instance) and a fixed-point argument. The argument is similar to those in the whole space $\Rm^5$ and somewhat standard nowadays. More details of these types of argument can be found in \cite{loc1, ls}.  

\paragraph{Perturbation theory} The continuous dependence of exterior solution on the initial data/error function immediately follows from the following lemma
\begin{lemma} \label{perturbation lemma}
 Let $M > 0$ be a constant. Then there exists two positive constants $\delta = \delta(M)$ and $C= C(M)$, such that if $v$ is a radial exterior solution to 
 \[
  \left\{\begin{array}{l} \partial_t^2 v - \Delta v = F(v) + e(x,t), \qquad (x,t)\in \Omega_R; \\
  (v,v_t)|_{t=0} = (v_0,v_1) \in \mathcal{H}(R) \end{array}\right.
 \]
 satisfying 
 \begin{align*}
  &\|\chi_R v\|_{Y(\Rm)} < M;& & \|\chi_R e(x,t)\|_{L^1 L^2(\Rm \times \Rm^5)} < \delta;
 \end{align*}
 and $(u_0,u_1)$ are a pair of radial initial data satisfying $\|(u_0,u_1)-(v_0,v_1)\|_{\mathcal{H}(R)} < \delta$, then the corresponding exterior solution $u$ to (CP1) in the exterior region $\Omega_R$ with initial data $(u_0,u_1)$ can be defined for all $t\in \Rm$ with 
 \begin{align*}
   \|\chi_R (u-v)\|_{Y(\Rm)} + \sup_{t\in \Rm} \|\vec{u}(t) -\vec{v}(t)\|_{\mathcal{H}(R+|t|)} \leq C\left(\|\chi_R e\|_{L^1 L^2(\Rm \times \Rm^5)} + \|(u_0,u_1)-(v_0,v_1)\|_{\mathcal{H}(R)}\right). 
 \end{align*}
 Here $R\geq 0$ is an arbitrary constant. 
\end{lemma}

The proof of Lemma \ref{perturbation lemma} is similar to the situation when the solution is defined in the whole space-time. Please see \cite{kenig, shen2}, for instance. More details about the exterior solutions and a similar perturbation theory can be found in Duyskaerts-Kenig-Merle \cite{oddtool}. 
 %\begin{itemize}
%  \item Well-posedness: Given any initial data $(u_0,u_1) \in \dot{H}^1 \times L^2(\Rm^3)$ and $R\geq 0$, there exists a unique exterior solution to (CP1) in the region $\Omega_R$ with a maximal lifespan $(-T_-,T_+)$.
%  \item Small data scattering: There exists a constant $\delta = \delta(\gamma) >0$, so that if the initial data satisfy $\|\chi_R \mathbf{S}_L (u_0,u_1)\|_{Y(\Rm)} < \delta$, then the corresponding exterior solution $u$ to (CP1) is defined for all $t\in\Rm$, so that $\|\chi_R u\|_{Y(\Rm)} \leq 2 \|\chi_R \mathbf{S}_L (u_0,u_1)\|_{Y(\Rm)}$.
 % \item Finite time blow-up criterion: If the exterior solution $u$ blow up in finite time, i.e. $T_+ < +\infty$, then $\|\chi_R u\|_{Y([0,T_+))} = +\infty$. The situation in negative time direction is similar. 
 % \item Scattering criterion: The exterior solution $u$ in $\Omega_R$ scatters in the positive time direction, i.e. $T_+ =  +\infty$ and there exists a linear free wave $u_L^+$ so that 
 % \[
 %  \lim_{t\rightarrow +\infty} \|\vec{u} - \vec{u}_L^+\|_{\dot{H}^1 \times L^2(\{x: |x|>|t|+R\})} = 0, 
 % \]
 % if and only if $\|\chi_R u\|_{Y([0,T_+))} < +\infty$. 
 %\end{itemize}

\subsection{Radiation fields of free waves} \label{sec: radiation fields} 

One of main tools of this work is the radiation field, which has a history of more than 50 years. Please see, Friedlander \cite{radiation1, radiation2} for instance. Generally speaking, radiation fields discuss the asymptotic behaviour of linear free waves. The following version of statement comes from Duyckaerts-Kenig-Merle \cite{dkm3}.

\begin{theorem}[Radiation fields] \label{radiation}
Assume that $d\geq 3$ and let $u$ be a solution to the free wave equation $\partial_t^2 u - \Delta u = 0$ with initial data $(u_0,u_1) \in \dot{H}^1 \times L^2(\Rm^d)$. Then ($u_r$ is the derivative in the radial direction)
\[
 \lim_{t\rightarrow \pm \infty} \int_{\Rm^d} \left(|\nabla u(x,t)|^2 - |u_r(x,t)|^2 + \frac{|u(x,t)|^2}{|x|^2}\right) {\rm d}x = 0
\]
 and there exist two functions $G_\pm \in L^2(\Rm \times \mathbb{S}^{d-1})$ such that
\begin{align*}
 \lim_{t\rightarrow \pm\infty} \int_0^\infty \int_{\mathbb{S}^{d-1}} \left|r^{\frac{d-1}{2}} \partial_t u(r\theta, t) - G_\pm (r\mp t, \theta)\right|^2 {\rm d}\theta {\rm d}r &= 0;\\
 \lim_{t\rightarrow \pm\infty} \int_0^\infty \int_{\mathbb{S}^{d-1}} \left|r^{\frac{d-1}{2}} \partial_r u(r\theta, t) \pm G_\pm (r\mp t, \theta)\right|^2 {\rm d}\theta {\rm d} r & = 0.
\end{align*}
In addition, the maps $(u_0,u_1) \rightarrow \sqrt{2} G_\pm$ are bijective isometries from $\dot{H}^1 \times L^2(\Rm^d)$ to $L^2 (\Rm \times \mathbb{S}^{d-1})$. 
\end{theorem}

\noindent In this work we call $G_\pm$ the radiation profiles of the linear free wave $u$, or equivalently, of the corresponding initial data $(u_0,u_1)$. Clearly the map/symmetry between radiation profiles $G_\pm$ is an isometry from $L^2(\Rm \times \mathbb{S}^{d-1})$ to itself. The symmetry of $G_\pm$ in 5-dimensional case can be given by (please see \cite{newradiation, shenradiation} for other dimensions, for example)
\begin{equation} \label{symmetry of G pm}
 G_+(s,\theta) = G_-(-s, -\theta). 
\end{equation}
Thus we may uniquely determine a linear free wave by its whole radiation profile in either time direction, or both two radiation profiles $G_\pm(s,\theta) \in L^2(\Rm^+ \times \mathbb S^4)$ but for positive $s$ only. It is not difficult to see that the free wave is radial if and only if its radiation profiles are independent of the angle $\theta$. The formula of a free wave in term of its radiation profile can also be given explicitly, see \cite{shenradiation}, for example. In this work we focus on the 5D radial case: 
\begin{equation} \label{free wave by radiation profile}
 u(r,t) = \frac{1}{r^3} \int_{t-r}^{t+r} (s-t) G_-(s) {\rm d}s.  
\end{equation}
A basic calculation gives the initial data in term of the radiation profile 
\begin{align} \label{initial data by radiation profile} 
 &u_0(r) =  \frac{1}{r^3} \int_{-r}^{r} s G_-(s) {\rm d}s; & & u_1(r) = \frac{1}{r^2} \left[G_-(r)+G_-(-r)\right] - \frac{1}{r^3}\int_{-r}^r G_-(s) {\rm d}s.
\end{align}
The following energy formula in the exterior region is useful in further argument. 
\begin{lemma} \label{tail G}
 Let $(u_0,u_1) \in \dot{H}^1\times L^2(\Rm^5)$ are radial initial data, whose radiation profile in the negative time direction is $G(s)$. Then we have
 \[
  \frac{1}{\sigma_4}\|(u_0,u_1)\|_{\mathcal{H}(R)}^2 = 2 \|G\|_{L^2(\{s: |s|>R\})}^2 + \frac{1}{R} \left(\int_{-R}^R G(s) {\rm d} s\right)^2 + \frac{3}{R^3} \left(\int_{-R}^R s G(s) {\rm d} s \right)^2 . 
 \]
 Here $\sigma_4$ is the area of unit sphere $\mathbb{S}^4$. By symmetry the same identity holds if we use the radiation profile in the positive directions. 
\end{lemma} 
\begin{proof}
By \eqref{initial data by radiation profile}, we may integrate by parts and obtain 
\begin{align*}
 \int_R^\infty (r^2 \partial_r u_0(r)) {\rm d} r & = \int_R^\infty [G(r)-G(-r)]^2 {\rm d} r + \frac{3}{R^3} \left(\int_{-R}^R s G(s) {\rm d} s \right)^2;\\
 \int_R^\infty (r^2 u_1(r))^2 {\rm d} r & = \int_R^\infty [G(r)+G(-r)]^2 {\rm d} r + \frac{1}{R} \left(\int_{-R}^R G(s) {\rm d} s\right)^2.
\end{align*}
A combination of these two identities immediately prove the lemma.
\end{proof}

Given a radial linear free, or its initial data, or equivalently its radiation profile $G$ (in the negative direction), we define the radiation residues $\vec{\tau}(R)= (\tau_1(R),\tau_2(R))$ at the radius $R$ by 
\begin{align} \label{characteristic profile}
 &\tau_1(R) = \frac{-1}{R^{1/2}} \int_{-R}^R G(s) {\rm d} s;& & \tau_2(R) = \frac{\sqrt{3}}{R^{3/2}} \int_{-R}^R s G(s) {\rm d} s. 
\end{align} 
Given $R>0$, it follows from the explicit formula \ref{free wave by radiation profile} and Lemma \ref{tail G} that the map defined by $(u_0,u_1)\rightarrow (\sqrt{2} G_-(s), \vec{\tau}(R))$ is an isometry from $\mathcal{H}(R)$ to $L^2(\{s: |s|>R\}) \times \Rm^2$. Radiation residues will be used frequently in this work. 

\begin{remark} \label{double tail G}
 A direct consequence of Lemma \ref{tail G} is 
 \[
  \sigma_4^{-1} \|(u_0,u_1)\|_{\mathcal{H}(R)}^2 \geq 2\|G_-\|_{L^2(\{s: |s|>R\})}^2, \qquad \forall R > 0. 
 \]
 It immediately follows that if $0\leq R_1 < R_2$, then 
 \[
  \sigma_4^{-1} \|(u_0,u_1)\|_{\mathcal{H}(R_1)}^2 \leq \sigma_4^{-1} \|(u_0,u_1)\|_{\mathcal{H}(R_2)}^2 + 2 \|G_-\|_{L^2(\{s: R_1 < |s|<R_2\})}^2 + \left|\vec{\tau}(R_1)\right|^2. 
 \]
\end{remark}

\subsection{Nonlinear radiation profiles} \label{sec: nonlinear profiles}

\begin{lemma} [Radiation fields of inhomogeneous equation] \label{scatter profile of nonlinear solution}
 Assume that $R\geq 0$. Let $u$ be a radial exterior solution to the wave equation
 \[
  \left\{\begin{array}{ll} \partial_t^2 u - \Delta u = F(t,x); & (x,t)\in \Omega_R; \\
  (u,u_t)|_{t=0} = (u_0,u_1) \in \dot{H}^1\times L^2. & \end{array} \right.
 \]
 If $F$ is a radial function satisfying $\|\chi_R F\|_{L^1 L^2(\Rm\times \Rm^5)}< +\infty$, then there exist unique radiation profiles $G_\pm \in L^2([R,+\infty))$ such that 
 \begin{align}
  \lim_{t\rightarrow +\infty} \int_{R+t}^\infty \left(\left|G_+(r-t) - r^2 u_t (r, t)\right|^2 + \left|G_+(r-t) + r^2 u_r (r, t)\right|^2\right) {\rm d}r & = 0; \label{positive ra}\\
  \lim_{t\rightarrow -\infty} \int_{R-t}^\infty \left(\left|G_-(r+t) - r^2 u_t(r,t)\right|^2 +  \left|G_-(r+t) - r^2 u_r(r,t)\right|^2\right) {\rm d} r & = 0. \label{negative ra}
 \end{align}
 In addition, the following estimates hold for $G_\pm$ given above and the radiation profiles $G_{0,\pm}$ of the initial data $(u_0,u_1)$:
 \begin{align*}
 \sqrt{2\sigma_4} \|G_- - G_{0,-}\|_{L^2([R,+\infty))} & \leq \|\chi_R F\|_{L^1 L^2((-\infty,0]\times \Rm^5)}; \\
  \sqrt{2\sigma_4}  \|G_+ - G_{0,+}\|_{L^2([R,+\infty))} & \leq \|\chi_R F\|_{L^1 L^2([0,+\infty)\times \Rm^5)}.
 \end{align*}
\end{lemma}
The proof of this lemma is the same as in the 3-dimensional case. Please refer to \cite{dynamics3d}. In this work we also need to calculate the radiation profiles explicitly. More precisely we have 

\begin{lemma} \label{explicit formula for G plus by F}
The radiation profile of radial initial data $(0,u_1) \in \dot{H}^1\times L^2(\Rm^5)$ can be given by the formula 
\[
 G_+ (s) =  \frac{1}{2}  s^2 u_1(s) - \frac{1}{2} \int_{s}^\infty r u_1(r) {\rm d} r, \qquad s>0.
\]
 In addition, the unique radiation profile $G_+$ defined in Lemma \ref{scatter profile of nonlinear solution} can be given by the formula
 \[
  G_+ (s) - G_{0,+} (s) = \frac{1}{2}\int_0^\infty (s+t)^2 F(t,t+s) {\rm d} t - \frac{1}{2}\int_0^\infty \int_{t+s}^\infty r F(t,r) {\rm d} r {\rm d} t, \qquad s>R.   
 \]
\end{lemma}
\begin{proof}
 It suffices to prove the first formula, since the second one immediately follows from the Duhamel's formula 
 \[
  u(T) = \mathbf{S}_L(T) (u_0,u_1) + \int_0^t \mathbf{S}_L(T-t) (0,F(t)) {\rm d} t. 
 \]
 By Hardy's inequality, the right hand side is a bounded operator from the space of radial $L^2 (\Rm^5)$ functions to $L^2(\Rm^+)$. Thus it suffices to consider smooth and compactly-supported functions $u_1$. By convention we use the notations $u(r,t)$ and $u_1(r)$ for the values of $u(x,t)$ and $u_1(x)$ with $|x|=r$. By using the explicit formula of solution s to the linear wave equation, we have  
 \begin{equation} \label{positive SLu1}
  u(|x|,t) = \frac{1}{4|x|^3} \int_{|x|-t}^{|x|+t} r \left[|x|^2 + r^2 - t^2 \right] u_1(r) {\rm d} r, \qquad |x|>t. 
 \end{equation}
 A direct calculation shows that 
 \begin{align*}
  |x|^2 u_t (|x|,t) =  \frac{1}{2} \left[(|x|+t)^2 u_1(|x|+t) + (|x|-t)^2 u_1 (|x|-t)\right] - \frac{1}{2|x|} \int_{|x|-t}^{|x|+t} t r u_1(r) {\rm d} r . 
 \end{align*}
 Thus 
 \begin{align*} 
 (t+s)^2 u_t (t+s,t) = \frac{1}{2} \left[(2t+s)^2 u_1(2t+s) + s^2 u_1 (s)\right] - \frac{1}{2(t+s)} \int_{s}^{2t+s} t r u_1(r) {\rm d} r . 
 \end{align*}
 Taking a limit, we have 
 \begin{align*}
  G_+(s)  = \lim_{t\rightarrow +\infty} (t+s)^2 u_t (t+s,t) = \frac{1}{2} s^2 u_1 (s) - \frac{1}{2} \int_{s}^{\infty}  r u_1(r) {\rm d} r,
 \end{align*}
 which finishes the proof. 
\end{proof}

\begin{corollary} \label{continuity of non-radiative} 
 Let $u$ be a radial non-radiative exterior solution in $\Omega_R$. Then $u$ must be $C^1$ in $\Omega_R$ whose initial data come with continuous radiation profiles $G_\pm (s)$ for $|s|>R$. 
\end{corollary}
\begin{proof}
 By the non-radiative assumption, we have 
 \begin{equation} \label{radiation profile initial non-radiative}
  G_+ (s) =  -\frac{1}{2}\int_0^\infty (s+t)^2 F(u(t,t+s)) {\rm d} t + \frac{1}{2}\int_0^\infty \int_{t+s}^\infty r F(u(t,r)) {\rm d} r {\rm d} t, \qquad s>R. 
 \end{equation}
 These integrands decay sufficiently fast and uniformly (see Lemma \ref{finite eq scattering}) near the infinity by the point-wise estimate $|u(r,t)| \lesssim r^{-3/2}$. The continuity of $G_+$ for $s>R$ then follows from this uniform decay and the continuity of $u$. By time symmetry we may deduce the continuity of $G_-(s)$ for $s>R$ in the same manner. Thus $G_\pm(s)$ are continuous for $|s|>R$, which implies that the initial data $(u_0,u_1)$ are $C^1$ for $|x|>R$. The $C^1$ of solution $u$ then follows. 
\end{proof}

\subsection{Asymptotically equivalent solutions}

Assume that $u, v \in \mathcal{C}(\Rm; \dot{H}^1\times L^2)$. We say that $u$ and $v$ are $R$-weakly asymptotically equivalent if
\[
 \lim_{t\rightarrow \pm \infty} \int_{|x|>R+|t|} |\nabla_{t,x} (u-v)|^2 {\rm d} x = 0.
\]
Here $R\geq 0$ is a constant. In particular, if $R=0$, then we say $u$ and $v$ are asymptotically equivalent to each other. Since the integral above only involves the values of $u, v$ in the exterior region, the definition above also applies to exterior solutions. A solution $u$ is called ($R$-weakly) non-radiative solution if and only if it is asymptotically equivalent to zero. Non-radiative solutions, which play an essential role in the channel of energy method, have been extensively studied in recent years. Let us consider two examples. We start by considering an $R$-weakly non-radiative radial free wave $u$. It is equivalent to saying that the radiation profiles $G_\pm (s) = 0$ for $s>R$, or $G_-(s) = 0$ for $|s|>R$. An application of the explicit formula of linear free wave in term of radiation profile immediately gives 
\[
 u(r,t) = \frac{1}{r^3} \int_{-R}^{R} s G_-(s) {\rm d}s - \frac{t}{r^3} \int_{-R}^R G_-(s) {\rm d} s, \qquad r>|t|+R. 
\]
This implies that all $R$-weakly non-radiative free waves form a two-dimensional linear space spanned by $1/r^3$ and $t/r^3$. The coefficients for these two basis are exactly constant multiples of radiation resides at the radius $R$. Next we consider all non-radiative solutions to (CP1). One specific example of non-radiative solution is exactly the ground state mentioned in the introduction section
\[
 W(x) = \left(1 + \frac{|x|^2}{15}\right)^{-3/2}. 
\] 
In fact this is the unique stationary solution up to a rescaling/sign symmetry. In other words, all nonzero radial stationary $\dot{H}^1$ solutions can be given by 
\begin{align*}
 &\pm W_\lambda (x);& &W_\lambda (x) = \frac{1}{\lambda^{3/2}} W\left(\frac{x}{\lambda}\right), \qquad \lambda > 0. 
\end{align*}
Before we conclude this section, we give a sufficient and necessary condition for an exterior solution to be asymptotically equivalent to some linear free wave.
\begin{lemma} \label{finite eq scattering}
 Let $u$ be a radial exterior solution to (CP1) defined in $\Omega_R$, then $u$ is $R$-weakly asymptotically equivalent to some finite-energy linear free wave $w_L$ if and only if $\|\chi_R u\|_{Y(\Rm)} < +\infty$. 
\end{lemma}
\begin{proof}
 If $\|\chi_R u\|_{Y(\Rm)} < +\infty$, then we have $\|\chi_R F(u)\|_{L^1 L^2(\Rm \times \Rm^5)} < +\infty$. The existence of asymptotically equivalent free waves has been given in the proof of Lemma \ref{scatter profile of nonlinear solution}. Conversely if $u$ is $R$-weakly asymptotically equivalent to a free wave $w_L$, then we have 
 \[
  \lim_{t\rightarrow +\infty} \left\|\vec{w}_L(t) - \vec{u}(t)\right\|_{\mathcal{H}(R+t)} = 0. 
 \]
 A combination of this limit with the finite speed of propagation and the fact 
 \[
  \lim_{t \rightarrow +\infty} \|w_L\|_{Y([t,+\infty))} = 0
 \]
 yields 
 \[
  \lim_{t\rightarrow +\infty} \|\chi_{R+t} \mathbf{S}_L(\vec{u}(t))\|_{Y(\Rm^+)} = 0. 
 \]
 The small data theory and the uniqueness of exterior solution then guarantees that 
 \[
  \|\chi_R u \|_{Y([t,+\infty))} < +\infty, \qquad \forall t \gg 1. 
 \]
 This yields the estimate $\|\chi_R u\|_{Y(\Rm^+)}<+\infty$, since the definition of exterior solutions guarantees that $\|\chi_R u\|_{Y([0,t])}<+\infty$ for all $t>0$. The negative time direction can be dealt with in a similar way. 
\end{proof}

\section{Comparison of asymptotically equivalent solutions}

In this section we introduce a few key observations, which help us compare two asymptotically equivalent solutions. We shall utilize these lemmata to investigate the properties of solutions with weak radiation, especially the non-radiative solutions. These lemmata can be proved in the same manner as in the 3-dimensional case. 

\begin{remark}
 The results in this section apply to wave equations with a more general nonlinear term. We only need to assume that the nonlinear term $F(u)$ satisfies the energy-critical assumption 
 \begin{align*}
  &F(0) = 0;& & |F(u)-F(v)| \lesssim_1 |u-v|\left(|u|^{4/3}+|v|^{4/3}\right). 
 \end{align*}
\end{remark}

\begin{lemma} \label{lemma inequality}
 Let $u$, $S$ be exterior solutions of (CP1) and $(\partial_t^2 - \Delta) S = F(S) + e(x,t)$ in $\Omega_{R_1}$, respectively, with 
 \begin{align*}
  \|\chi_{R_1} u\|_{Y(\Rm)}, \|\chi_{R_1} S\|_{Y(\Rm)}, \|\chi_{R_1} e(x,t)\|_{L^1 L^2(\Rm\times \Rm^5)} < +\infty.  
 \end{align*}
 Let $w = u - S$ and $G$, $\vec{\tau}$ be the radiation profile and radiation residue of the initial data $\vec{w}(0)$. There exists an absolute constant $C_1\geq 1$ such that the following inequality holds for any $R_2 > R_1\geq 0$:
 \begin{align*}
  \|\chi_{R_1} w\|_{Y(\Rm)} &\leq C_1 \left(|\vec{\tau}(R_1)| + \|G\|_{L^2(\{s: R_1<|s|<R_2\})} + \|\vec{w}(0)\|_{\mathcal{H}(R_2)}\right)\\
  & \qquad + C_1\left(\left\|\chi_{R_2}\left(F(u)-F(S) - e(x,t))\right)\right\|_{L^1 L^2(\Rm \times \Rm^5)} + \|\chi_{R_1,R_2}e(x,t)\|_{L^1 L^2(\Rm \times \Rm^5)}\right)\\
  & \qquad \qquad + C_1\left(\|\chi_{R_1,R_2} w\|_{Y(\Rm)}^{4/3} + \|\chi_{R_1,R_2}S\|_{Y(\Rm)}^{4/3} \right)\|\chi_{R_1,R_2} w\|_{Y(\Rm)}.
 \end{align*}
  Here $\chi_{R_1,R_2}$ is the characteristic function of the region 
 \[
  \Omega_{R_1,R_2} = \{(x,t) \in \Rm^5 \times \Rm: |t|+R_1<|x|<|t|+R_2\}. 
 \]
 In addition, the inequality $|F(x+y)-F(y)|\leq C_1 |x|(|x|^{4/3}+|y|^{4/3})$ holds for all numbers $x, y$. 
\end{lemma}
\begin{proof}
 It is sufficient to prove the first inequality, because the second inequality clearly holds for a sufficiently large constant $C_1$. First of all, we may apply Strichartz estimates, as well as Remark \ref{double tail G}, and obtain 
 \begin{align*}
  \|\chi_{R_1} w_L\|_{Y(\Rm)} & \lesssim_1 \|(w(0), w_t(0))\|_{\mathcal{H}(R_1)} \\
  & \lesssim_1 |\vec{\tau}(R_1)| + \|G\|_{L^2(\{s: R_1<|s|<R_2\})} + \|\vec{w}(0)\|_{\mathcal{H}(R_2)}. 
 \end{align*}
 Here $w_L$ is the linear free wave with initial data $\vec{w}(0)$. Since $w$ satisfies the equation $(\partial_t^2 - \Delta) w = F(u) - F(S) - e(x,t)$, we have 
 \begin{align*}
  \|\chi_{R_1} w\|_{Y(\Rm)} & \lesssim_1 \|\chi_{R_1} w_L\|_{Y(\Rm)} + \left\|\chi_{R_1}\left(F(u) - F(S) - e(x,t)\right)\right\|_{L^1 L^2} \\
  & \lesssim_1 \|\chi_{R_1} w_L\|_{Y(\Rm)} + \left\|\chi_{R_2}\left(F(u) - F(S) - e(x,t)\right)\right\|_{L^1 L^2} \\
  & \qquad + \left\|\chi_{R_1,R_2}\left(F(u) - F(S)\right)\right\|_{L^1 L^2} + \left\|\chi_{R_1,R_2}e(x,t)\right\|_{L^1 L^2}. 
 \end{align*}
 Finally H\"{o}lder inequality gives 
 \begin{align*}
  \left\|\chi_{R_1,R_2}\left(F(u) - F(S)\right)\right\|_{L^1 L^2} \lesssim_1 \left(\|\chi_{R_1,R_2} w\|_{Y(\Rm)}^{4/3} + \|\chi_{R_1,R_2}S\|_{Y(\Rm)}^{4/3}\right)\|\chi_{R_1,R_2} w\|_{Y(\Rm)}. 
 \end{align*}
 A combination of these inequalities finishes the proof. 
\end{proof} 

\begin{lemma} \label{lemma connection} 
 There exists absolute positive constants $\varepsilon_1$, $\beta$, $\eta$ such that if $0\leq R_1<R_2$ and
 \begin{itemize} 
  \item $u$ is an exterior solution to (CP1) and $S$ is an exterior solution to the equation 
 \[
  (\partial_t^2 -\Delta) S= F(S) + e(x,t),
 \] 
 both in the region $\Omega_{R_1}$, with $\|\chi_{R_1} u\|_{Y(\Rm)}, \|\chi_{R_1} S\|_{Y(\Rm)}, \|\chi_{R_1} e(x,t)\|_{L^1 L^2} < +\infty$; 
 \item both $u$, $S$ are asymptotically equivalent to each other in $\Omega_{R_1}$; 
 \item $u$, $S$, $w=u-S$ and the radiation residue $\vec{\tau}(r)$ of $\vec{w}(0)$ satisfy the following inequalities
 \begin{align*}
  \varepsilon \doteq \|(w(\cdot,0),w_t(\cdot,0))\|_{\mathcal{H}(R_2)}  + \|\chi_{R_1,R_2} e(x,t)\|_{L^1 L^2(\Rm \times \Rm^5)}  & \\
 + \|\chi_{R_2}\left(F(u) - F(S)- e(x,t)\right)\|_{L^1 L^2(\Rm \times \Rm^5)} & \leq \varepsilon_1; \\
   \|\chi_{R_1,R_2} S\|_{Y(\Rm)} & \leq \eta; \\
  \sup_{R_1\leq r\leq R_2} \left|\vec{\tau}(r)\right| & \leq \beta;
 \end{align*}
 \end{itemize}
 then we have 
 \begin{align*}
  \|\chi_{R_1} w\|_{Y(\Rm)} + \sup_{t\in \Rm} \|\vec{w}(t)\|_{\mathcal{H}(R_1+|t|)} \lesssim_1 \left|\vec{\tau}(R_1)\right| +\varepsilon.
 \end{align*}
\end{lemma}
\begin{proof}
 Let $w_L$ and $G$ be the linear free wave and radiation profile with initial data $(w(\cdot,0), w_t(\cdot,0))$. By Lemma \ref{lemma inequality}, we obtain for any $R\in [R_1,R_2)$ that
 \begin{align*}
  \|\chi_{R} w\|_{Y(\Rm)} \leq C_1 \left(\left|\vec{\tau}(R)\right| + \|G\|_{L^2(\{s: R<|s|<R_2\})} + \|\chi_{R,R_2} w\|_{Y(\Rm)}^{7/3}+ \eta^{4/3} \|\chi_{R,R_2} w\|_{Y(\Rm)} +\varepsilon \right).
 \end{align*}
 We choose $\eta$ to be a sufficiently small number such that $C_1 \eta^{4/3} < 1/(4C_1) < 1/2$, thus 
  \begin{align} \label{recurrence w Y} 
  \|\chi_{R,R_2} w\|_{Y(\Rm)} \leq 2C_1 \left(\left|\vec{\tau}(R)\right| + \|G\|_{L^2(\{s: R<|s|<R_2\})} + \|\chi_{R,R_2} w\|_{Y(\Rm)}^{7/3}+\varepsilon \right).
 \end{align}
 We choose small constants $\varepsilon_1 = \beta$ such that 
 \begin{align*}
  2C_1(8C_1\beta)^{4/3} < \frac{1}{4C_1} < \frac{1}{4} \quad \Longrightarrow \quad  8C_1 \beta > 2C_1(3\beta + (8C_1\beta)^{7/3}). 
 \end{align*}
 As a result, if $\|G\|_{L^2(\{s: R<|s|<R_2\})} \leq \beta$, then a continuity argument in $R$ shows that 
 \[
  \|\chi_{R,R_2} w\|_{Y(\Rm)} < 8 C_1 \beta. 
 \]
 Inserting this to \eqref{recurrence w Y} and using the choice of $\beta$, we obtain 
 \[
  \|\chi_{R,R_2} w\|_{Y(\Rm)} \leq 2C_1 \left(\left|\vec{\tau}(R)\right| + \|G\|_{L^2(\{s: R<|s|<R_2\})} +\varepsilon \right) + \frac{1}{4} \|\chi_{R,R_2} w\|_{Y(\Rm)},
 \]
 which implies 
 \[
   \|\chi_{R,R_2} w\|_{Y(\Rm)} \leq \frac{8}{3}C_1 \left(\left|\vec{\tau}(R)\right| + \|G\|_{L^2(\{s: R<|s|<R_2\})} +\varepsilon \right).
 \]
 An application of the nonlinear radiation profile shows (we apply Lemma \ref{scatter profile of nonlinear solution} on $w$ and recall the choice of $\beta, \eta$)
 \begin{align*}
 \sqrt{2\sigma_4} \|G\|_{L^2(\{s:|s|>R\})} & \leq \|\chi_{R} (F(u)-e(x,t)-F(S))\|_{L^1 L^2 (\Rm \times \Rm^3)}\nonumber \\
 & \leq \|\chi_{R,R_2} (F(w+S)-F(S))\|_{L^1 L^2(\Rm \times \Rm^3)} + \varepsilon \nonumber\\
 & \leq C_1\left(\|\chi_{R,R_2} w\|_{Y(\Rm)}^{7/3} + \|\chi_{R,R_2} S\|_{Y(\Rm)}^{4/3} \|\chi_{R,R_2} w\|_{Y(\Rm)}\right) + \varepsilon \\
 & \leq C_1 \left((8C_1 \beta)^{4/3} + \eta^{4/3}\right) \|\chi_{R,R_2} w\|_{Y(\Rm)} + \varepsilon \nonumber \\
 & \leq \frac{3}{8C_1}\cdot \frac{8}{3}C_1 \left(\left|\vec{\tau}(R)\right| + \|G\|_{L^2(\{s: R<|s|<R_2\})} +\varepsilon \right) + \varepsilon\\
 & \leq \left|\vec{\tau}(R)\right| + \|G\|_{L^2(\{s: R<|s|<R_2\})} +2\varepsilon.
\end{align*}
This immediately gives (please note that $\sqrt{2\sigma_4} > 5$)
\[
 \|G\|_{L^2(\{s:R<|s|<R_2\})} \leq \frac{1}{4}\left|\vec{\tau}(R)\right| + \frac{1}{2} \varepsilon \leq \frac{3}{4}\beta.
\]
A continuity argument in $R$ shows that $\|G\|_{L^2(\{s: R_1<|s|<R_2\})} \leq 3\beta/4$. Thus the inequalities above hold for all $R\in [R_1,R_2)$. A combination of these inequalities, the Strichartz estimates with Remark \ref{double tail G} finishes the proof. 
\end{proof}

\begin{lemma} \label{lemma connection 2} 
 Let $\eta$ be the constant in Lemma \ref{lemma connection}. There exists an absolute positive constant $\varepsilon_2$ such that if $19 R_2/20 \leq R_1<R_2$ and
 \begin{itemize} 
  \item $u$ is an exterior solution to (CP1) and $S$ is an exterior solution to the equation 
 \[
  (\partial_t^2 -\Delta) S= F(S) + e(x,t). 
 \] 
 with $\|\chi_{R_1} u\|_{Y(\Rm)}, \|\chi_{R_1} S\|_{Y(\Rm)}, \|\chi_{R_1} e(x,t)\|_{L^1 L^2} < +\infty$. 
 \item Solutions $u$, $S$ are asymptotically equivalent to each other in $\Omega_{R_1}$. 
 \item $u$, $S$ and $w=u-S$ satisfy the following inequalities
 \begin{align*}
  \varepsilon \doteq \|(w(\cdot,0),w_t(\cdot,0))\|_{\mathcal{H}(R_2)}  + \|\chi_{R_1,R_2} e(x,t)\|_{L^1 L^2(\Rm \times \Rm^3)}  & \\
 + \|\chi_{R_2}\left(F(u) - F(S)- e(x,t)\right)\|_{L^1 L^2(\Rm \times \Rm^3)} & \leq \varepsilon_2; \\
   \|\chi_{R_1,R_2} S\|_{Y(\Rm)} & \leq \eta; 
 \end{align*}
 \end{itemize}
 Then we have 
 \begin{align*}
  \|\chi_{R_1} w\|_{Y(\Rm)} + \sup_{t\in \Rm} \left\|\vec{w}(t)\right\|_{\mathcal{H}(R_1+|t|)} \lesssim_1 \varepsilon. 
 \end{align*}
\end{lemma}

\begin{proof}
 The proof is similar to Lemma \ref{lemma connection}. Let $w_L$, $G$ and $\vec{\tau}(R)$ be the linear free wave, radiation profile and radiation residue with initial data $(w(\cdot,0), w_t(\cdot,0))$. First of all, we assume that $R_1 \geq (19/20) R_2$ and obtain for $R\in [R_1,R_2)$ that
 \begin{align*}
  |\tau_1(R)| & = \left|R^{-1/2} \int_{-R}^{R} G(s) {\rm d} s\right| \\
  & \leq R^{-1/2} \left|\int_{-R_2}^{R_2} G(s) {\rm d} s\right| + R^{-1/2} \int_{R<|s|<R_2} |G(s)| {\rm d} s\\
  & \leq \left(\frac{20}{19}\right)^{1/2} |\tau_1(R_2)|+ \left(\frac{2(R_2 - R)}{R}\right)^{1/2} \|G\|_{L^2(\{s: R<|s|<R_2\})}\\
  & \leq C \|\vec{w}(0)\|_{\mathcal{H}(R_2)} + (1/3)\|G\|_{L^2(\{s: R<|s|<R_2\})},
 \end{align*}
 and that 
 \begin{align*}
   |\tau_2(R)| & = \left|\frac{\sqrt{3}}{R^{3/2}} \int_{-R}^{R} s G(s) {\rm d} s\right| \\
  & \leq \frac{\sqrt{3}}{R^{3/2}} \left|\int_{-R_2}^{R_2} s G(s) {\rm d} s\right| + \frac{\sqrt{3}}{R^{3/2}} \int_{R<|s|<R_2} |sG(s)| {\rm d} s\\
  & \leq \left(\frac{20}{19}\right)^{3/2} |\tau_2(R_2)| + \left(\frac{6(R_2 - R)R_2^2}{R^3}\right)^{1/2} \|G\|_{L^2(\{s: R<|s|<R_2\})}\\
  & \leq C \|\vec{w}(0)\|_{\mathcal{H}(R_2)} + (2/3) \|G\|_{L^2(\{s: R<|s|<R_2\})}. 
 \end{align*}
 Here $C$ is an absolute constant. Combining these with Lemma \ref{lemma inequality}, we obtain 
  \begin{align*}
  \|\chi_{R} w\|_{Y(\Rm)} \leq C_1 \left( 2\|G\|_{L^2(\{s: R<|s|<R_2\})} + \|\chi_{R,R_2} w\|_{Y(\Rm)}^{7/3}+ \eta^{4/3} \|\chi_{R,R_2} w\|_{Y(\Rm)} +C_2 \varepsilon \right).
 \end{align*}
Here $C_2> 1$ is an absolute constant. By choosing the same constants $\eta$, $\beta$ as in Lemma \ref{lemma connection} and applying a continuity argument in $R$, we obtain that if $\|G\|_{L^2(\{s: R<|s|<R_2\})} \leq \beta$ and $\varepsilon \leq \varepsilon_2 \doteq C_2^{-1} \beta$, then 
 \[
   \|\chi_{R,R_2} w\|_{Y(\Rm)} \leq \frac{8}{3}C_1 \left(2\|G\|_{L^2(\{s: R<|s|<R_2\})} +C_2\varepsilon \right) \leq 8C_1 \beta. 
 \]
 As in the proof of Lemma \ref{lemma connection}, an application of the nonlinear radiation profile shows
 \begin{align*}
 \sqrt{2\sigma_4} \|G\|_{L^2(\{s:|s|>R\})} & \leq  C_1\left(\|\chi_{R,R_2} w\|_{Y(\Rm)}^{7/3} + \|\chi_{R,R_2} S\|_{Y(\Rm)}^{4/3} \|\chi_{R,R_2} w\|_{Y(\Rm)}\right) + \varepsilon \\
 & \leq C_1 \left((8C_1 \beta)^{4/3} + \eta^{4/3}\right) \|\chi_{R,R_2} w\|_{Y(\Rm)} + \varepsilon \nonumber \\
 & \leq \frac{3}{8C_1}\cdot \frac{8}{3}C_1 \left(2\|G\|_{L^2(\{s: R<|s|<R_2\})} + C_2 \varepsilon \right) + \varepsilon\\
 & \leq 2\|G\|_{L^2(\{s: R<|s|<R_2\})} +2 C_2 \varepsilon. 
\end{align*}
It immediately follows that 
\[
 \|G\|_{L^2(\{s: R<|s|<R_2\})} \leq \frac{2}{3} C_2 \varepsilon \leq \frac{2}{3}\beta. 
\]
A continuity argument then yields that $\|G\|_{L^2(\{s: R_1<|s|<R_2\})} \leq 2\beta/3$. As a result, the estimates given above hold for $R=R_1$. A combination of these estimates with Remark \ref{double tail G} finishes the proof. 
\end{proof}

\section{Non-radiative solutions}

In this section we investigate the properties of all possible non-radiative solutions to (CP1). The results in this section help us classify all nontrivial non-radiative solutions, up to a rescaling. The results in the first four subsection holds for radial non-radiative solutions to wave equation with more general nonlinear terms. In fact we only need to assume the energy-critical/symmetric conditions on the nonlinear term $F$, i.e. 
\begin{align*}
 &F(0) = 0;& & |F(u)-F(v)| \lesssim_1 (|u|^{4/3}+|v|^{4/3})|u-v|;& &F(-u) = -F(u). 
\end{align*}

\subsection{Extension and maximal domain}

We first show that the non-radiative extension of a given non-radiative solution is unique and describe the possible behaviour of the extension with a maximal domain. 

\begin{lemma} \label{uniqueness of extension}
 Let $u$ and $\tilde{u}$ be two non-radiative exterior solutions defined in $\Omega_r$ and $\Omega_{\tilde{r}}$ respectively. If there exists a radius $R > \max\{r,\tilde{r}\}$ such that $u \equiv \tilde{u}$ in the exterior region $\Omega_R$. Then we must have 
 \[
  u(x,t) = \tilde{u}(x,t), \qquad |x|>|t|+\max\{r,\tilde{r}\}. 
 \]
\end{lemma}
\begin{proof}
 It is a direct consequence of Lemma \ref{lemma connection 2}. Please note that Lemma \ref{finite eq scattering} guarantees that $\|\chi_r u\|_{Y(\Rm)} < +\infty$ and $\|\chi_{\tilde{r}}\tilde{u}\|_{Y(\Rm)} < +\infty$. 
\end{proof}

\begin{lemma} \label{extension with finite Y} 
 Let $u$ be a non-radiative exterior solution in $\Omega_R$ to (CP1) with $R>0$. Then we may extend it to a non-radiative exterior solution defined in $\Omega_r$ with $r<R$.
\end{lemma}
\begin{proof}
 We first extend the domain $u$ to a slightly large domain $\Omega_{r_1}$, which is not necessarily a non-radiative solution. Let $\tilde{u}$ be the solution to the linear wave equation 
 \[
  \left\{\begin{array}{l} \partial_t^2 \tilde{u} - \Delta \tilde{u} = \chi_R F(u), \quad (x,t) \in \Rm^5 \times \Rm;\\
  (u,u_t)|_{t=0} = (u_0,u_1).\end{array}\right.
 \]
 Here $(u_0,u_1)\in \dot{H}^1\times L^2$ matches the initial data of $u$ in the exterior region $\{x: |x|>R\}$. The finite of propagation implies that $u\equiv \tilde{u}$ in $\Omega_R$. By the Strichartz estimate we also have $\|\tilde{u}\|_{Y(\Rm)}<+\infty$. Thus if $r_1$ is slightly smaller than $R$, than $\tilde{u}$ solves the approximated equation 
 \[
  \partial_t^2 \tilde{u} - \Delta \tilde{u} = F(\tilde{u}) - \chi_{r_1,R} F(\tilde{u})
 \]
 in the exterior region $\Omega_{r_1}$ with $\|\chi_{r_1,R}F(\tilde{u})\|_{L^1 L^2} \leq \|\chi_{r_1,R}\tilde{u}\|_{Y(\Rm)}^{7/3} \ll 1$. By the perturbation theory, the exterior solution to (CP1) with initial data $(u_0,u_1)$ in $\Omega_{r_1}$ exists for all $t\in \Rm$. It coincides with $u$ in $\Omega_R$ by finite speed of propagation. We still call it $u$. It clearly satisfies $\|\chi_{r_1} u\|_{Y(\Rm)} < +\infty$. Let its nonlinear radiation profiles be $G_{\pm}^u \in L^2([r_1,+\infty))$. Denote $M=\|\chi_{R} u\|_{Y(\Rm)}$ and let $\varepsilon = \varepsilon(M)$ be a small positive constant to be determined later. We choose $r\in [\max\{r_1,R/2\},R)$ to be a radius slightly smaller than $R$ such that 
 \begin{align*}
  &\|G_\pm^u\|_{L^2([r,R])} \leq \varepsilon;& &\|\chi_{r,R} u\|_{Y(\Rm)} \leq \varepsilon;
 \end{align*}
 and let $X = \{G\in L^2(\{s: r<|s|<R\}: \|G\|_{L^2}\leq 4\varepsilon)$ to a complete metric space. Given any $G\in X$, we extend it to a function in $L^2(\Rm)$ by defining $G(s) =0$ for $s>R$ and defining $G(s)$ for $|s|<r$ in a way such that 
 \begin{align*}
  &\int_{-R}^R G(s) {\rm d}s = 0;& & \int_{-R}^R s G(s) {\rm d} s = 0; & & \|G\|_{L^2(\Rm)} \lesssim_1 \|G\|_X \lesssim_1 \varepsilon. 
 \end{align*}
 Let $(v_0,v_1)\in \dot{H}^1\times L^2$ be initial data such that the radiation profile of $(v_0,v_1)-(u_0,u_1)$ is exactly $G$. Clearly $(u_0,u_1)$ coincide with $(v_0,v_1)$ in the exterior region $\{x: |x|>R\}$. By the perturbation theory, if $\varepsilon$ is sufficiently small, then the corresponding exterior solution $v$ of (CP1) in $\Omega_r$ with initial data $(v_0,v_1)$ is defined for all $t\in \Rm$ and satisfies
 \[
  \|\chi_r(u-v)\|_{Y(\Rm)} \lesssim_M \varepsilon. 
 \]  
 Let $\tilde{G}_\pm$ be the nonlinear radiation profile of the solution to the inhomogeneous wave equation $(\partial_t^2-\Delta) w = F(v) - F(u)$ in the exterior region $\Omega_r$ with zero initial data. Then we have 
 \[
  \|\tilde{G}_\pm \|_{L^2([r,+\infty))} \lesssim_1 \|\chi_r (F(v)-F(u))\|_{L^1 L^2}. 
 \]
 By finite speed of wave propagation, $u$ and $v$ coincide in the region $\Omega_R$. Thus 
 \[
  \|\tilde{G}_\pm \|_{L^2([r,+\infty))} \lesssim_1 \|\chi_{r,R}(F(u)-F(v))\|_{L^1 L^2} \lesssim_M \varepsilon^{7/3}. 
 \]
 By linearity the nonlinear radiation profile $G_\pm$ of $v$ is given by 
 \[
  \left\{\begin{array}{ll} G_+(s) = G_+^u (s) + G(-s) + \tilde{G}_+(s), & s > r; \\ G_-(s) = G_-^u(s) + G(s) + \tilde{G}_-(s), & s>r. \end{array}\right.
 \]
 Our assumptions and finite speed of propagation implies that all the radiation profiles are zero for $r>R$. Thus $v$ is a non-radiative solution if and only if $G(s)$ is a fixed point of the map $\mathbf{T}: X \rightarrow L^2(\{s: r<|s|<R\})$ defined by
 \[
  \mathbf{T} G = \left\{\begin{array}{ll} -G_-^u(s) - \tilde{G}_-(s), & r<s<R; \\ -G_+^u(-s) - \tilde{G}_+(-s), & -R<s<-r. \end{array}\right.
 \] 
 Thus in order to complete the proof, its suffices to verify that $\mathbf{T}$ is a contraction map from $X$ to itself if $\varepsilon$ is sufficiently small. Indeed we have 
 \[
  \|\mathbf{T} G\|_{L^2(\{s: r<|s|<R\})} \leq \sum_\pm \left(\|G_\pm^u\|_{L^2([r,R])} + \|\tilde{G}_\pm\|_{L^2([r,R])}\right) \leq 2\varepsilon + C(M) \varepsilon^{7/3}< 4\varepsilon,
 \]
 as long as $\varepsilon$ is sufficiently small. Thus $\mathbf{T} X \subset X$. In addition, if $G_1, G_2 \in X$ and $v^j$, $\tilde{G}_\pm^{j}$ are their corresponding solutions/nonlinear radiation profiles as defined above, then
 \begin{align*}
  \|\mathbf{T} (G_1-G_2)\|_{L^2(\{s: r<|s|<R\})} & \leq \sum_{\pm} \|\tilde{G}_\pm^{1} - \tilde{G}_\pm^2\|_{L^2([r,R])}\\
  & \lesssim_1 \|\chi_r (F(v^1)-F(v^2))\|_{L^1 L^2(\Rm \times \Rm^5)} \\
  & \lesssim_1 \|\chi_{r,R} (F(v^1)-F(v^2))\|_{L^1 L^2(\Rm \times \Rm^5)} \\
  & \lesssim_M \varepsilon^{4/3} \|\chi_{r,R} (v^1 - v^2)\|_{Y(\Rm)}. 
 \end{align*}
 By the Strichartz estimates, we also have 
 \begin{align*}
  \|\chi_{r,R} (v^1 - v^2)\|_{Y(\Rm)} &\leq \|\chi_r (v_L^1 - v_L^2)\|_{Y(\Rm)} + C_1 \|\chi_r (F(v^1)-F(v^2))\|_{L^1 L^2(\Rm \times \Rm^5)}\\
  & \leq C_2 \|G_1-G_2\|_{L^2(\{s: r<|s|<R\})} + C_3(M) \varepsilon^{4/3} \|\chi_{r,R} (v^1 - v^2)\|_{Y(\Rm)}. 
 \end{align*}
 Here $C_1, C_2$ are absolute constants and $C_3(M)$ is a constant determined by $M$; $v_L^j$'s are the linear free wave with the same initial data as $v^j$. Therefore if $\varepsilon$ is sufficiently small, then 
 \[
  \|\chi_{r,R} (v^1 - v^2)\|_{Y(\Rm)} \lesssim_1 \|G_1-G_2\|_{L^2(\{s:r<|s|<R\})}. 
 \]
 In summary we have
 \[
  \|\mathbf{T} (G_1-G_2)\|_{L^2(\{s: r<|s|<R\})} \lesssim_M \varepsilon^{4/3}\|G_1-G_2\|_{L^2([r,R])} < \frac{1}{2} \|G_1-G_2\|_{L^2(\{s:r<|s|<R\})},
 \]
 as long as $\varepsilon$ is small enough, and finish the proof. 
\end{proof}

\begin{proposition}[Maximal domain] \label{p maximal domain}
 Let $u$ be a radial non-radiative exterior solution defined in $\Omega_R$ (with initial data in $\mathcal{H}(R)$). Then we may extend the domain of $u$ to $\Omega_r$ such that either of the following holds 
 \begin{itemize}
  \item $r=0$ and $u$ becomes a non-radiative exterior solution defined in $\Omega_0$ with $\|\chi_0 u\|_{Y(\Rm)} < +\infty$. We call it a scattering non-radiative solution; 
  \item $u$ satisfies $\|\chi_r u\|_{Y(\Rm)} = +\infty$; in addition, $u$ is a non-radiative exterior solution in $\Omega_{r_1}$ for any $r_1>r$. We call it a blow-up non-radiative solution. In this case the initial data of $u$ are not necessary contained in $\mathcal{H}(r)$ but are contained in $\mathcal{H}(r_1)$ for any $r_1>r$. 
 \end{itemize}
\end{proposition}
\begin{proof}
 We consider all possible extensions $(\tilde{u}, \Omega_{\tilde{r}})$ of $u$ such that $\tilde{u}$ is still a radial non-radiative exterior solution. Lemma \ref{uniqueness of extension} implies that any two such extensions coincide with each other in the overlapping region of their domains. As a result, we may naturally extend the domain of $u$ to $\Omega_r$, where $r$ is the greatest lower bound of all such radii $\tilde{r}$. It remains to show that if $\|\chi_r u\|_{Y(\Rm)} < +\infty$, then 
 \begin{itemize}
  \item[(I)] $r$ must be zero;
  \item[(II)] $u$ is a non-radiative exterior solution in $\Omega_0$ with initial data in $\dot{H}^1\times L^2$. 
 \end{itemize}
 We first show that if $\|\chi_r u\|_{Y(\Rm)} < +\infty$, then $u$ is a non-radiative exterior solution in $\Omega_r$. Our assumption guarantees that there exists a sequence $r_k \rightarrow r^+$, such that $u$ is a non-radiative exterior solution in $\Omega_{r_k}$. The initial data of the restriction to $\Omega_{r_k}$ are uniquely determined in the space $\mathcal{H}(r_k)$. In other words, both the radiation profile $G(s)$ for $|s|>r_k$ and the radiation residue $\vec{\tau}(r_k)$ are uniquely determined. In addition, these data for different radii $r_k$ must match each other. More precisely, $G(s)$ can be consistently defined for $|s|>r$ and 
 \begin{align*}
 &r_k^{1/2} \tau_1(r_k) - r_j^{1/2} \tau_1(r_j) = -\int_{r_j<|s|<r_k} G(s) {\rm d} s,& & r_j < r_k;\\
  & r_k^{3/2} \tau_2(r_k) - r_j^{3/2} \tau_2(r_j) = \sqrt{3} \int_{r_j<|s|<r_k} s G(s) {\rm d} s,& & r_j<r_k.  
 \end{align*}
 An application of Lemma \ref{scatter profile of nonlinear solution} shows that $\|G\|_{L^2(\{s:|s|>r\})} < +\infty$. This also implies that $r_k^{1/2} \tau_1(r_k)$ and $r_k^{3/2} \tau_2(r_k)$ converge as $k\rightarrow +\infty$, which guarantees that the initial data are contained in $\mathcal{H}(r)$ as long as $r>0$. On the other hand, we claim that if $r=0$, then 
 \begin{align*}
   \lim_{j\rightarrow +\infty}  \vec{\tau} (r_j) = 0;
 \end{align*}
 This implies that the initial data are exactly the element in $\dot{H}^1\times L^2$ with radiation profile $G(s)$. Now let us verify the claim. First of all, the Strichartz estimates and our assumption $\|\chi_0 u\|_{Y(\Rm)} < +\infty$ implies that the free wave $u_L$ with the same initial data as $u$ satisfies 
 \[
  \limsup_{j\rightarrow +\infty} \|\chi_{r_j} u_L\|_{Y([-r_j,r_j])} = 0. 
 \]
 On the other hand, if we let $\tilde{u}_L$ be the free wave with radiation profile $G(s) \in L^2(\Rm)$, then we may recall \eqref{free wave by radiation profile} and rewrite $w_L$ in the following form, as long as $r>|t|+r_j$
 \[
  u_L(r,t) = \tilde{u}_L(r,t) + \frac{1}{r^3} \left[\frac{r_j^{3/2} \tau_2(r_j)}{\sqrt{3}} - \int_{-r_j}^{r_j} s G(s) {\rm d} s\right] + \frac{t}{r^3} \left[r_j^{1/2} \tau_1(r_j) + \int_{-r_j}^{r_j} G(s) {\rm d} s\right].
 \]
 The decay of $Y$ norm implies that 
 \begin{align*}
  &\lim_{j\rightarrow +\infty} \left[\frac{1}{\sqrt{3}}\tau_2(r_j) - r_j^{-3/2} \int_{-r_j}^{r_j} sG(s) {\rm d} s\right] = 0& &\Longrightarrow & &\lim_{j\rightarrow +\infty}\tau_2(r_j) = 0; \\
  &\lim_{j\rightarrow +\infty} \left[\tau_1(r_j) + r_j^{-1/2} \int_{-r_j}^{r_j} G(s) {\rm d} s\right] = 0& &\Longrightarrow & &\lim_{j\rightarrow +\infty}\tau_1(r_j) = 0.
 \end{align*}
 In summary, the initial data are always contained in the space $\mathcal{H}(r)$. Finite speed of propagation then verifies that $u$ is exactly the exterior solution to (CP1) with these pair of initial data. Our assumption $\|\chi_r u\|_{Y(\Rm)} < +\infty$ implies that $u$ scatters to some linear free wave in the exterior region $\Omega_r$, which implies that $u$ is non-radiative. Finally Lemma \ref{extension with finite Y} implies that $r=0$, otherwise we may extend the domain to a larger exterior region, which contradicts with the fact that $r$ is the greatest lower bound. This eventually finishes the proof. 
\end{proof}

\begin{remark}
 A similar result to Proposition \ref{p maximal domain} for more general energy-critical wave equations was proved in Collot-Duyckaerts-Kenig-Merle \cite{classNR}. They proved that a blow-up non-radiative solution defined in a maximal domain $\Omega_R$ satisfies at least one of the followings: I. the Strichartz norm blows up in $\Omega_R$; II. the $\mathcal{H}(R)$ norm of the initial data blows up. In this particular case of classic focusing wave equation in 5D, Proposition \ref{p maximal domain} implies that the blow-up of non-radiative solutions always comes from the blow-up of the Strichartz norm. 
\end{remark}

\subsection{Radiation residues}

In this subsection we discuss the radiation residues of initial data for non-radiative solutions. 

\begin{lemma} \label{small tau}
 Let $u$ be a radial non-radiative solution to (CP1) and $\vec{\tau}(r)$ be the radiation residue of its initial data. If $|\vec{\tau}(r)| < \beta$ for all $r\geq R$, where $\beta$ is the positive constant in Lemma \ref{lemma connection}, then 
 \[
  \|\vec{u}(0)\|_{\mathcal{H}(R)} + \|\chi_R u\|_{Y(\Rm)} \lesssim_1 \left|\vec{\tau}(R)\right|. 
 \]
\end{lemma}
\begin{proof}
 Let us compare the non-radiative solution $u$ with the zero solution $S\equiv 0$. For large radius $R_2$, we have
 \[
  \lim_{R_2\rightarrow +\infty} \left(\left\|\vec{u}(0)-\vec{S}(0)\right\|_{\mathcal{H}(R_2)} + \|\chi_{R_2}(F(u)-F(S))\|_{L^1 L^2(\Rm \times \Rm^5)}\right) = 0. 
 \]
 Thus we may apply Lemma \ref{lemma connection} in the exterior region $\Omega_R$, make $R_2 \rightarrow +\infty$ and finish the proof. 
\end{proof}

\begin{lemma} \label{small tau zero}
 Let $u$ be a radial non-radiative solution of (CP1) with a maximal domain $\Omega_{R_0}$ and $\vec{\tau}(R)$ be the radiation residues of its initial data. If $|\vec{\tau}(R)| < \beta$ for all $R>R_0$, where $\beta$ is the positive constant in Lemma \ref{lemma connection}, then $u\equiv 0$.
\end{lemma}
\begin{proof}
 We start by applying Lemma \ref{small tau} and deduce
 \begin{equation} \label{upper bound of exterior region trivial nonradi} 
  \|\vec{u}(0)\|_{\mathcal{H}(R)} + \|\chi_R u\|_{Y(\Rm)} \lesssim_1 \left|\vec{\tau}(R)\right|, \qquad \forall R>R_0. 
 \end{equation}
 The uniform boundedness of $\|\chi_R u\|_{Y(\Rm)}$ then implies that $u$ must be a scattering non-radiative solution. Thus we must have that $R_0 = 0$ and that $|\vec{\tau}(R)| \rightarrow 0$ as $R\rightarrow 0^+$. Inserting this into \eqref{upper bound of exterior region trivial nonradi}, we immediately conclude that $u \equiv 0$. 
\end{proof}

\begin{remark} \label{no zero small}
 The same argument shows that if $R_1>R_0$ and 
 \begin{align*}
  &\vec{\tau}(R_1) = 0;& &\vec{\tau}(r) < \beta, \quad \forall r\in (R_1,+\infty);
 \end{align*}
 then $u$ must be zero as well. 
\end{remark}

\begin{lemma} \label{identification tau}
 Let $u$, $u_\ast$ be two radial non-radiative solution of (CP1) and $\vec{\tau}(r)$, $\vec{\tau}_\ast(r)$ be the radiation residues of their initial data. There exists a constant $\rho_2 > 0$, such that if $\vec{\tau}$, $\vec{\tau}_\ast$ and a radius $R>0$ satisfies
 \begin{align*}
  & \vec{\tau}(R) = \vec{\tau}_\ast(R); & & |\vec{\tau} (r)|, |\vec{\tau}_\ast (r)| < \rho_2, \quad \forall r\geq R; 
 \end{align*}
 then $u = u_\ast$. 
\end{lemma}
\begin{proof}
By Lemma \ref{small tau}, we have $\|\chi_R u\|_{Y(\Rm)} \lesssim_1 |\vec{\tau}(R)| \ll 1$ when $\rho_2$ is sufficiently small.  An application of Lemma \ref{lemma connection} on $u_\ast$ and $u$ then completes the proof. 
\end{proof}

\begin{lemma} \label{monotonicity of tau}
 Let $u$ be a non-trivial radial non-radiative solution to (CP1) and $\vec{\tau}(r)$ be the radiation residue of its initial data. There exists a small constant $\rho_0 > 0$, such that if $|\vec{\tau}(r)| < \rho_0$ for all $r\in [R,+\infty)$, then $|\vec{\tau}(r)|$ must be a strictly decreasing function of $r$ in the interval $[R,+\infty)$.  
\end{lemma}
\begin{proof}
 Let $\rho_0 < \beta$ be a small constant to be determined later. Lemma \ref{small tau} immediately gives 
 \[
  \|\vec{u}(0)\|_{\mathcal{H}(r_0)} + \|\chi_{r_0} u\|_{Y(\Rm)} \lesssim_1 |\vec{\tau}(r_0)|, \qquad r_0\geq R. 
 \]
 Since $u$ is non-trivial, by Lemma \ref{small tau zero} we must have that $R>0$. The uniform boundedness of $\|\chi_{r_0} u\|_{Y(\Rm)}$ also implies that the maximal domain of $u$ is larger than $\Omega_R$. A combination of the small data theory and the point-wise estimate $|v(r)| \lesssim_1 r^{-3/2} \|v\|_{\dot{H}^1}$ for radial $\dot{H}^1(\Rm^5)$ functions then yields 
 \begin{equation} \label{pointwise u small tau}
  |u(r,t)| \lesssim_1 r^{-3/2} \|u(\cdot,t)\|_{\mathcal{H}(|t|+r_0)} \lesssim_1 r^{-3/2} |\vec{\tau}(r_0)|, \qquad r\geq r_0+|t|.
 \end{equation}
Corollary \ref{continuity of non-radiative} and the formula \eqref{radiation profile initial non-radiative} then gives the continuity of the radiation profile $G_\pm(s)$ of $\vec{u}(0)$ for $|s|\geq R$ and the following upper bound 
 \begin{equation} \label{decay of G non-ra}
  |G_\pm (s)| \lesssim_1 s^{-1/2} |\vec{\tau}(r_0)|^{7/3}, \qquad \forall s\geq r_0\geq R. 
 \end{equation}
 Now we consider the derivatives
 \begin{align*}
  \tau'_1(r) & = \frac{\rm d}{{\rm d} r} \left(\frac{-1}{r^{1/2}} \int_{-r}^r G_-(s) {\rm d} s\right) = - \frac{1}{2r} \tau_1(r) - \frac{1}{r^{1/2}} \left[G_-(r) + G_+(r)\right];\\
  \tau'_2(r) & = \frac{\rm d}{{\rm d} r} \left(\frac{\sqrt{3}}{r^{3/2}} \int_{-r}^r s G_-(s) {\rm d} s\right) = -\frac{3}{2r} \tau_2(r) + \frac{\sqrt{3}}{r^{1/2}} \left[G_-(r) - G_+(r)\right].
 \end{align*}
 It immediately follows that 
 \begin{align} \label{derivative of tau square}
  \frac{\rm d}{{\rm d} r} |\vec{\tau}(r)|^2  = 2 \vec{\tau}(r) \cdot \frac{\rm d}{{\rm d}r}\vec{\tau}(r)  \leq \frac{1}{r} \left[-\tau_1^2 (r)-3\tau_2^2(r) + C |\vec{\tau}(r)|^{10/3}\right], \quad r \geq R. 
 \end{align}
 Here $C$ is an absolute constant. The derivative is clearly negative as long as $|\vec{\tau}(r)| < \rho_0$ is sufficiently small. Please note that $\vec{\tau}(r)$ must be nonzero for $r\geq R$, according to Remark \ref{no zero small}. 
\end{proof}

Next we further discuss the asymptotic behaviour of the radiation residue $\vec{\tau}(r)$ by using the symmetric property $F(u) = - F(-u)$. 

\begin{lemma} \label{detailed asymptotic behaviour}
  Let $u$ be a non-trivial radial non-radiative solution to (CP1) and $G$, $\vec{\tau}(r)$ be the radiation profile and residue of its initial data, respectively. There exists a small constant $\rho_0 > 0$, such that if $|\vec{\tau}(r)| < \rho_0$ for all $r\in [R,+\infty)$, then $|\tau_j (r)|$ must be either identically zero or a strictly decreasing function of $r$ in the interval $[R,+\infty)$. In addition, both $\int_{-r}^r G(s) {\rm  d}s$ and $\int_{-r}^r s G(s) {\rm d} s$ converge as $r\rightarrow +\infty$ with
 \begin{align} \label{equivalence asymptotic tau}
  &\int_{-r}^r G(s) {\rm d} s \simeq_1 \int_{-R}^R G(s) {\rm d} s;& & \int_{-r}^r s G(s) {\rm d} s \simeq_1 \int_{-R}^R s G(s) {\rm d} s.
 \end{align}
\end{lemma}
\begin{proof}
 We first recall \eqref{derivative of tau square} and deduce that if $\rho_0$ is sufficiently small, then 
 \[ 
  \frac{\rm d}{{\rm d} r} |\vec{\tau}(r)|^2 \leq - \frac{9}{10 r} |\vec{\tau}(r)|^2, \qquad r\geq R.  
 \]
 It follows that $|\vec{\tau}(r)|^2 \leq (r/R)^{-9/10} |\vec{\tau}(R)|^2$. Combining this with \eqref{derivative of tau square}, we obtain 
 \[ 
  \frac{\rm d}{{\rm d} r}\left(r |\vec{\tau}(r)|^2 \right) \lesssim_1 |\vec{\tau}(r)|^{10/3} \lesssim_1 (r/R)^{-3/2} |\vec{\tau}(R)|^{10/3}, \qquad r\geq R. 
 \]
 Thus we have 
 \[
  r |\vec{\tau}(r)|^2 \lesssim_1 R |\vec{\tau}(R)|^2 \qquad \Longrightarrow \qquad |\vec{\tau}(r)| \lesssim_1 (r/R)^{-1/2} |\vec{\tau}(R)|. 
 \]
 A combination of this with \eqref{pointwise u small tau} yields that 
 \begin{equation} \label{pointwise u upgrade}
  |u(r,t)| \lesssim_1 r^{-3/2} (r-|t|)^{-1/2} R^{1/2} |\vec{\tau}(R)|, \qquad r\geq R+|t|.
 \end{equation}
 Now let us use the notations 
 \begin{align*}
  &\bar{u}(x,t) = u(x,-t);& & u_e = \frac{u+\bar{u}}{2};& & u_o = \frac{u-\bar{u}}{2}.
 \end{align*}
 There are all non-radiative solutions and solve the following wave equation, respectively. 
 \begin{align*}
  &\square \bar{u} = F(\bar{u});& & \square u_e = \frac{F(u)+F(\bar{u})}{2}; & &\square u_o = \frac{F(u)-F(\bar{u})}{2}. 
 \end{align*}
 For convenience we use the notation $\square = \partial_t^2 -\Delta$. The corresponding radiation profiles of initial data of $u_e$ and $u_o$, i.e. $(u_0,0)$ and $(0,u_1)$, are given by 
 \begin{align*}
  & G_e (s) = \frac{G(s)-G(-s)}{2}; & & G_o (s) = \frac{G(s)+G(-s)}{2}. 
 \end{align*}
 Please note that $G_e$, $G_o$ are odd and even functions of $s$, respectively. These are inconsistent with their subscripts. The corresponding radiation residues are given by 
 \begin{align*}
  &\vec{\tau}_e (r) = (0, \tau_2(r)); & &\vec{\tau}_o (r) = (\tau_1(r),0). 
 \end{align*}
 As a result, we have 
 \begin{align*}
  \sigma_4^{-1} \|\vec{u}_e(0)\|_{\mathcal{H}(r)}^2 = 4 \int_r^\infty |G_e(s)|^2 {\rm d} s + |\tau_2 (r)|^2; \\
  \sigma_4^{-1} \|\vec{u}_o(0)\|_{\mathcal{H}(r)}^2 = 4 \int_r^\infty |G_o(s)|^2 {\rm d} s + |\tau_1 (r)|^2. 
 \end{align*}
 Next let us consider 
 \[
  M(r) \doteq \sup_{t\in \Rm} \|\vec{u}_e(t)\|_{\mathcal{H}(|t|+r)} + \|\chi_r u_e\|_{Y(\Rm)}, \qquad r\geq R;
 \]
 which must be finite, because $u$ is a small solution in $\Omega_R$. It then follows from the Strichartz estimate and Lemma \ref{tail G} that 
 \begin{align*}
  M(r) &\lesssim_1 \|\vec{u}_e(0)\|_{\mathcal{H}(r)} +  \frac{1}{2} \|\chi_r (F(u)+F(\bar{u}))\|_{L^1 L^2(\Rm \times \Rm^5)} \\
  & \lesssim_1 \|G_e\|_{L^2([r,+\infty))} + |\tau_2(r)| + \|\chi_r u_e(|u|^{4/3}+|\bar{u}|^{4/3})\|_{L^1 L^2(\Rm \times \Rm^5)} \\
  & \lesssim_1 \|G_e\|_{L^2([r,+\infty))} + |\tau_2(r)| + M(r) \|\chi_r u\|_{Y(\Rm)}^{4/3}\\
  & \lesssim_1 \|G_e\|_{L^2([r,+\infty))} + |\tau_2(r)| + M(r) |\vec{\tau}(r)|^{4/3}. 
 \end{align*}
 When $\rho_0$ is sufficiently small, the last term in the right hand side can be absorbed by the left hand side. Thus 
 \[
  M(r) = \sup_{t\in \Rm} \|\vec{u}_e(t)\|_{\mathcal{H}(|t|+r)} + \|\chi_r u_e\|_{Y(\Rm)} \lesssim_1  \|G_e\|_{L^2([r,+\infty))} + |\tau_2(r)|. 
 \]
 We then recall that $u_e$ is a non-radiative solution and apply Lemma \ref{scatter profile of nonlinear solution} to deduce 
 \begin{align*}
  \|G_e\|_{L^2[r,+\infty)} & \lesssim_1 \|\chi_R (F(u)+F(\bar{u}))\|_{L^1 L^2((-\infty,0] \times \Rm^5)} \lesssim_1 M(r) |\vec{\tau}(r)|^{4/3}\\
  & \lesssim_1 \left(\|G_e\|_{L^2([r,+\infty))} + |\tau_2(r)|\right) |\vec{\tau}(r)|^{4/3}. 
 \end{align*}
 When $\rho_0$ is sufficiently small, we obtain $\|G_e\|_{L^2[r,+\infty)} \lesssim_1 |\tau_2(r)| |\vec{\tau}(r)|^{4/3}$. Therefore
 \begin{equation} \label{even estimate u}
  \|G_e\|_{L^2[r,+\infty)} + \sup_{t\in \Rm} \|\vec{u}_e(t)\|_{\mathcal{H}(|t|+r)} + \|\chi_r u_e\|_{Y(\Rm)} \lesssim_1 |\tau_2(r)|, \qquad r\geq R. 
 \end{equation}
 A similar argument shows that 
  \begin{equation} \label{odd estimate u}
  \|G_o\|_{L^2[r,+\infty)} + \sup_{t\in \Rm} \|\vec{u}_o(t)\|_{\mathcal{H}(|t|+r)} + \|\chi_r u_o\|_{Y(\Rm)} \lesssim_1 |\tau_1(r)|, \qquad r\geq R. 
 \end{equation}
 Now we are ready to give a point-wise estimate of $G_e$ and $G_o$. We recall that $u_e$ is a non-radiative solution, utilize the symmetry of radiation profiles and apply Lemma \ref{explicit formula for G plus by F} to deduce
 \[
  |G_e(s)| \leq \frac{1}{4}\int_0^\infty (s+t)^2 \left|F(u(s+t,t))+F(u(s+t,-t))\right|{\rm d} t + \frac{1}{4}\int_0^\infty \int_{t+s}^\infty r \left|F(u)+F(\bar{u})\right| {\rm d} r {\rm d} t.
 \]
We utilize the point-wise estimate of radial $\dot{H}^1$ functions and \eqref{pointwise u upgrade}, \eqref{even estimate u} to give an point-wise upper bound ($r\geq |t|+s$, $s\geq R$)
\begin{align*}
 \left|F(u(r,t)) + F(u(r,-t))\right| & \lesssim_1 |u_e| \left(|u|^{4/3} + |\bar{u}|^{4/3}\right) \\
 & \lesssim_1 r^{-3/2} \|u_e(t)\|_{\mathcal{H}(|t|+s)} \cdot r^{-2} (r-|t|)^{-2/3} R^{2/3} |\vec{\tau}(R)|^{4/3}\\
 & \lesssim_1 r^{-7/2} (r-|t|)^{-2/3} |\tau_2(s)| R^{2/3}  |\vec{\tau}(R)|^{4/3}.
\end{align*}
Inserting this upper bound, we obtain 
\begin{equation} \label{upper bound of Ge p}
 |G_e(s)| \lesssim_1 R^{2/3} s^{-7/6} |\tau_2(s)| |\vec{\tau}(R)|^{4/3}, \qquad s \geq R. 
\end{equation}
Similarly we have 
\begin{equation} \label{upper bound of Go p}
 |G_o(s)| \lesssim_1 R^{2/3} s^{-7/6} |\tau_1(s)| |\vec{\tau}(R)|^{4/3}, \qquad s \geq R.
\end{equation} 
A direct calculation shows that 
\begin{align*}
 \frac{\rm d}{{\rm d} r} |\tau_2(r)|^2 & = -\frac{3}{r} |\tau_2(r)|^2 + \frac{4\sqrt{3}}{r^{1/2}} G_e(r) \tau_2(r) = -\frac{3}{r} |\tau_2(r)|^2\left[1 + O\left(\frac{R^{2/3}}{r^{2/3}} |\vec{\tau}(R)|^{4/3}\right) \right]\leq 0; \\
 \frac{\rm d}{{\rm d} r} |\tau_1(r)|^2 & = -\frac{1}{r} |\tau_1(r)|^2 - \frac{4}{r^{1/2}} G_o(r) \tau_1(r) = -\frac{1}{r} |\tau_1(r)|^2\left[1 + O\left(\frac{R^{2/3}}{r^{2/3}}  |\vec{\tau}(R)|^{4/3}\right)\right] \leq 0.
\end{align*}
As a result, $\tau_j(r)$ satisfies either of the following
\begin{itemize}
 \item $\tau_j(r)$ never change its sign and $|\tau_j(r)|$ is strictly decreasing. 
 \item $\tau_j(r) \equiv 0$ for $r\geq 1$. 
\end{itemize}
Now we prove the equivalence relationship \eqref{equivalence asymptotic tau}. It suffices to consider the case $\tau_j(r) \neq 0$. By \eqref{upper bound of Go p} we have 
\begin{align*}
\left| \frac{\rm d}{{\rm d} r} \int_{-r}^r G(s) {\rm d} s \right| = 2\left| G_o(r) \right| \lesssim_1 |\vec{\tau}(R)|^{4/3} R^{2/3} r^{-5/3} \left|\int_{-r}^r G(s) {\rm d} s\right|. 
\end{align*}
It immediately follows that 
\[
 \left|\ln \left|\int_{-r_2}^{r_2} G(s) {\rm d} s\right| - \ln \left|\int_{-r_1}^{r_1} G(s) {\rm d} s\right|\right| \lesssim_1 |\vec{\tau}(R)|^{4/3} R^{2/3} \int_{r_1}^{r_2} r^{-5/3} {\rm d} r, \qquad r_2 > r_1 \geq R.
\]
This implies that $\ln \left|\int_{-r}^r G(s) {\rm d}s\right|$, thus $\int_{-r}^r G(s) {\rm d}s$ converges as $r\rightarrow +\infty$. In addition we may choose $r_2 = r$, $r_1 =R$ and obtain
\[
 \left|\ln \left|\int_{-r}^r G(s) {\rm d} s\right| - \ln \left|\int_{-R}^R G(s) {\rm d} s\right|\right| \lesssim_1 |\vec{\tau}(R)|^{4/3}\qquad \Rightarrow \qquad \left|\int_{-r}^r G(s) {\rm d} s\right| \simeq_1 \left|\int_{-R}^R G(s) {\rm d} s\right|.
\]
This finishes the proof. The case $j=2$ is similar, since we have
\begin{align*}
\left| \frac{\rm d}{{\rm d} r} \int_{-r}^r s G(s) {\rm d} s \right| = 2r \left| G_e(r) \right| \lesssim_1 |\vec{\tau}(R)|^{4/3} R^{2/3} r^{-5/3} \left|\int_{-r}^r s G(s) {\rm d} s\right|. 
\end{align*}
\end{proof}

\begin{remark} \label{decay of GoGe} 
 Let $u$ be as in Lemma \ref{detailed asymptotic behaviour}. The conclusion of Lemma \ref{detailed asymptotic behaviour} implies that ($r\geq R$)
 \begin{align*}
  &|\tau_1(r)| \simeq_1 (r/R)^{-1/2} |\tau_1(R)|& & |\tau_2(r)| \simeq_1 (r/R)^{-3/2} |\tau_2(R)|.
 \end{align*}
 Combining these with \eqref{upper bound of Ge p} and \eqref{upper bound of Go p}, we obtain ($s\geq R$)
 \begin{align*}
  & |G_e(s)| \lesssim_1 R^{13/6} s^{-8/3} |\vec{\tau}(R)|^{7/3};& &|G_o(s)| \lesssim_1 R^{7/6} s^{-5/3} |\vec{\tau}(R)|^{7/3}. 
 \end{align*}
\end{remark}

\begin{corollary} \label{decay of initial data non-ra}
 Under the same assumption of Lemma \ref{detailed asymptotic behaviour}, we have 
 \[
  \|(u_0,u_1)\|_{\mathcal{H}(r)} \simeq_1 (r/R)^{-1/2} |\tau_1(R)| + (r/R)^{-3/2} |\tau_2(R)|, \qquad \forall r\geq R. 
 \]
\end{corollary}
\begin{proof}
 By Lemma \ref{tail G} and the estimate given in the proof of Lemma \ref{detailed asymptotic behaviour}, we have 
 \begin{align*}
  \|(u_0,u_1)\|_{\mathcal{H}(r)} & \simeq_1 |\tau_1(r)|+|\tau_2(r)|+\|G\|_{L^2(\{s:|s|>r\})}\\
  & \simeq_1 |\tau_1(r)|+|\tau_2(r)| + \|G_o\|_{L^2([r,+\infty))} + \|G_e\|_{L^2([r,+\infty))} \\
  & \simeq_1 |\tau_1(r)|+|\tau_2(r)|\\
  & \simeq_1 (r/R)^{-1/2} |\tau_1(R)| + (r/R)^{-3/2} |\tau_2(R)|. 
 \end{align*}
\end{proof}

\begin{lemma} \label{existence of class}
 There exists a small constant $\rho_1>0$, such that if $\rho \in (0,\rho_1)$ and $\omega \in \mathbb{S}^1 \subset \Rm^2$, then there exists a radial non-radiative solution $u$ to (CP1) whose initial data comes with a radiation residue $\vec{\tau}(r)$ satisfying $\vec{\tau}(1) = \rho \omega$ and $|\vec{\tau}(r)|<\rho$ for all $r>1$. 
\end{lemma}

\begin{proof}
 Let us consider the complete metric space (with distance given by the $L^2$ norm of difference)
 \[
  X = \left\{G\in L^2(\{s: |s|>1\}): \|G\|_{L^2(\{s: |s|>1\})} \leq c \rho^{7/3}\right\}. 
 \]
 Here $c$ is a large constant to be determined later. We choose $\rho_1 = c^{-3/4} \ll 1$. Given any $G \in X$, we extend its domain to the whole real line by setting 
 \begin{align*}
  \left(-\int_{-1}^1 G(s) {\rm d} s, \sqrt{3} \int_{-1}^1 s G(s) {\rm d} s\right) = \rho \omega. 
 \end{align*}
 Let $(u_0,u_1)$ be initial data with radiation profile $G$ (in the negative time direction). Then our assumption above guarantees that 
 \[
  \|(u_0,u_1)\|_{\mathcal{H}(1)} \lesssim_1 c \rho^{7/3} + \rho \lesssim_1 \rho \ll 1.
 \] 
 Small data theory then gives an exterior solution $u$ defined in the whole exterior region $\Omega_1$ with 
 \[
  \|\chi_1 u\|_{Y((\Rm)} \lesssim_1 \rho.
 \]
 Next we let $G_\pm$ be the nonlinear radiation profile of the exterior solution to the wave equation 
 \[
  \left\{\begin{array}{l} \partial_t^2 w - \Delta w = F(u), \qquad (x,t) \in \Omega_1; \\ \vec{w}(0) = 0. \end{array}\right.
 \]
 It immediately follows from Lemma \ref{scatter profile of nonlinear solution} that 
 \[
  \|G_\pm\|_{L^2([1,+\infty))} \lesssim_1 \rho^{7/3}. 
 \]
 We define a map $\mathbf{T} : X \rightarrow X$ by 
 \[
  (\mathbf{T} G)(s) = \left\{\begin{array}{ll} - G_-(s), & s>1; \\ -G_+(-s), & s<-1.\end{array}\right.
 \]
 A large value of $c$ guarantees that $\mathbf{T} X \subset X$. It is not difficult to see that $u$ is a non-radiative solution in $\Omega_1$ if and only if $\mathbf{T} G = G$. Next we show that a fixed-point of $\mathbf{T}$ exists, which follows from the contraction map theorem. Indeed, given two $G,G_\ast\in X$, we recall Lemma \ref{scatter profile of nonlinear solution} and apply the perturbation theory of small solution to deduce
 \begin{align*}
  \|\mathbf{T} G -\mathbf{T} G_\ast\|_{L^2(\{s:|s|>1\})} & \lesssim_1 \|\chi_1(F(u)-F(u_\ast))\|_{L^1 L^2(\Rm \times \Rm^5)} \\
  & \lesssim_1 \rho^{4/3} \|\chi_1(u-u_\ast)\|_{Y(\Rm)}\\
  & \lesssim_1 \rho^{4/3} \|\vec{u}(0)-\vec{u}_\ast(0)\|_{\mathcal{H}(1)}\\
  & \lesssim_1 \rho^{4/3} \|G-G_\ast\|_{L^2(\{s: |s|>1\})}. 
 \end{align*}
 Here $u_\ast$ is the corresponding exterior solution for $G_\ast$. When $c$ is sufficiently large, or equivalently speaking, $\rho_1$ is sufficiently small, $\mathbf{T}$ becomes a contraction map. The way we constructed $G$ implies that the corresponding radiation residue $\vec{\tau}(1)$ of the non-radiative solution $u$ is exactly $\rho \omega$. Finally we show that $|\vec{\tau}(r)|<\rho$ as long as $\rho<\rho_1$ is sufficiently small. Indeed, we have 
 \[
  |\vec{\tau}(r)| \lesssim_1 \|(u_0,u_1)\|_{\mathcal{H}(r)} \lesssim_1 \rho, \qquad \forall r\geq 1.
 \]
 Thus we may apply Lemma \ref{monotonicity of tau} and finish the proof if $\rho_1$ is sufficiently small. 
 \end{proof}

\subsection{Two classifications of non-radiative solutions}

\paragraph{Characteristic angle} If we fix a small positive constant $\rho$, then Lemma \ref{existence of class} gives a radial non-radiative solution $N_\omega$ for each $\omega \in \mathbb S^1$, such that its corresponding radiation residue $\vec{\tau}$ satisfies 
\begin{align*}
 &\vec{\tau}(1) = \rho \omega;& &|\vec{\tau}(r)| < \rho, \quad \forall r>1. 
\end{align*}
The scaling invariance of radiation residue implies that any two different non-radiative solutions $N_{\omega_1}$ and $N_{\omega_2}$ can never coincide with each other even if a dilation can be applied. In addition, any non-trivial non-radiative solution must be one of these $N_\omega$'s, up to a dilation, thanks to Lemma \ref{small tau zero} and Lemma \ref{identification tau}. Therefore we may define $\mathcal{N}_\omega$ (called non-radiative class) to be the set containing all dilations of $N_\omega$ and write the set of all non-radiative solutions in the form 
\[
 \mathcal{N} = \{\mathbf{0}\} \cup \left(\bigcup_{\omega\in \mathbb S^1} \mathcal{N}_\omega\right). 
\]
By applying Lemma \ref{lemma connection} and Lemma \ref{lemma connection 2}, we may compare two near-by non-radiative solutions and obtain the following continuity
\begin{lemma}
 Let $N_{\omega_0}$ be a non-radiative solution defined above and $R> 0$ be a radius such that $\|\chi_R N_{\omega_0}\|_{Y(\Rm)} < +\infty$. Then 
 \[
  \lim_{\omega \rightarrow \omega_0} \left(\left\|\vec{N}_{\omega}(0)- \vec{N}_{\omega_0}(0)\right\|_{\mathcal{H}(R)}+\|\chi_R (N_{\omega}-N_{\omega_0})\|_{Y(\Rm)}\right) = 0. 
 \]
\end{lemma}
For convenience we also use the notations $N_\theta$ and $\mathcal{N}_\theta$ for the corresponding non-radiative solution and class for $\omega = (\cos \theta, \sin \theta) \in \mathbb{S}^1$. We call $\theta$ (or $\omega$) the characteristic angle of a non-radiative solution. Please note that the non-radiative class $\mathcal{N}_\theta$ depends on the choice of $\rho$, except for four special cases, i.e. $\theta = 0, \pi/2, \pi, -\pi/2$. In the rest of this work we shall fix a sufficiently small $\rho$ and utilize the non-radiative classes associated to it. More details about the four special classes $\theta = 0, \pi/2, \pi, -\pi/2$ will be discussed later in this section. 

\paragraph{Asymptotic numbers} Given any non-radiative solution $u$, the limits 
\begin{align*}
 &\alpha_1 \doteq -\lim_{r\rightarrow +\infty} \int_{-r}^r G(s) {\rm d} s = -\int_{-\infty}^\infty G_o(s) {\rm d} s;& & \alpha_2 \doteq \lim_{r\rightarrow +\infty} \int_{-r}^r s G(s) {\rm d} s = \int_{-\infty}^\infty s G_e(s) {\rm d} s
\end{align*}
are both well-defined, according to Lemma \ref{detailed asymptotic behaviour}. Please note that the decay of $G_o$ and $G_e$ given in Remark \ref{decay of GoGe} guarantees that the improper integrals of them converge absolutely. We call $\alpha_1,\alpha_2$ the first and second asymptotic number, respectively. Please note that the zero solution also comes with asymptotic numbers $(0,0)$. All other non-radiative solutions come with at least one nonzero asymptotic number. The asymptotic numbers can also be determined by observing the asymptotic behaviour of initial data, as the name indicates. Let us use the notations in the proof of Lemma \ref{detailed asymptotic behaviour}. Consider the initial data $(\alpha_2 |x|^{-3},0)$ in the exterior region $\{x: |x|>r\}$, whose radiation profile $G_\ast$ satisfies 
\begin{align*}
 &G_\ast (s) = 0, |s|>r; & &\int_{-r}^r G_\ast (s) {\rm d} s = 0;& &\int_{-r}^r  s G_\ast (s) {\rm d} s = \alpha_2. 
\end{align*}
By Lemma \ref{tail G} and Remark \ref{decay of GoGe}, we have for large $r$ that
\begin{align*}
 \left\|(u_0,0) - (\alpha_2 |x|^{-3},0)\right\|_{\mathcal{H}(r)} & \lesssim_1 \|G_e\|_{L^2(\{s: |s|>r\})} + \frac{1}{r^{3/2}}  \left|\alpha_2 - \int_{-r}^r s G_e(s) {\rm d} s\right|\\
 & \lesssim_1 \|G_e\|_{L^2([r,+\infty))} + \frac{1}{r^{3/2}}  \left|\int_r^\infty s G_e(s) {\rm d} s\right| \lesssim r^{-13/6}. 
\end{align*}
The implicit constant above is independent of sufficiently large radius $r$, but may depend on the specific non-radiative solution. In other words, we have 
\[
 \left\|u_0 - \alpha_2 |x|^{-3}\right\|_{\dot{H}^1(\{x:|x|>r\})} \lesssim  r^{-13/6}. 
\]
Similarly we have
\begin{align*}
 \left\|u_1 - \alpha_1 |x|^{-3}\right\|_{L^2(\{x: |x|>r\})} & \lesssim_1 \|G_o\|_{L^2([r,+\infty))} + \frac{1}{r^{1/2}} \left|\int_r^\infty G_o(s) {\rm d} s\right| \lesssim r^{-7/6}.
\end{align*}
Since $\||x|^{-3}\|_{\dot{H}^1(\{x:|x|>r\})} \simeq r^{-3/2}$ and $\||x|^{-3}\|_{L^2(\{x:|x|>r\})} \simeq r^{-1/2}$, the asymptotic behaviour given above uniquely determines the values of $\alpha_1$ and $\alpha_2$. Next we prove that the asymptotic numbers uniquely determine a non-radiative solution. 

\begin{lemma}
 Let $u$ and $u^\ast$ be two radial radiative solutions of (CP1) with the same first and second asymptotic numbers. Then $u\equiv u^\ast$.
 \end{lemma}
\begin{proof}
 It suffices to show $u=u^\ast$ in the exterior region $\Omega_R$ for a very large number $R$. Let $G$, $G^\ast$ be the corresponding radiation profiles of their initial data. We use the same notations $\tau_j$, $\tau_j^\ast$, $G_e$, $G_e^\ast$, $G_o$, $G_o^\ast$ as given in the proof of Lemma \ref{detailed asymptotic behaviour}. 
We choose a sufficiently large number $R\gg 1$ such that
\begin{align*}
 \|G_e\|_{L^2([r,+\infty))}, \|G_e^\ast\|_{L^2([r,+\infty))} & \leq r^{-2}, & & r\geq R;\\
 \|G_o\|_{L^2([r,+\infty))}, \|G_o^\ast\|_{L^2([r,+\infty))} & \leq r^{-1}, & & r\geq R;\\
 \|\chi_r u\|_{Y(\Rm)}, \|\chi_r u^\ast\|_{Y(\Rm)} & \leq  c r^{-3/8}. & & r\geq R.
\end{align*} 
Here $c$ is a small absolute constant to be determined later. The first two lines hold for sufficiently large $R$ because higher decay rates ($13/6$ and $7/6$) can be verified by Remark \ref{decay of GoGe}. The last line holds for sufficiently large $R$ by Corollary \ref{decay of initial data non-ra} and the small data theory. We will apply an induction argument. We assume that the following inequalities hold for some positive number $\kappa \geq 1/2$. 
\begin{align*}
 \|G_e - G_e^\ast\|_{L^2([r,+\infty))} & \leq 2 r^{-\max\{1/2+\kappa,2\}}, & & r\geq R; \\
 \|G_o - G_o^\ast\|_{L^2([r,+\infty))} & \leq 2 r^{-1/2-\kappa}, & & r\geq R.
\end{align*} 
Clearly this holds for $\kappa =1/2$. By our assumption on the asymptotic numbers, we have 
\begin{align*}
 |\tau_1(r) - \tau_1^\ast (r)| & \lesssim_1 \frac{1}{r^{1/2}} \int_r^\infty |G_o(s)-G_o^\ast (s)| {\rm d} s \lesssim_1 r^{-1/2 -\kappa}, \qquad r\geq R; \\
 |\tau_2(r) - \tau_2^\ast (r)| & \lesssim_1 \frac{1}{r^{3/2}} \int_r^\infty |sG_e(s)-sG_e^\ast (s)| {\rm d} s \lesssim_1 r^{-\max\{1/2+\kappa,2\}}, \qquad r\geq R
\end{align*}
It follows that $\|\vec{u}(0)-\vec{u}^\ast(0)\|_{\mathcal{H}(r)} \lesssim_1 r^{-1/2-\kappa}$. By the Lipchitz continuity of small data solution, we obtain 
\[
 \|\chi_r (u-u^\ast)\|_{Y(\Rm)} \lesssim_1 r^{-1/2-\kappa}. 
\] 
Applying Lemma \ref{scatter profile of nonlinear solution}, we obtain 
\begin{align*}
 \|G-G^\ast\|_{L^2(\{s: |s|>r\})} & \lesssim_1 \|\chi_r(F(u)-F(u^\ast))\|_{L^1 L^2(\Rm \times \Rm^5)}\\
 & \lesssim_1 \|\chi_r(u-u^\ast)\|_{Y(\Rm)} \left(\|\chi_r u\|_{Y(\Rm)}^{4/3} + \|\chi_r u^\ast\|_{Y(\Rm)}^{4/3}\right)\\
 & \lesssim_1 c^{4/3} r^{-1-\kappa}. 
\end{align*}
As a result, we have
\[
  \|G_e - G_e^\ast\|_{L^2([r,+\infty))} +  \|G_o - G_o^\ast\|_{L^2([r,+\infty))} \lesssim_1 c^{4/3} r^{-1-\kappa}. 
\]
When $c$ is a sufficiently small, we immediately obtain
\begin{align*}
 \|G_e - G_e^\ast\|_{L^2([r,+\infty))} & \leq 2 r^{-\max\{1+\kappa,2\}}, & & r\geq R; \\
 \|G_o - G_o^\ast\|_{L^2([r,+\infty))} & \leq 2 r^{-1-\kappa}, & & r\geq R.
\end{align*} 
The argument above makes the decay rate $\kappa$ increase by $1/2$. Repeating this process, we obtain that the inequality $\|\chi_R (u-u^\ast) \|_{Y(\Rm)} \lesssim_1 R^{-1/2-\kappa}$ holds for all sufficiently large $\kappa$. This implies that the norm must be zero thus $u=u^\ast$ in the exterior region $\Omega_R$. 
\end{proof}

\paragraph{Special non-radiative classes} Now let us consider two special non-radiative solutions, i.e. the ground states $\pm W(x)$. By time symmetry, the radiation profile $G_-$ of $(W,0)$ is an odd function of $s$, thus the characteristic profile $\tau_1 (r) = 0$ for all $r>0$. In addition, we recall the explicit formula 
\[
 \tau_2 (r) = \frac{\sqrt{3}}{r^{3/2}} \int_{-r}^r s G_-(s) {\rm d} s = \sqrt{3} r^{3/2}  W(r) > 0, \qquad r > 0;
\]
and obtain that $W \in \mathcal{N}_{\pi/2}$. Similarly $-W(x) \in \mathcal{N}_{-\pi/2}$. Next we consider another special case $\mathcal{N}_0$. This one can not be given by an explicit formula but we may analyze its property. Since $\vec{\tau}(1) = (\rho,0)$ for $N_0$, the second radiation residue has to remain zero for all $r > 1$. It follows that the initial value $(u_0,u_1)$ satisfy $u_0(r) = 0$ for all $r\geq 1$. Thus $N_0$ must be an odd function of $t$ in the exterior $\Omega_1$, i.e. $N_0(r, -t) = - N_0(r,t)$ for $r>|t|+1$. Because both $N_0(r,-t)$ and $-N_0(r,t)$ are non-radiative solutions with the same domain of $N_0$, they must be the same in the whole domain. In other words, $N_0$ must be an odd function of time $t$ in its whole domain. In addition, the first asymptotic number of $N_0$ must be a positive number, according to Lemma \ref{detailed asymptotic behaviour}. Conversely, if $u$ is an odd non-radiative solution with a positive first asymptotic number, then its second radiation residue must be zero and its first radiation residue must be positive near infinity, which guarantees that $u \in \mathcal{N}_0$. It is not difficult to see $\mathcal{N}_\pi = -\mathcal{N}_0$. Please note that the non-radiative classes $\mathcal{N}_{\pi/2}, \mathcal{N}_{-\pi/2}, \mathcal{N}_0, \mathcal{N}_\pi$ are independent of the choice of small positive constant $\rho$. More behaviour of the non-radiative class $\mathcal{N}_0$ will be discussed in later sections, which plays an important role in the radiation theory of global solutions to (CP1).

\subsection{Time evolution of non-radiative solutions}

\paragraph{Time translation relationship} Let us consider a non-radiative solution $u$ defined in $\Omega_R$ with the asymptotic numbers $\alpha_j$. By the asymptotic behaviour its initial data $(u_0,u_1)$ satisfy
\begin{align*}
 &\left\|u_0 -\alpha_2 |x|^{-1}\right\|_{\dot{H}^1(\{x: |x|>r\})} \lesssim r^{-13/6};& &\left\|u_1 -\alpha_1 |x|^{-1}\right\|_{L^2(\{x: |x|>r\})} \lesssim r^{-7/6}.
\end{align*}
Clearly $v_{t_0} = u(\cdot, \cdot+t_0)$ is also a non-radiative solution defined in the region $\Omega_{R+|t_0|}$. Now let us consider the asymptotic numbers of $v_{t_0}$. By perturbation theory, we may compare $u$ with the approximated solution $\alpha_2 |x|^{-3} + \alpha_1 t |x|^{-3}$ to deduce  
\[
 \sup_{t\in \Rm} \left\|\vec{u}(t) - \left((\alpha_2+t\alpha_1)|x|^{-3},\alpha_1 |x|^{-3}\right)\right\|_{\mathcal{H}(r+|t|)} \lesssim r^{-7/6}, \qquad r\gg 1.
\]
By integrating in $t$, we also have for a fixed $t_0$ that
\begin{align*}
& \left\|u(\cdot,t_0) - (\alpha_2 +t_0 \alpha_1)|x|^{-3}\right\|_{L^2(x: |x|>|t_0|+r)}\\
  & \leq \left\|u(\cdot,0) - \alpha_2 |x|^{-3}\right\|_{L^2(x: |x|>|t_0|+r)} + |t_0| \sup_{t\in [0,t_0]} \left\|u_t(\cdot,t)-\alpha_1 |x|^{-3}\right\|_{L^2(\{x: |x|>|t_0|+r\})} \\
 & \lesssim \sum_{k=0}^\infty 2^k (|t_0|+r) \left\|\frac{u(\cdot,0) - \alpha_2 |x|^{-3}}{|x|}\right\|_{L^2(\{x: |x|>2^k(|t_0|+r)\})} + r^{-7/6} \\
 & \lesssim \sum_{k=0}^\infty 2^k (|t_0|+r) \left\|u(\cdot,0) - \alpha_2 |x|^{-3}\right\|_{\dot{H}^1(\{x: |x|>2^k(|t_0|+r)\})} + r^{-7/6}\\
 & \lesssim \sum_{k=0}^\infty \left(2^k (|t_0|+r)\right)^{1-13/6} + r^{-7/6}\\
 & \lesssim r^{-7/6}. 
\end{align*}
Here the implicit constant is independent of $r\gg 1$. On the other hand, if $(\alpha_1(t_0), \alpha_2(t_0))$ are the asymptotic numbers of $v_{t_0}$, then its initial data $\vec{v}_{t_0}(0) = \vec{u}(t_0)$ satisfy
\begin{align*}
 \left\|\vec{u}(t_0) - (\alpha_2(t_0) |x|^{-3}, \alpha_1(t_0) |x|^{-3})\right\|_{L^2(\{x: |x|>r\})} \lesssim r^{-7/6}. 
\end{align*}
Again we apply dyadic decomposition and Hardy's inequality here to deduce the $L^2$ estimate of $u(\cdot,t_0)$ by the $\dot{H}^1$ estimate. Since $\||x|^{-3}\|_{L^2(x:|x|>r)}\simeq r^{-1/2}$, we may compare the asymptotic behaviours given above and obtain  
\begin{align*}
 &\alpha_1(t_0) = \alpha_1;& &\alpha_2(t_0) = \alpha_2 + \alpha_1 t_0. 
\end{align*}
\paragraph{Completeness of asymptotic numbers} Since we have already shown that the ground states are the non-radiative solutions with zero first asymptotic number and nonzero second asymptotic number; and the odd non-radiative solutions in $\mathcal{N}_0$ come with nonzero first asymptotic number and zero second asymptotic number. The time-translation relationship actually means that there exists a non-radiative solution with any given pair of asymptotic numbers $(\alpha_1, \alpha_2) \in \Rm^2$. 

\paragraph{Time evolution of non-radiative classes} If $u$ is a ground state, then clearly $v_{t_0}$ defined above are the same ground state for all time $t$. Now let us consider other situations. By a time translation/sign change and dilation it suffices to consider the case $u \in \mathcal{N}_0$ with a first characteristic number $\alpha_1 = 1$. Thus the characteristic numbers are given by $(\alpha_1(t), \alpha_2(t)) = (1,t)$. The asymptotic numbers of time-translated solutions also imply that $\vec{u}(r,t+t_1)$ and $\vec{u}(r,t+t_2)$ can never be dilation of each other unless $t_1 = t_2$. In other words, $u(\cdot,\cdot+t)$ falls in different non-radiative classes for different times $t$. In addition, the positive first asymptotic number implies that $u(\cdot,\cdot+t)$ always falls in a non-radiative class $\mathcal{N}_{\theta}$ with $\theta \in (-\pi/2,\pi/2)$. Now let us estimate the characteristic angles $\theta(t')$ and the scale $R(t')$ of $v_{t'}$. Here the scale is the radius $R$ such that $\vec{\tau}(R) = \rho (\cos \theta, \sin \theta)$. An application of Lemma \ref{detailed asymptotic behaviour} yields
\begin{align*}
 R(t)^{1/2} \rho \cos \theta(t) & \simeq_1 \alpha_1(t) = 1; \\
 R(t)^{3/2} \rho \sin \theta(t) & \simeq_1 \alpha_2(t) = t. 
\end{align*}
When $t>0$ is sufficiently large, we obtain that $\theta$ is roughly $\pi/2$ with
\begin{align*}
 &R(t) \simeq t^{2/3};& &\frac{\pi}{2} -\theta \simeq t^{-1/3}. 
\end{align*}
Similarly when $t<0$ is a large negative number, $\theta$ is roughly $-\pi/2$ with
\begin{align*}
 &R(t) \simeq t^{2/3};& &\frac{\pi}{2} + \theta \simeq |t|^{-1/3}. 
\end{align*}
Next we fix a time $t_0$, recall Lemma \ref{monotonicity of tau} and deduce that there exists a radius $R<R(t_0)$, such that the norm of radiation residue $\vec{\tau}(r)$ of $\vec{u}(t_0)$ is strictly decreasing in $[R,+\infty)$. By the uniform boundedness 
\[
 \sup_{r\geq R} |\vec{\tau}(r)| \lesssim_1 \|\vec{u}\|_{\mathcal{H}(R)},
\]
we deduce that both $R(t)$ and $\theta(t) \in \mathbb{S}^1$ are continuous functions of $t$. We summarize all of the above and conclude that the map from the time $t$ to the non-radiative classes is a homeomorphism from $\Rm$ to the right half of the circle $\mathbb{S}^1$. Or equivalently speaking, the map from $t\in \Rm$ to $\theta(t)$ is a strictly increasing homeomorphism form $\Rm$ to $(-\pi/2,\pi/2)$. 

\begin{remark}
 The argument above shows that the non-radiative class always approaches the classes of ground states as the time tends to infinity, which indicates that the soliton resolution contains ground states only. 
\end{remark}

\begin{remark} \label{ut point-wise estimate}
 Let $u$ be the odd non-radiative solution above and $t\in (0,1)$ be a small time. We use the notation $G$ and $\vec{\tau}$ for the radiation profile and residue of $\vec{u}(t)$. Please note that the scale $R(t)$ is uniformly bounded for all $t\in [0,1]$. For a large radius $\displaystyle r > \sup_{t\in (0,1)} R(t)$ we have
 \begin{align*}
  u(r,t) & = \frac{1}{r^3} \int_{-r}^r s G(s) {\rm d} s = \frac{1}{r^3} \left[t - 2\int_r^\infty sG_e(s) {\rm d} s\right] \\
  & = \frac{t}{r^3} \left[1+O(\rho^{4/3} R(t)^{2/3} r^{-2/3})\right]. 
 \end{align*}
 Here we utilize notation $G_e$ and the estimate \eqref{upper bound of Ge p} given in the proof of Lemma \ref{detailed asymptotic behaviour}, as well as the conclusion of Lemma \ref{detailed asymptotic behaviour}; the implicit constant for $O(\rho^{4/3}R(t)^{2/3} r^{-2/3})$ is an absolute constant. 
 A point-wise estimate $u_t(x,0) = |x|^{-3}(1+O(|x|^{-2/3}))$ immediately follows for large $x$. 
\end{remark}

\begin{corollary} \label{asymptotic of non-ra 23}
 Let $v$ be a radial non-radiative solution to (CP1) with a maximal domain $\Omega_{R_0}$. Then $v$ may extend to a $\mathcal{C}^1$ function in a small neighbourhood of any given point $(x,t)$ on the boundary point of $\Omega_{R_0}$, as long as $|t|$ is sufficiently large. In addition, we have
 \[
  \|\vec{v}(t)\|_{\mathcal{H}(|t|+R_0)} \lesssim |t|^{-1/6}. 
 \]
\end{corollary}
\begin{proof}
 This is clear true for any ground state or the zero solution. By a dilation/sign symmetry, we may let $v = v_{t_0}$ be the non-radiative solution defined above for some time $t_0$. Given a time $t_1$, we have $\vec{v}(t_1) = \vec{u}(t_0+t_1)$. Recalling the scale function $R(t)$ introduce above, we may extend the domain of $v$ (without loss of non-radiative property) to 
 \[
  \{(x,t): |x|>R(t_0+t_1) + |t-t_1|\},
 \]
 which contains all the boundary points $\{(x,t_1): |x|=|t_1|+R_0\}\subset \partial \Omega_{R_0}$, as long as $|t_1|$ is sufficiently large, by the asymptotic behaviour $R(t) \simeq t^{2/3}$ for large time.  The $\mathcal{C}^1$ continuity of non-radiative solutions then finishes the proof of continuous extension. Finally we recall Corollary \ref{decay of initial data non-ra} and deduce 
 \[
  \|\vec{v}(t_1)\|_{\mathcal{H}(|t_1|+R_0)} \lesssim \left(\frac{|t_1|+R_0}{R(t_0+t_1)}\right)^{-1/2} \lesssim |t_1|^{-1/6}, \qquad |t_1| \gg 1. 
 \]
 The implicit constant in the inequality above does not depend on $t_1$. 
\end{proof}

\subsection{Global existence and universal cover}

\begin{figure}[h]
 \centering
 \includegraphics[scale=0.85]{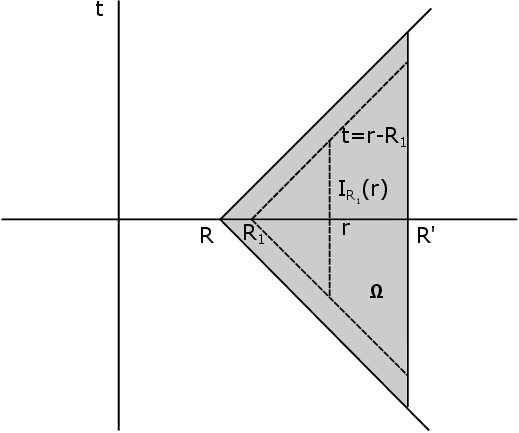}
 \caption{Integral path and the region $\Omega$} \label{figure global}
\end{figure}

In this subsection we show that all the radial non-radiative solutions can are defined in the region $\Omega_0$, although some of them may be blow-up solutions. Let $u$ be a radial non-radiative solution with a maximal region $\Omega_R$. If $u$ is s stationary solution, then $u$ is either zero or a ground state, thus a scattering solution defined in $\Omega_0$. Thus it suffice to consider non-stationary solutions. We rewrite the equation (CP1) in the form of 
\[
 u_{rr} - u_{tt} + |u|^{4/3} u = -\frac{4}{r} u_r. 
\] 
Given any $R_1>R$, let us define 
\[
 I_{R_1}(r) = \int_{-r+R_1}^{r-R_1} \left(\frac{1}{2} |u_r(r,t)|^2 + \frac{1}{2} |u_t(r,t)|^2 + \frac{3}{10} |u(r,t)|^{10/3} \right){\rm d}t, \qquad r>R_1. 
\]
It follows that ($r>R_1>R$)
\begin{align*}
 \frac{\rm d}{{\rm d}r} I_{R_1}& = -\frac{4}{r} \int_{-r+R_1}^{r-R_1} |u_r(r,t)|^2 {\rm d} t + \left(\frac{1}{2}|(u_t+u_r)(r,r-R_1)|^2 + \frac{3}{10} |u(r,r-R_1)|^{10/3}\right)\\
  &\qquad + \left(\frac{1}{2}|(u_t-u_r)(r,R_1-r)|^2 + \frac{3}{10} |u(r,R_1-r)|^{10/3}\right) \geq -\frac{8}{r} I_{R_1}(r). 
\end{align*}
Here we apply smooth approximation techniques and utilize the fact $u \in \mathcal{C}^1$. It follows from the backward Gronwall's inequality that 
\begin{equation} \label{Gronwall ineq}
 I_{R_1} (r) \leq \left(r/R'\right)^{-8} I_{R_1} (R'), \qquad R_1<r<R'. 
\end{equation}
According to Corollary \ref{asymptotic of non-ra 23} and the small data theory, we may fix a large radius $R'\gg R+1$, such that $I_{R_1}(R')$ is uniformly bounded with respect to $R_1\in (R,R')$ and that 
\[
  \|\chi_R u\|_{Y(\{t: |t|>R'-R-1\})} < +\infty.  
\]
If $R$ were not zero, then \eqref{Gronwall ineq} implies that $I_{R_1}(r)$ would be uniformly bounded for all $(R_1,r)$ with $R<R_1<r<R'$. By the expression of derivative given above, we also conclude that the integrals 
\begin{align*}
 \int_{R_1}^{R'} \left(\frac{1}{2}|(u_t+u_r)(r,r-R_1)|^2 + \frac{3}{10} |u(r,r-R_1)|^{10/3}\right) {\rm d} r; \\
 \int_{R_1}^{R'} \left(\frac{1}{2}|(u_t-u_r)(r,R_1-r)|^2 + \frac{3}{10} |u(r,R_1-r)|^{10/3}\right) {\rm d} r
\end{align*}
are also uniformly bounded for all $R_1\in (R,R')$. Please note that $(u_t+u_r)(r,r-R_1)$ is the derivative of $u(r,r-R_1)$ with respect to $r$ and that $(-u_t+u_r)(r,R_1-r)$ is the derivative of $u(r,R_1-r)$ with respect to $r$. It immediately follows that $u(r,t)$ are uniformly bounded in the region (please see figure \ref{figure global})
\[
 \Omega = \left\{(r,t): R-r<t<r-R, R<r<R'\right\}. 
\]
This implies that 
\[
 \|\chi_R u\|_{Y(\Rm)} \leq \|\chi_{R+1} u\|_{Y(\Rm)} + \|\chi_\Omega u\|_{Y(\Rm)} + \|\chi_R u\|_{Y(\{t: |t|>R'-R+1\})} < +\infty. 
\]
This is a contradiction. Thus we must have $R = 0$. 

\begin{remark}
 A similar argument (as well as Corollary \ref{asymptotic of non-ra 23} and radial Sobolev embedding) shows that any radial non-radiative solution $u$ must satisfies 
\[
 \sup_{|x|>\max\{|t|,r\}} |u(x,t)| < +\infty, \qquad \forall r>0. 
\]
%Finite speed of propagation and the Strichartz estimates then give 
%\[
% \sup_{t\in \Rm} \|\vec{u}(t)\|_{\mathcal{H}(\max\{|t|,r\})} < +\infty, \qquad \forall r>0. 
%\]
\end{remark}

\paragraph{Universal cover} Let us fix a non-radiative solution $u \in \mathcal{N}_0$. Without loss of generality, we may apply a dilation and assume that $u$ and its radiation residue $\vec{\tau}(r)$ satisfy
\begin{align*}
 %&|\vec{\tau}(r)| < \rho, \qquad r\geq 1;\\ 
 \sup_{t\in \Rm} \left(\|\vec{u}(t)\|_{\mathcal{H}(1+|t|)} + \|u(\cdot,t)\|_{L^{10/3}(\{x: |x|>1+|t|\})} + \||x|^{-1} u(x,t)\|_{L^2(\{x: |x|>1+|t|\})} \right) \ll 1. 
\end{align*}
The Sobolev embedding and the continuity then implies 
\begin{equation} \label{lower energy U}
 \inf_{t\in [-1/10,1/10]} \int_{|x|>3} \left(\frac{1}{2}|\nabla U(x,t)|^2 + \frac{1}{2} |U_t(x,t)|^2 - \frac{3}{10} |U(x,t)|^\frac{10}{3}\right) {\rm d} x > 0. 
\end{equation}
Here $\rho$ is the small constant in the classification of non-radiative solutions. 
%Lemma \ref{detailed asymptotic behaviour}. 
Given $t_0 \in \Rm$, we may extend the domain of the non-radiative solution $u(x,t+t_0)$ to $\Omega_0$ by the global existence result given above. In other words, we may extend the domain of $u$ to $\{(x,t): |x|>|t-t_0|\}$. This process can be done for each $t_0 \in \Rm$ and defines a solution $U(x,t)$ in $(\Rm^5 \setminus \{0\})\times \Rm$. We say $U$ is a solution to (CP1) in the following sense: Given any $t_0\in \Rm^+\times \Rm$, $U(x,t+t_0)$ is a non-radiative solution to (CP1) with a maximal domain $\Omega_0$. Please note that although the values of $U(x,t)$ may be defined for multiple times, finite speed propagation guarantees that these definitions are consistent with each other. Given any non-stationary non-radiative exterior solution $v$ to (CP1), there exists a sign $\zeta\in \{+1,-1\}$, a scale $\lambda>0$ and a time $t_0 \in \Rm$, such that 
\[
 v(x,t) = \zeta \frac{1}{\lambda^{3/2}} U \left(\frac{x}{\lambda}, \frac{t+t_0}{\lambda}\right), \qquad |x|>|t|.
\]
In addition, given a small neighbourhood of $\{\pm \pi/2\}$, there exists a constant $T_0$, such that if the characteristic angle of a non-radiative solution $v$ is not contained in the neighbourhood mentioned above, then the corresponding $(\zeta, \lambda, t_0)$ satisfies $|t_0|/\lambda \leq T_0$. 
\section{Odd non-radiative solutions}

In this section we investigate the global behaviour of the non-radiative solutions in the non-radiative class $\mathcal{N}_0$. 

\paragraph{Self-similar solutions} We first make a brief review on the self-similar solutions of (CP1). Let us consider the exterior solutions in the form of 
\[
 u(x,t) = |x|^{-3/2} \varphi(t/|x|), \qquad (x,t) \in \Omega_0.
\]
A direct calculation shows that $u$ is a solution to (CP1) if and only if the function $\varphi$ solves the following ordinary differential equation
\[
 (1-y^2) \varphi''(y) - y \varphi'(y) + \frac{9}{4} \varphi - |\varphi|^{4/3} \varphi = 0, \qquad y\in (-1,1). 
\]
We are particularly interested in the odd solutions $\varphi$ with initial data 
\begin{align*}
 &\varphi(0) = 0;& & \varphi'(0) = \nu > 0. 
\end{align*}
The following quantity is conserved for all $y$ in the interval of existence
\[
 \frac{1}{2}(1-y^2)|\varphi'(y)|^2 + \frac{9}{8} |\varphi(y)|^2 - \frac{3}{10} |\varphi(y)|^{10/3} = \frac{1}{2} \nu^2. 
\]
A careful analysis and the classic theory of ordinary differential equations indicates that the fate of solution falls in the following two categories 
\begin{itemize} 
 \item[(i)] The solution blows up at $y_+ \in (0,1]$ with $|\varphi|\rightarrow +\infty$;  
 \item[(ii)] The solution is defined for all $t\in (0,1)$ and both $(1-y^2)|\varphi'(y)|^2$ and $\varphi(y)$ converge as $y\rightarrow 1^-$;
\end{itemize} 
In addition, the solution is strictly increasing for large values of $\nu > 1.8$. We conduct a numerical simulation for small values $\nu \in (0,2]$. It turns that the solution $\varphi_\nu$ satisfies 
\begin{itemize}
 \item[(a)] The solution $\varphi_\nu$ is in category (ii); 
 \item[(b)] The solution is aways positive and satisfies that $\varphi_\nu(y)< \nu y$ for $y\in (0,1/2]$.
\end{itemize}
Please note that the corresponding solution $u = |x|^{-3/2} \varphi_\nu(t/|x|)$ is always a $C^2$ (classic) solution to (CP1) in the region $\Omega_0$ but $u$ is NOT contained in the Strichartz space $L^{7/3} L^{14/3}$. 

\paragraph{Main result} Our main result of this section is 

\begin{proposition} \label{proposition odd non-ra}
 Let $N(x,t) \in \mathcal{N}_0$ be an odd non-radiative solution and $\nu_0 = 1.86$, $\nu_1 = 0.05$. Then $N$ and its initial velocity satisfy the inequality
 \begin{align*}
  &0<N(x,t)<|x|^{-3/2} \varphi_{\nu_0} (t/|x|),& & |x|>t>0;\\
  &0<N_t(x,0)<\nu_0 |x|^{-5/2},& & x \in \Rm^5\setminus\{0\}. 
 \end{align*} 
 In addition, there exists a radius $R$, such that 
 \begin{align*}
  &N(x,t)\geq |x|^{-3/2} \varphi_{\nu_1} (t/|x|),& & R>|x|>t>0;\\
  &N_t(x,0)\geq \nu_1 |x|^{-5/2},& & 0<|x|<R. 
 \end{align*} 
 Please note that the radius $R$ depends on the scaling of $N \in \mathcal{N}_0$. 
\end{proposition}
The major part of this section will be devoted to the proof of this proposition. One of the most important observation is that the radial linear propagation operator $u = \mathbf{S}_L(0,u_1)$ (as given in \eqref{positive SLu1}) is a positive operator of $u_1$ in the exterior region $\Omega_0$, i.e. $u$ is nonnegative as long as $u_1$ is. In addition, if $u_1$ is positive, then $u$ is also positive in the affected region. 

\begin{remark}
 In this section we will frequently utilize numerical simulation. All the related codes and data are available upon request. Please contact the author by Email. 
\end{remark}

\begin{remark}
 We conjecture that when $r\rightarrow 0^+$, any odd non-radiative solution $N \in \mathcal{N}_0$ satisfies the following asymptotic behaviour in the exterior region $\Omega_0$:
 \begin{align*}
  &r^{5/2} N_t(r,0) \rightarrow \nu_2;& &\frac{ N(r,t)}{r^{-3/2} \varphi_{\nu_2}(t/r)} \rightarrow 1.  
 \end{align*}
 Here $\nu_2 \approx 1.575$ is the constant making $(1-y^2)|\varphi'_{\nu_2}(y)|^2 \rightarrow 0$ as $y\rightarrow 1^-$. This is equivalent to saying that
 \[
  r^2 \partial_t \left(r^{-3/2} \varphi_{\nu_2}(t/r)\right) = r^{-1/2} \varphi'_{\nu_2}(t/r) = (r-t)^{-1/2} \cdot \left(1-\frac{t}{r}\right)^{1/2} \varphi'_{\nu_2} \left(\frac{t}{r}\right)
 \]
 converges to zero along a line $r = t + \hbox{Const}$, i.e. the self-similar solution $r^{-3/2} \varphi_{\nu_2}(t/r)$ becomes a ``non-radiative'' solution. This is also consistent with a rough numerical simulation conducted by the author. The strict proof of this asymptotic behaviour seems to be more difficult. But Proposition \ref{proposition odd non-ra} is sufficient for our application. 
\end{remark}

\subsection{Proof of the first part}

The asymptotic behaviour of odd non-radiative solutions have been discussed in the previous sections. According to Remark \ref{ut point-wise estimate}, we have the point-wise estimate
\[
 N_t (x,0) \simeq |x|^{-3}, \qquad |x|\gg 1. 
\]
Thus the inequality $0<N_t(x,0)<\nu_0 |x|^{-5/2}$ always holds near infinity, more precisely, in the exterior region $\{x: |x|>R\}$ for some large radius $R$. It immediately follows from the positive property of $\mathbf{S}_L(0,u_1)$ that 
\begin{equation} \label{inner property N0}
 0< N(x,t) < |x|^{-3/2} \varphi_{\nu_0} (t/|x|), \qquad |x|\geq t+R, \; t>0; 
\end{equation}
where the latter is exactly the solution to (CP1) with initial data $(0,\nu_0 |x|^{-5/2})$. Now let us push $R$ to zero and prove by a continuity argument. According to Lemma \ref{continuity of non-radiative}, $N_t(x,0)$ must be a continuous function of $x$ for $x\neq 0$. By continuity, if the inequality $0<N_t(x,0)<\nu_0 |x|^{-5/2}$ failed for some $x\neq 0$, then there would exist a radius $R$, such that the inequality still holds for $|x|>R$ but either $N_t(x,0)=0$ or $N_t(x,0)=\nu_0 |x|^{-5/2}$ for $|x| = R$. Again the positive property of $\mathbf{S}_L(0,u_1)$ implies that \eqref{inner property N0} still holds for this radius $R$. By dilation we may assume $R=1$, without loss of generality. Now we recall that $N(r,t)$ satisfies the equation 
\[
 (\partial_t^2 - \partial_r^2) (r^2 N(r,t)) = r^2 |N(r,t)|^{4/3} N(r,t) - 2 N(r,t). 
\]
We integrate along the line $r = t+1$ and utilize the non-radiative assumption to deduce that 
\[
  N_t(1,0) = \int_1^\infty \left[2 N(r,r-1) - r^2 |N(r,r-1)|^{4/3} N(r,r-1)\right] {\rm d} r. 
\]
Now let us consider the function $g(z) = 2z - |z|^{4/3} z$ and rewrite the identity above in the form of 
\begin{equation} \label{expression of R2Nt}
  N_t(1,0) = \int_1^\infty r^{-3/2} g(r^{3/2} N(r,r-1)) {\rm d} r. 
\end{equation}
We observe that $g(z)$ satisfies the following properties
\begin{itemize}
 \item $g(z)>0$ for $z \in (0,z_0)$ and $g(z) < 0$ for $z>z_0$. Here the zero point $z_0 = 2^{3/4} \approx 1.681793$. 
 \item $g(z)$ increases in the interval $[0,z_{\rm max}]$ and decreases in the interval $[z_{\rm max}, +\infty)$. Here the maximum point is $z_{\rm max} = (6/7)^{3/4}\approx 0.890820$ and the maximal value is $g_{\rm max} = (8/7) z_{\max} \approx 1.018080$.
\end{itemize}
By \eqref{inner property N0}, we have 
\begin{equation} \label{upper bound r32 N}
 0<r^{3/2} N(r,r-1)<\varphi_{\nu_0} (1-1/r), \qquad r>1. 
\end{equation}
We first show that $N_t(1,0) < \nu_0$. It suffices to give a suitable upper bound of the right hand of \eqref{expression of R2Nt}. A combination of the monotonicity of $g$ and \eqref{upper bound r32 N} yields
\begin{align*}
 N_t(1,0) \leq \int_1^\infty r^{-3/2} g\left(\min\{\varphi_{\nu_0}(1-1/r), z_{\rm max}\}\right) {\rm d} r.
\end{align*}
A change of variable gives 
\[
 N_t(1,0) \leq \int_0^1 g\left(\min\{\varphi_{\nu_0}(y), z_{\rm max}\}\right) (1-y)^{-1/2} {\rm d} y.
\]
A numerical simulation shows that the integral in the right hand is roughly $1.85024 < \nu_0$. We still need to show $N_t(1,0)>0$. According to \eqref{expression of R2Nt}, this is true if $r^{3/2} N(r,r-1) < z_0$ holds for all $r>1$. We still need to deal with the case when there exists $r>1$ such that $N(r,r-1)\geq z_0$. We will show that in this case the integral in the right hand side of \eqref{expression of R2Nt} is still positive. We first give a lemma. 

\begin{lemma} \label{derivative estimate N0}
 Let $y_0\in (0,1)$ satisfy $\varphi_{\nu_0}(y_0) = z_0$, $\kappa_0 = (1-y_0)/(1+y_0)$ and $g_{-} = -g(\sup \varphi_{\nu_0})$. Then 
 \[
  \frac{\rm d}{{\rm d} r} (r^2 N(r,r-1)) \leq (2r-1)^{-1/2} \nu_0 + (\kappa_0 g_-) r^{-1/2}, \qquad r>1. 
 \]
\end{lemma}
\begin{proof}
 Let $w(r,t) = r^2 N(r,t)$. It satisfies $(\partial_t^2-\partial_r^2) w = r^2 |N(r,t)|^{4/3} N(r,t) - 2 N(r,t)$. An integration along the line $r+t = 2r_0 - 1$ gives 
 \begin{align*}
   \frac{\rm d}{{\rm d} r} (r^2 N(r,r-1))|_{r=r_0} & = w_t (r_0,r_0-1) + w_r (r_0,r_0-1) \\
   & = w_t(2r_0-1,0) + w_r(2r_0-1,0) \\
   &\qquad + \int_{r_0}^{2r_0-1} \left[r^2 N(r,2r_0-1-r)^{7/3}  - 2 N(r,2r_0-1-r)\right] {\rm d} r \\
   & \leq (2r_0-1)^{-1/2} \nu_0 - \int_{r_0}^{2r_0-1} r^{-3/2} g\left(r^{3/2} N(r,2r_0-1-r)\right) {\rm d} r. 
 \end{align*}
 Since we have $0< r^{3/2} N(r,2r_0-1-r) < \varphi_{\nu_0} ((2r_0-1-r)/r) \leq \sup \varphi_{\nu_0}$, the inequality 
 \[
  g\left(r^{3/2} N(r,2r_0-1-r)\right) \geq -g_-
 \]
 always holds. There are two cases. If $r_0 +\kappa_0 r_0 \geq 2r_0 - 1$, then we have 
 \begin{align*}
  \frac{\rm d}{{\rm d} r} (r^2 N(r,r-1))|_{r=r_0} &\leq (2r_0-1)^{-1/2} \nu_0 + \int_{r_0}^{2r_0-1} r^{-3/2} g_- {\rm d} r\\
  & \leq (2r_0-1)^{-1/2} \nu_0 + g_- r_0^{-3/2} (r_0 - 1) \\
  & \leq (2r_0-1)^{-1/2} \nu_0 + (k_0 g_-) r_0^{-1/2}. 
 \end{align*}
 On the other hand, if $r_0+\kappa_0 r_0 < 2r_0-1$, then we observe that $r>r_0+\kappa_0 r_0$ implies 
 \[
  \frac{2r_0-1-r}{r} < \frac{2r_0-r}{r} < \frac{1-\kappa_0}{1+\kappa_0} = y_0. 
 \]
 It follows for these $r$'s that 
 \[
  0< r^{3/2} N(r,2r_0-1-r) \leq \varphi_{\nu_0} \left(\frac{2r_0-1-r}{r}\right) < \varphi_{\nu_0} (y_0) = z_0.
 \]
 Namely, the value of $g\left(r^{3/2} N(r,2r_0-1-r)\right)$ is positive for these $r$. Thus 
 \begin{align*}
  \frac{\rm d}{{\rm d} r} (r^2 N(r,r-1))|_{r=r_0} & \leq (2r_0-1)^{-1/2} \nu_0 - \int_{r_0}^{r_0+\kappa_0 r_0} r^{-3/2} g\left(r^{3/2} N(r,2r_0-1-r)\right) {\rm d} r\\
   & \leq (2r_0-1)^{-1/2} \nu_0 + \int_{r_0}^{r_0+\kappa_0 r_0} r^{-3/2} g_- {\rm d} r\\
   & \leq (2r_0-1)^{-1/2} \nu_0 + (k_0 g_-) r_0^{-1/2}. 
 \end{align*}
 This finishes the proof. 
\end{proof}

Now we come back to the integral in \eqref{expression of R2Nt}. We consider a decreasing sequence of values: 
\[
 z_0\approx 1.681793, z_1=1.6, z_2=1.5, z_3=1.4, \cdots, z_{16}=0.1
\]
where $z_0$ is the zero point of $g(z)$, and a corresponding sequence of intervals $J_k$ defined in the following way: $J_0 = (a_0, +\infty)$ with 
\[
 a_0 = \min\{r\geq 1: r^{3/2} N(r,r-1) = z_0\}.
\]
For $1\leq k\leq 16$, we let $J_k = (a_k, b_k)$ with 
\begin{align*}
 &b_k = \min\{r\geq 1: r^{3/2} N(r,r-1) = z_{k-1}\};& & a_k = \max\{r \in [1,b_k): r^{3/2} N(r,r-1) = z_k\}. 
\end{align*}
The definition implies that $b_k \leq a_{k-1}$ and 
\begin{align*}
 &z_k <r^{3/2} N(r,r-1)<z_{k-1}, & & r\in J_k, \; k\geq 1;\\
 &r^{3/2} N(r,r-1) < z_0, & & 1\leq r < a_0. 
\end{align*}
Next we give an upper bound estimate on the ratio $a_k /b_k$ by the lemma above. For any $r\in J_k$, we have
\[
 z_k < r^{3/2} N(r,r-1) \leq \varphi_{\nu_0} (1-1/r)  
\]
Let $y_k = \varphi_{\nu_0}^{-1} (z_k)$. It follows that $y_k < 1-1/r$ and that 
\[
 (2r-1)^{-1/2} = r^{-1/2} (2-1/r)^{-1/2} < r^{-1/2} (1+y_k)^{-1/2}. 
\]
A combination of this with Lemma \ref{derivative estimate N0} yields 
\[
  \frac{\rm d}{{\rm d} r} (r^2 N(r,r-1)) \leq \left[(1+y_k)^{-1/2} \nu_0 + (\kappa_0 g_-)\right] r^{-1/2}, \qquad r\in J_k. 
\]
For convenience we let $\lambda_k = (1+y_k)^{-1/2} \nu_0 + (\kappa_0 g_-) > 1$ be the absolute constants. Integrating from $a_k$ to $b_k$, we deduce from the inequality above that
\[
 b_k^{1/2} z_{k-1} - a_k^{1/2} z_k \leq 2\lambda_k b_k^{1/2}  - 2 \lambda_k a_k^{1/2} \qquad \Rightarrow \qquad \frac{a_k}{b_k} \leq \gamma_k \doteq \left(\frac{2\lambda_k - z_{k-1}}{2\lambda_k - z_k}\right)^2 < 1. 
\]
Please note that $\gamma_k$ are also absolute constants. Now we are able to give a lower bound of the integral in \eqref{expression of R2Nt}
\begin{align*}
 N_t(1,0) & = \int_1^\infty r^{-3/2} g(r^{3/2} N(r,r-1)) {\rm d} r \\
 & \geq \sum_{k=1}^{16} \int_{a_k}^{b_k} r^{-3/2} g(r^{3/2} N(r,r-1)) {\rm d} r + \int_{a_0}^\infty r^{-3/2} g(r^{3/2} N(r,r-1)) {\rm d} r.
\end{align*}
Here we recall the facts that $J_k$'s are disjoint and that the integrand is positive for all $r<a_0$. Next we let 
\[
 m_k = \min\{g(z): z_k \leq z \leq z_{k-1}\} = \min\{g(z_{k-1}), g(z_k)\}, \qquad 1\leq k\leq 16. 
\]
It follows from the range of $r^{3/2} N(r,r-1)$ in each interval $J_k$ and the inequality $a_k/b_k \leq \gamma_k$ that
\begin{align*}
 N_t(1,0) &\geq  \sum_{k=1}^{16} \int_{a_k}^{b_k} r^{-3/2} m_k {\rm d} r - \int_{a_0}^\infty r^{-3/2} g_- {\rm d} r\\
 & \geq \sum_{k=1}^{16} \left[2m_k (a_k^{-1/2} - b_k^{-1/2})\right] - 2 a_0^{-1/2} g_- \\
 & \geq \sum_{k=1}^{16} \left[2m_k b_k^{-1/2} (\gamma_k^{-1/2}-1)\right] - 2 a_0^{-1/2} g_-.
\end{align*} 
By the fact that $b_k \leq a_{k-1}$ we also have 
\[
 b_k \leq a_{k-1} \leq \gamma_{k-1} b_{k-1} \leq \cdots \leq a_0 \prod_{j=1}^{k-1} \gamma_j. 
\]
Inserting this upper bound, we obtain 
\begin{align*}
 N_t(1,0) &\geq \sum_{k=1}^{16} \left[2m_k a_0^{-1/2} \left(\prod_{j=1}^k \gamma_j^{-1/2} - \prod_{j=1}^{k-1} \gamma_j^{-1/2} \right)\right] - 2 a_0^{-1/2} g_-\\
 & = \sum_{k=1}^{16} \left[2m_k a_0^{-1/2} \left(\prod_{j=1}^k \frac{2\lambda_j-z_j}{2\lambda_j-z_{j-1}} - \prod_{j=1}^{k-1} \frac{2\lambda_j-z_j}{2\lambda_j-z_{j-1}} \right)\right] - 2 a_0^{-1/2} g_-.
 \end{align*}
 In order to show $N_t(1,0)>0$, we only need to verify the following inequality 
 \[
  \sum_{k=1}^{16} \left[m_k \left(\prod_{j=1}^k \frac{2\lambda_j-z_j}{2\lambda_j-z_{j-1}} - \prod_{j=1}^{k-1} \frac{2\lambda_j-z_j}{2\lambda_j-z_{j-1}} \right)\right] > g_-. 
 \]
 Both sides are absolute constants. The author conducts a numerical simulation and verifies this inequality. The summary of data is displayed in the table \ref{summary of data}. 
 
 \begin{table}[h]
\caption{Summary of Data}
\begin{center}
\begin{tabular}{|c|c|c|c|c|c|c|}\hline
 k & range & $y_k = \varphi_{\nu_0}^{-1}(z_k)$ & $\lambda_k$ & $\min g$ & $\displaystyle \prod_{j=1}^k \frac{2\lambda_j\!-\!z_j}{2\lambda_j\!-\!z_{j-1}}$ & Contribution \\
 \hline
 $1$ & $1.6-z_0$ & $0.928249$ & $1.349242$ & $0.000000$ & $1.080450$ & $0.000000$\\ 
 \hline
 $2$ 	& $1.5-1.6$ & $0.874605$ & $1.368272$ & $0.205806$ & $1.175515$ &	$0.019565$\\ 
 \hline
 $3$ & $1.4-1.5$ & $0.814979$ & $1.390407$ & $0.424393$ & $1.267293$ &	$0.038950$\\ 
 \hline
 $4$ & $1.3-1.4$ & $0.752686$ & $1.414727$ & $0.607370$ & $1.355949$ & $0.053847$\\ 
 \hline
 $5$ & $1.2-1.3$ & $	0.689656$ & $1.440692$ & $0.755546$ & $1.441694$ & $0.064784$\\ 
 \hline
 $6$ & $1.1-1.2$ & $0.626943$ & $1.468010$ & $0.869772$ & $1.524739$ & $0.072231$\\ 
 \hline
 $7$ & $1.0-1.1$ & $0.565085$ & $1.496548$ & $0.950941$ & $1.605282$ & $0.076591$\\ 
 \hline
 $8$ 	& $0.9-1.0$ & $0.504328$ & $1.526275$ & $1.000000$ & $1.683491$ & $0.078209$\\ 
 \hline
 $9$ & $0.8-0.9$ & $0.444745$ & $1.557229$ & $1.005877$ & $1.759513$ & $0.076469$\\ 
 \hline
$10$ & $0.7-0.8$ & $0.386316$ & $1.589504$ & $0.964927$ & $1.833473$ & $0.071366$\\ 
\hline
$11$ & $0.6-0.7$ & $0.328960$ & $1.623233$ 	& $0.896364$ & $1.905474$ & $0.064539$\\ 
\hline
$12$ & $0.5-0.6$ & $0.272566$ & $1.658595$ & $0.801575$ & $1.975601$ & $0.056212$\\ 
\hline
$13$ & $0.4-0.5$ & $0.217003$ & $1.695814$ & $0.682111$ & $2.043922$ & $0.046603$\\ 
\hline
$14$ & $0.3-0.4$ & $0.162126$ & $1.735163$ & $0.539751$ & $2.110492$ & $0.035931$\\ 
\hline
$15$ & $0.2-0.3$ & $0.107779$ & $1.776980$ & $0.376608$ & $2.175352$ & $0.024426$\\ 
\hline
$16$ & $0.1-0.2$ & $0.053795$ & $1.821679$ & $0.195358$ & $2.238527$ & $0.012342$\\ 
  \hline
  Total & $0.1-z_0$ & - & - & - & - & $\mathbf{0.792065}$\\
  \hline
  \multicolumn{2}{|c|}{$\sup \varphi_{\nu_0}\approx 1.860262$} & \multicolumn{2}{c|}{$y_0 \approx 0.964141$} & \multicolumn{2}{c|}{$\kappa_0 \approx 0.018257$} & $g_-\approx \mathbf{0.535522}$ \\
  \hline
\end{tabular}
\end{center}
\label{summary of data}
\end{table}

In summary, we always have $0<r^{5/2} N_t(r,0) < \nu_0$. The inequality $0<r^{3/2} N(r,t)<\varphi_{\nu_0}(t/r)$ then follows from the monotonicity. 

\subsection{Proof of part two}

Now we show that $r^{5/2} N_t(r,0)$ is bounded from the below by $\nu_1$ when $r$ is small. The inequality $r^{3/2} N(r,t) \leq \varphi_{\nu_1}(t/r)$ for small numbers $t<r$ then follows from the monotonicity and finite speed of propagation. We start by a lemma 

\begin{lemma} \label{two bounds lemma}
 Let $N \in \mathcal{N}_0$. Then either (I) or (II) of the following holds
 \begin{itemize}
  \item [(I)] The inequality $r^{3/2} N(r,r-1) \leq 0.6$ holds for all $r>1$;
  \item [(II)] There exists a radius $R_0 > 1$, such that the inequality $r^{3/2} N(r,r-1) \leq 0.6$ holds for $r\in (1,R_0]$ and 
  \begin{align*}
     &\int_{R_0}^\infty r^{-3/2} g(r^{3/2} N(r,r-1)) {\rm d} r \geq 0;& &N_t (1,0) \geq C_1 R_0^{-1/2}.& 
  \end{align*}
  Here $C_1\approx 0.184221$ is an absolute constant. 
 \end{itemize}
\end{lemma}
\begin{proof}
 We consider the number 
 \[
  M = \sup_{r > 1}\left(r^{3/2} N(r,r-1)\right)
 \]
 and divide the proof into three cases. The first case $M > z_0$ is the most complicated one. Since we have already proved $0<r^{5/2} N_t(r,0) < \nu_0$, the argument in the last subsection still works. We may choose intervals $J_k = (a_k, b_k)$ as usual and obtain $a_k/b_k \leq \gamma_k$. We have 
 \begin{align*}
  \int_{b_{12}}^\infty  r^{-3/2} g(r^{3/2} N(r,r-1)) {\rm d} r & \geq \sum_{k=1}^{11} \int_{a_k}^{b_k} r^{-3/2} g(r^{3/2} N(r,r-1)) {\rm d} r - \int_{a_0}^\infty g_- r^{-3/2} {\rm d} r \\
  & \geq \sum_{k=1}^{11}  \left[2m_k a_0^{-\frac{1}{2}} \left(\prod_{j=1}^k \frac{2\lambda_j-z_j}{2\lambda_j-z_{j-1}} - \prod_{j=1}^{k-1} \frac{2\lambda_j-z_j}{2\lambda_j-z_{j-1}} \right)\right] - 2a_0^{-\frac{1}{2}} g_- \\
  & \geq 0.
 \end{align*}
 Please note that our numerical simulation shows that the first eleven terms are already big enough to neutralize the negative term, thus we may still reserve the last five terms. It follows that 
 \begin{align*}
 N_t(1,0) & \geq \int_{1}^{b_{12}} r^{-3/2} g(r^{3/2} N(r,r-1)) {\rm d} r\\
  &\geq \sum_{k=12}^{16} \int_{a_k}^{b_k} m_k r^{-3/2} {\rm d} r\\
  &\geq \sum_{k=12}^{16}  \left[2m_k b_{12}^{-1/2} \left(\prod_{j=12}^k \frac{2\lambda_j-z_j}{2\lambda_j-z_{j-1}} - \prod_{j=12}^{k-1} \frac{2\lambda_j-z_j}{2\lambda_j-z_{j-1}} \right)\right]\\
  & \approx 0.184221 b_{12}^{-1/2}. 
 \end{align*}  
Our way to choose $b_{12}$ guarantees that $r^{3/2} N(r,r-1) < 0.6$ for all $r<b_{12}$. This verifies (II) for $R_0 = b_{12}$. The second case is $0.6 < M\leq z_0$. We may still choose $J_k$ for $12 \leq k\leq 16$ as usual and follow the same argument as above to verifies (II). Finally if $M \leq 0.6$, then (I) clearly holds. 
\end{proof}

The following lemma is the key tool of this subsection. 

\begin{lemma} \label{lemma left push up}
 Let $N \in \mathcal{N}_0$, $R>0$ and $c \leq \nu_1 = 0.05$ be a small positive constant. If $r^{5/2} N_t(r,0) \geq c$ holds for $r \in [R, 17 R]$, then $R^{5/2} N_t(R,0) > 11c/10$.   
\end{lemma}
\begin{proof}
 By dilation it suffices to consider the case $R=1$. Now we apply Lemma \ref{two bounds lemma}. If (II) holds with $R_0\leq 9$, then we have 
 \[ 
  N_t(1,0)\geq C_1\cdot 9^{-1/2} > 11\nu_1/10 \geq 11c/10.
 \]
 Thus it suffices to consider case (I) and case (II) with $R_0>9$. In either case, we have 
 \begin{align*}
  &0\leq r^{3/2} N(r,r-1) \leq 0.6, \quad 1\leq r\leq 9;& &\int_9^\infty r^{-3/2} g(r^{3/2} N(r,r-1)) {\rm d} r \geq 0. 
 \end{align*}
 Thus we have 
 \[
  N_t(1,0) \geq \int_1^9 r^{-3/2} g(r^{3/2} N(r,r-1)) {\rm d} r.
 \]
 By finite speed of propagation, our assumption on the initial data implies that $r^{3/2} N(r,r-1) \geq \varphi_c(1-1/r)$ for $r\in (1,9)$. For convenience we consider the solution $\varphi_\ast$ to the linear ordinary differential equation 
 \[
  \left\{\begin{array}{l} (1-y^2)\varphi''(y) - y \varphi'(y) + \frac{9}{4} y = 0, \quad y\in (-1,1);\\
  y(0) = 0; \\
  y'(0) = 1. \end{array} \right.
 \]
 The classic Wronskian techniques show that $\varphi'_c(y) \varphi_\ast(y) > \varphi'_\ast(y) \varphi_c(y)$ for all $y\in (0,1)$, which implies that $\varphi_c(y)/\varphi_\ast(y)$ is an increasing function of $t$. As a result, we have $\varphi_c(y) > c \varphi_\ast(y)>0$ for all $y \in (0,1)$. It follows that 
 \[
  0<c \varphi_\ast(1-1/r) < r^{3/2} N(r,r-1) \leq 0.6, \qquad r\in (1,9). 
 \]
 By the monotonicity of $g$ in $[0,0.6]$ we have
  \[
  N_t(1,0) \geq \int_1^9 r^{-3/2} g(c \varphi_\ast(1-1/r)) {\rm d} r
 \]
 By the conserved quantity 
 \[
 \frac{1}{2}(1-y^2)|\varphi_\ast'(y)|^2 + \frac{9}{8} |\varphi_\ast(y)|^2  = \frac{1}{2}, 
 \]
 we have $\varphi_\ast(y)\leq 2/3$. Thus 
 \[
  g(c\varphi_\ast(y)) = 2 c \varphi_\ast(y) \left(1 - c^{4/3}|\varphi_\ast(y)|^{4/3}/2\right) \geq \frac{99}{50} c \varphi_\ast(y), \qquad y\in (0,1). 
 \]
 Therefore we have 
 \begin{align*}
  N_t(1,0) & \geq \frac{99c}{50}\int_1^9 r^{-3/2}  \varphi_\ast(1-1/r) {\rm d} r = \frac{99c}{50} \int_0^{8/9} (1-y)^{-1/2} \varphi_\ast(y) {\rm d} y. 
 \end{align*}
 Finally a numerical method shows that the value of the integral is roughly $0.604556$ thus $N_t(1,0) > 1.1c$. 
\end{proof}

Finally we are able to prove the second part of Proposition \ref{proposition odd non-ra}. According to the first part, we may find a positive number $c_0$, such that 
\[
 r^{5/2} N_t(r,0) > c_0, \qquad r\in [1,17]. 
\]
Without loss of generality, we may assume that $c_0 < \nu_1$. An application of Lemma \ref{lemma left push up} yields that $N_t(1,0) > c_1 \doteq (11/10)c_0 $. A continuity argument, as well as Lemma \ref{lemma left push up}, then shows that $r^{5/2} N_t(r,0) > c_1$ holds for all $r\in (0,1]$. Repeating this process, we may show that 
\[
 r^{5/2} N_t(r,0) > \left(\frac{11}{10}\right)^k c_0, \qquad \forall r \leq 17^{1-k};
\]
as long as $(11/10)^{k-1} c_0 \leq \nu_1$. Finally we let $k_0$ be the minimal positive integer such that $(11/10)^k c_0 > \nu_1$ and conclude
\[
 r^{5/2} N_t(r,0)  > \left(\frac{11}{10}\right)^{k_0} c_0 > \nu_1, \qquad \forall r < 17^{1-k_0}. 
\] 
Finite speed of propagation and monotonicity then shows that 
\[
 r^{3/2} N(r,t) >  \varphi_{\nu_1} (t/r), \qquad 0<t<r<17^{1-k_0}/2. 
\]

\section{Energy norm and radiation of global solutions}

In this section we discuss the energy norms and radiation part of a global solution to (CP1). 

\subsection{Energy norm estimate}

In this subsection we discuss the upper bound of the least energy norm in a long time interval. The main result of this section is 
\begin{lemma} \label{local boundedness}
Let $L>1$ be a constant. There exists a small absolute constant $\varepsilon_0 > 0$ and a large constant $K_0= K_0(L)\gg 1$, such that if 
\begin{itemize}
 \item $u$ is a radial solution to (CP1) defined in a maximal time interval $(-T_-,T_+)$ with an energy $E\geq 0$;
 \item The initial data $(u_0,u_1)$ satisfy $\|(u_0,u_1)\|_{\mathcal{H}(R)} < \varepsilon_0$;  
\end{itemize} 
then for any time $T\in [R, T_+/L)$, there exists a time $t \in [T, LT]$ satisfying
\[
 \|\vec{u}(\cdot,t)\|_{\dot{H}^1 \times L^2}^2 \leq 20 E+K_0. 
\]
\end{lemma}
\begin{remark}
 The proof is based on the virial identity. This argument dates back to Levine \cite{negativeenergy}. Levine showed that any solution with a negative energy must blow up in finite time. It has been proved by a similar argument that if $T_+ = +\infty$, then 
 \[
  \liminf_{t\rightarrow+\infty} \|\vec{u}(t)\|_{\dot{H}^1\times L^2}^2 \leq 5 E. 
 \]
 The upper bound given in Lemma \ref{local boundedness} is larger but applies to a finite (but long) time interval. Please note that this argument does not depend on the radial assumption. Since the proof in 5D is similar to the 3D case, we skip some details. 
\end{remark}
\begin{proof}
We assume that $\|\vec{u}(t)\|_{\dot{H}^1 \times L^2}^2 > 20 E+K_0$ for all $t\in [T,LT]$ and deduce a contradiction, as long as $K_0$ is sufficiently large. Let $\varphi: \Rm \rightarrow [0,1]$ be a smooth cut-off function satisfying 
 \[
  \varphi(s) = \left\{\begin{array}{ll} 1, & s\leq 2\\ 0, & s\geq 3. \end{array}\right.
 \]
 and $\phi(s) = \varphi^2(s)$. We then define 
 \[
  J(t) = \int_{\Rm^5} |u(x,t)|^2 \phi(|x|/t) {\rm d} x.
 \]
 A straight-forward calculation yields 
 \[
  J'(t) = 2 \int_{\Rm^5} u u_t \phi(|x|/t) {\rm d} x - \int_{\Rm^5} |u|^2 \phi'(|x|/t) \frac{|x|}{t^2} {\rm d} x; 
 \]
 and
 \begin{align*}
  J''(t) & = 2 \int_{\Rm^5} |u_t|^2 \phi(|x|/t) {\rm d} x + 2 \int_{\Rm^5} u u_{tt} \phi(|x|/t) {\rm d} x - 4 \int_{\Rm^5} u u_t \phi'(|x|/t) \frac{|x|}{t^2} {\rm d} x \\
  & \qquad + \int_{\Rm^5} |u|^2 \phi''(|x|/t) \frac{|x|^2}{t^4} {\rm d} x + 2\int_{\Rm^5} |u|^2 \phi'(|x|/t) \frac{|x|}{t^3} {\rm d} x. 
 \end{align*}
 Inserting the equation $u_{tt} = \Delta u + |u|^{4/3} u$ and integrating by parts, we obtain 
 \begin{align*}
  J''(t) & = 2 \int_{\Rm^5} (|u_t|^2-|\nabla u|^2 +|u|^{10/3}) \phi(|x|/t) {\rm d} x - 2\int_{\Rm^5} \phi'(|x|/t)u \nabla u \cdot \frac{x}{|x|t} {\rm d} x \\
  & \quad - 4 \int_{\Rm^5} u u_t \phi'(|x|/t) \frac{|x|}{t^2} {\rm d} x + \int_{\Rm^5} |u|^2 \phi''(|x|/t) \frac{|x|^2}{t^4} {\rm d} x + 2\int_{\Rm^5} |u|^2 \phi'(|x|/t) \frac{|x|}{t^3} {\rm d} x.
 \end{align*}
 By the finite speed of propagation, the small data theory, Hardy's inequality, we have 
 \begin{equation} \label{exterior estimate ft}
  \int_{|x|>|t|+R} \left(|\nabla u(x,t)|^2 + |u_t(x,t)|^2 + \frac{|u(x,t)|^2}{|x|^2} + |u(x,t)|^{10/3}\right) {\rm d} x \lesssim_1 \varepsilon_0^2. 
 \end{equation}
 Combining this with the facts 
 \begin{itemize}  
  \item $\phi(|x|/t)-1$ is nonzero only for $|x|>2 t \geq t+R$;
  \item $\phi'(|x|/t)$ and $\phi''(|x|/t)$ are nonzero only for $t+R\leq 2t < |x| < 3t$;
 \end{itemize} 
 we may write 
 \begin{align*}
  J''(t) & =  \int_{\Rm^5} \left(2|u_t|^2 - 2|\nabla u|^2 + 2|u|^{10/3} \right) {\rm d} x  + O(\varepsilon_0^2)\\
  & = \int_{\Rm^5} \left(5 |u_t|^2 + |\nabla u|^2\right) {\rm d} x + \frac{1}{3}\|\vec{u}(t)\|_{\dot{H}^1 \times L^2}^2 - \frac{20}{3} E + O(\varepsilon_0^2) > 0. 
 \end{align*}
 Here the error term $O(\varepsilon_0^2)$ satisfies $|O(\varepsilon_0^2)| \lesssim_1 \varepsilon_0^2$. As a result, if $\varepsilon_0$ is sufficiently small, we have 
 \begin{align}
  |J'(t)|^2 & \leq \frac{9}{2} \left(\int_{\Rm^5} u u_t \phi(|x|/t) {\rm d} x\right)^2 + 9 \left(\int_{\Rm^5} |u|^2 \phi'(|x|/t) \frac{|x|}{t^2} {\rm d} x\right)^2 \nonumber\\
  & \leq \frac{9}{2} \left(\int_{\Rm^5} |u_t|^2 {\rm d} x\right)\left(\int_{\Rm^5} |u|^2 \phi^2 (|x|/t) {\rm d} x\right) + \left(\int_{\Rm^5} |u|^2 \phi(|x|/t) {\rm d} x\right) O(\varepsilon_0^2)\nonumber\\
  & \leq \frac{9}{10}J''(t) J(t). \label{JppJJp}
 \end{align}
 It is not difficult to see that $J(t) \in \mathcal{C}^2([T,LT])$. Now we assume that $J''(t)$ takes its minimum value $M_0$ at time $t_0$ in the time interval $[\frac{7L+8}{15}T,\frac{8L+7}{15}T]$. By the expression of $J''(t)$, we have 
 \[
  M_0 \geq \int_{\Rm^5} \left(5 |u_t(x,t_0)|^2 + |\nabla u(x,t_0)|^2\right) {\rm d} x \geq \|\vec{u}(t_0)\|_{\dot{H}^1\times L^2}^2 \geq 20 E + K_0. 
 \]
 By the energy conservation law and the assumption $E\geq 0$, we also have 
 \[
  \int_{\Rm^5} |u(x,t_0)|^{10/3} {\rm d} x \leq \frac{5}{3} \|\vec{u}(t_0)\|_{\dot{H}^1 \times L^2}^2 \leq 5M_0/3. 
 \]
 By H\"{o}lder inequality and $|J'(t)|^2 \leq J''(t) J(t)$ we have 
 \begin{align*}
  J(t_0) & \leq \int_{|x|<3t_0} |u(x,t_0)|^2 {\rm d} x \lesssim_L T^2 M_0^{3/5}; \quad \Longrightarrow \quad |J'(t_0)| \lesssim_L T M_0^{4/5}. 
 \end{align*}
There are two cases, i.e. $t_0\leq \frac{L+1}{2} T$ and $t_0> \frac{L+1}{2} T$. Let us consider the first case. The second case can be dealt with in the same manner. The only difference is that we need to consider the blow up in the negative time direction instead. Let $t_1 = t_0 + \frac{L-1}{30} T$. 
 \begin{align*}
   J'(t_1) & = J'(t_0) + \int_{t_0}^{t_1} J''(t) {\rm d} t;\\
   J(t_1) & = J(t_0) + \frac{L-1}{30} T J'(t_0) + \int_{t_0}^{t_1} \left(t_1-t\right) J''(t) {\rm d} t.
 \end{align*}
 Since $J''(t) \geq M_0 > 0$ for $t\in [t_0,t_1]\subset [\frac{7L+8}{15}T,\frac{8L+7}{15}T]$, the following inequalities hold:
 \begin{align*}
  \int_{t_0}^{t_1} J''(t) {\rm d} t & \gtrsim_L T M_0; \\
  \int_{t_0}^{t_1} \left(t_1-t\right) J''(t) {\rm d} t & \gtrsim_L T^2 M_0. 
 \end{align*}
 This implies that if $K_0 \geq K_0(L)$ is sufficiently large (thus $M_0$ is sufficiently large), then the integral part is the dominating term in the expression of $J(t_1)$ and $J'(t_1)$. Thus $J(t_1), J'(t_1) > 0$. In addition, it is clear that 
 \[
  \int_{t_0}^{t_1} J''(t) {\rm d} t  \geq \frac{30}{(L-1)T} \int_{t_0}^{t_1} \left(t_1-t\right) J''(t) {\rm d} t.
 \]
 Therefore the following inequality holds as long as the constant $K_0 > K_0(T)$ is sufficiently large:
 \begin{align} \label{lower bound of QJJp}
  &\frac{J'(t_1)}{J(t)} \geq \frac{20}{(L-1)T};& &J(t_1), J'(t_1) > 0. 
 \end{align}
 We define $Q(t) = J'(t)/J(t)$ for $t\geq t_1$. By \eqref{JppJJp} and \eqref{lower bound of QJJp} we have
 \begin{align*}
  &Q'(t) = \frac{J''(t) J(t) - [J'(t)]^2}{J(t)^2} \geq \frac{1}{9} Q(t)^2;& &Q(t_1) \geq \frac{20}{(L-1)T}. 
 \end{align*}
 An integration shows that 
 \[
  \frac{(L-1)T}{20} > \frac{1}{Q(t_1)} - \frac{1}{Q(t)} \geq \frac{1}{9} (t-t_1), \qquad  t_1\leq t < LT. 
 \]
 This gives a contradiction since $t_1 \leq \frac{8L+7}{15}T$. 
 \end{proof}

\subsection{Extraction of the radiation part}

In this subsection we extract the radiation part of a global solution to (CP1). A similar result has been proved in various dimensions, even in non-radial case. Please see Duyckaerts-Kenig-Merle\cite{dkmradial, dkm3} and C\^{o}te-Kenig-Lawrie-Schlag\cite{4dprofile}, for instance. But most of the proof for $d> 3$ given in the literature above assumes that the solution is of type II, i.e. 
\[
 \limsup_{t\rightarrow +\infty} \|\vec{u}(t)\|_{\dot{H}^1\times L^2} < +\infty. 
\]
For reason of completeness, we give a short proof in the radial setting without this assumption. The same proof applies to all dimensions $3\leq d \leq 6$, with the following $Y$ norm 
\[
 \|u\|_{Y(J)} = \|u\|_{L^\frac{d+2}{d-2} L^\frac{2(d+2)}{d-2} (J\times \Rm^d)}. 
\] 
\begin{lemma} \label{radiation part}
 Let $u$ be a radial solution to (CP1) defined for all $t\geq 0$. There exists a radial free wave $v_L$ with a finite energy, such that 
 \[
  \lim_{t\rightarrow +\infty} \int_{|x|>t-A} |\nabla_{t,x}(u-v_L)(x,t)|^2 {\rm d} x = 0, \qquad \forall A \in \Rm. 
 \]
\end{lemma}
\begin{proof} 
 First of all, the energy $E$ must be nonnegative, otherwise the solution would blow up in finite time, according to Levine \cite{negativeenergy}. Let $\varepsilon \in (0,\varepsilon_0)$ be a small constant to be determined later. Here $\varepsilon_0$ is the absolute constant in Lemma \ref{local boundedness}. Given a solution to (CP1) defined for all time $t\geq 0$, it is clear that the following inequality holds for any sufficiently large $R>0$:
 \[
  \|\vec{u}(0)\|_{\mathcal{H}(R)} < \varepsilon. 
 \]
 By the small data theory and finite speed of propagation, it immediately follows that 
 \[
  \|\vec{u}(t)\|_{\mathcal{H}(R+t)} < 2\varepsilon, \qquad \forall t>0.
 \]
 According to Lemma \ref{local boundedness}, there exists a time $t_1 \in [R\varepsilon^{-14/3}, 3R\varepsilon^{-14/3}]$, such that 
 \[
  \|\vec{u}(t_1)\|_{\dot{H}^1\times L^2} \leq C(E). 
 \]
 Here $C(E)$ is a constant depending on $E$ only. Now let us consider the linear free wave $u_L = \mathbf{S}_L (t-t_1) \vec{u}(t_1)$. By the Strichartz estimates and finite speed of propagation, we have 
 \[
  \|\chi_{R} u_L\|_{Y([t_1,+\infty))} \lesssim_1 \|\vec{u}(t_1)\|_{\mathcal{H}(t_1+R)} \lesssim_1 \varepsilon. 
 \]
 In addition, if $(x,t)\in \Omega \doteq \{(x,t): t-R<|x|<t+R, t>t_1\}$, then the following point-wise estimate holds
 \begin{align*}
  |u_L(x,t)| & \lesssim_1 |x|^{-3/2} \|u_L(t)\|_{\dot{H}^1(\Rm^5)} \lesssim_1 C(E) |x|^{-3/2}. 
 \end{align*}
 A direct calculation shows that 
 \begin{align*}
  \|\chi_\Omega u_L\|_{Y([t_1,+\infty))} & = \left(\int_{t_1}^\infty \left(\int_{t-R<|x|<t+R}|u_L(x,t)|^{14/3}{\rm d} x\right)^{1/2}{\rm d} t\right)^{3/7}\\
  & \lesssim_1 C(E)  \left(\int_{t_1}^\infty \left(R t^{-3}\right)^{1/2}{\rm d} t\right)^{3/7}\\
  & \lesssim_1 C(E) \varepsilon. 
 \end{align*}
 In summary, we have $\|\chi_{-R} u_L\|_{Y([t_1,+\infty))} \lesssim_E \varepsilon$. Therefore as long as $\varepsilon < \varepsilon(E)$ is sufficiently small, we may apply the small data theory and the finite speed of propagation to conclude that 
 \begin{align*}
  &\|\chi_{-R} u\|_{Y([t_1,+\infty))} \lesssim_E \varepsilon,& & \|\chi_{-R} F(u)\|_{L^1 L^2([t_1,+\infty)\times \Rm^5)}\lesssim_E \varepsilon^{7/3};
 \end{align*}
 which implies that $u$ scatters to some linear free wave in the exterior region $\{x: |x|>t-R\}$ as time tends to positive infinity. 
 
 Since the argument above works for any sufficiently large $R>0$, we may find a sequence $A_k \rightarrow \infty$ and a sequence of free waves $v_k$, such that 
 \[
  \lim_{t\rightarrow +\infty} \int_{|x|>t-A_k} |\nabla_{t,x}(u-v_k)(x,t)|^2 {\rm d} x = 0. 
 \]
By the theory of radiation profiles, there exists a sequence of $G_k \in L^2([-A_k, +\infty))$, such that 
\[
 \lim_{t\rightarrow +\infty} \int_{t-A_k}^\infty \left(|r^2 u_r(r,t)+G_k(r-t)|^2 + |r^2 u_t(r,t) - G_k(r-t)|^2 \right) {\rm d} r = 0. 
\]
It is not difficult to see $G_j(s) = G_{k}(s)$ holds almost everywhere for $s > - \min\{A_j,A_k\}$. In addition, we also have 
\[
 \sup_k \|G_k\|_{L^2([-A_k,+\infty))} < +\infty,
\]
thanks to Lemma \ref{local boundedness}. Thus we may define $G \in L^2(\Rm)$ such that the restriction of $G$ on $[-A_k,+\infty)$ is exactly $G_k$ for all $k$. It is clear that 
\[
 \lim_{t\rightarrow +\infty} \int_{t-A}^\infty \left(|r^2 u_r(r,t)+G_k(r-t)|^2 + |r^2 u_t(r,t) - G_k(r-t)|^2 \right) {\rm d} r = 0, \qquad \forall A \in \Rm. 
\]
Finally we let $v_L$ be the linear free wave with radiation profile $G$ (in the positive time direction) and the final conclusion follows from the limit above and the basic properties of radiation fields.  
\end{proof}

\begin{remark} \label{upper bound of EvL and G}
 According to Lemma \ref{local boundedness}, the energy $E(v_L)$ of the linear free wave $v_L$ and the norm of the radiation profile $G$ can be dominated by the energy $E$ of the global solution $u$:
 \[
  E(v_L) \simeq_1 \|G\|_{L^2(\Rm)}^2 \lesssim_1 E + 1. 
 \]
\end{remark}

\section{Soliton resolution of almost non-radiative solutions} 

In this section, we prove that the approximated soliton resolution holds at a large time $t_1$ as long as the radiation in the exterior region $\{(x,t): |x|>|t-t_2|\}$ is very weak in term of Strichartz norm for any time $t_2$ in a long interval $[\ell^{-1} t_1, \ell t_1]$. Here $\ell \gg 1$ is a large constant. This is a major step toward the proof of the main theorem. 

\begin{definition}
 Let $u$ be a radial solution to (CP1) defined for all $t>0$. Given two constants $\ell \gg 1$, $\delta\ll 1$, we say that $t$ satisfies the local weak radiation condition $\mathcal{L}(\ell,\delta)$ if and only if the nonlinear radiation profile $G_+$ introduced in \eqref{radiation profile u G} and the initial data $(u_0,u_1)$ of $u$ satisfy the following conditions: 
 \begin{align*}
  &\|G_+\|_{L^2([-\ell t,-\ell^{-1}t])} \leq \delta;& & \|(u_0,u_1)\|_{\mathcal{H}(\ell^{-1} t)} < \delta. 
 \end{align*}
\end{definition}

\begin{remark}
 If $u$ and $t_0$ satisfy the definition above, then the small data theory implies that $u$ can also be defined in the exterior region $\Omega_{\ell^{-1} t_0}$. As a result, $u$ is a well-defined exterior solution in any exterior region $\{x: |x|>|t-t_1|\}$ for $t_1\geq \ell^{-1} t_0$. 
\end{remark}

\begin{lemma} \label{vL lemma}
Let $E> E(W,0)$. If $u$ is a radial solution to (CP1) defined for all $t>0$ with an energy smaller than $E$ and $t_0$ is a time satisfying $\mathcal{L}(\ell, \delta)$, then there exists a linear free wave $v_{t_1,L}$ for each $t_1\in [t_0/2, 3t_0/2]$, such that 
\begin{itemize}
 \item $\displaystyle \lim_{t\rightarrow \pm \infty} \int_{|x|>|t-t_1|} \left|\nabla_{t,x} (u-v_{t_1,L})(x,t)\right|^2 {\rm d} x = 0$; 
 \item $\displaystyle \|\chi_{|x|>|t-t_1|} v_{t_1,L}\|_{Y(\Rm)} + \|v_{t_1,L}\|_{Y([t_1,+\infty))} \lesssim_1 \delta + \ell^{-1/2} E^{1/2}$;
 \item $\|v_{t_1,L}(\cdot,t_1)\|_{L^{10/3}(\Rm^5)} \lesssim_1 \delta + \ell^{-1/2} E^{1/2}$. 
 \item $\displaystyle \int_{|x|<t-\ell^{-1}t_0} \left(|\nabla_{t,x} v_{t_1,L}(x,t)|^2 + \frac{|v_{t_1,L}(x,t)|^2}{|x|^2}\right){\rm d} x \lesssim_1 \delta^2$ holds for all $t\geq t_1$;
 \item $\displaystyle \int_{|x|>t_1+\ell^{-1}t_0} |\nabla_{t,x} v_{t_1,L}(x,t_1)|^2 {\rm d} x \lesssim_1 \delta^2 + \ell^{-1} E$. 
\end{itemize}
\end{lemma}
\begin{proof}
Let $G_\pm$ be the radiation profiles of the nonlinear solution $u$. Our assumption on the initial data implies that 
 \[
  \|G_\pm\|_{L^2([\ell^{-1}t_0, +\infty)} \lesssim_1 \delta. 
 \]
 Now we choose $v_{t_1,L}$ to be the radial free wave with radiation profile (in the negative time direction)
 \[
  G_{t_1} (s) = \left\{\begin{array}{ll} G_-(s), & s>t_1; \\ G_+(-s), & s<t_1.\end{array}\right.
 \]
 Clearly $v_{t_1,L}$ satisfies the asymptotic equivalence condition above. We only need to verify the inequalities. For this purpose we split the radiation profile into two parts: $G_{t_1} = G_{t_1}^1 + G_{t_1}^2$, where $G_{t_1}^j$ are the restriction of $G_{t_1}$ on $\{s: |s|>\ell^{-1} t_0\}$ and $(-\ell^{-1} t_0, \ell^{-1} t_0)$ respectively. We then split $v_{t_1,L}$ into two parts accordingly. Our assumption guarantees that $\|G_{t_1}^1\|_{L^2(\Rm)} \lesssim_1 \delta$. The estimates of this part then follows from the Sobolev, Strichartz and energy estimates. Now let us consider the second part $v_{t_1,L}^2$, which can be given by the explicit formula 
 \[
  v_{t_1,L}^2 (r,t) = \frac{1}{r^3} \int_{t-r}^{t+r} (s-t) G_{t_1}^2 (s) {\rm d}s.
 \]
 It immediately follows that if $r<t-\ell^{-1}t_0$, then $v_{t_1,L}^2 (r,t) = 0$. This immediately gives the energy estimate in the interior region. In addition, an application of Cauchy-Schwarz shows that 
 \[
  |v_{t_1,L}^2 (r,t)| \leq \frac{1}{r^3} \int_{t-r}^{t+r} r \left|G_{t_1}^2 (s)\right| {\rm d} s \lesssim_1 \frac{\ell^{-1/2} t_0^{1/2}}{r^2} \|G_{t_1}^2\|_{L^2}  \lesssim_1 \frac{\ell^{-1/2}(t_0 E)^{1/2}}{r^2}. 
 \]
Here we utilize the upper bound given in Remark \ref{upper bound of EvL and G}. The upper bounds of Strichartz norms and $L^{10/3}$ norm then follow from the support and upper bound given above. Finally if $r>t+\ell^{-1} t_0$ and $t>0$ then we have 
 \[
  v_{t_1,L}^2 (r,t) = \frac{1}{r^3} \int_{-\ell^{-1} t_0}^{\ell^{-1} t_0} (s-t) G_{t_1}^2 (s) {\rm d}s.
 \]
 Thus 
 \[
  \vec{v}_{t_1,L}^2 (r,t_1) = \left(A r^{-3}, B r^{-3}\right),\qquad r>t_1 + \ell^{-1} t_0;
 \]
 where the constants $A$ and $B$ satisfy 
 \begin{align*}
  |A| & = \left|\int_{-\ell^{-1} t_0}^{\ell^{-1} t_0} (s-t_1) G_{t_1}^2 (s) {\rm d}s\right| \lesssim_1 \ell^{-1/2} t_0^{3/2} E^{1/2};\\
  |B| & = \left| - \int_{-\ell^{-1} t_0}^{\ell^{-1} t_0} G_{t_1}^2 (s) {\rm d}s\right| \lesssim_1 \ell^{-1/2} t_0^{1/2} E^{1/2}. 
 \end{align*}
 The energy estimate in the exterior region then follows a direct calculation. 
\end{proof}

\begin{remark} \label{small Strichartz norm}
 A similar argument to the proof of Lemma \ref{vL lemma} shows that any radial free wave $v$ satsifes
 \begin{align*}
  \|\chi_{|x|>|t-t_0|}v\|_{Y(\Rm)} + \|v\|_{Y([t_0,+\infty))} \lesssim_1 \ell^{-1/2}\|G_-\|_{L^2(\{s: |s|<\ell^{-1} t_0\})} + \|G_-\|_{L^2(\{s: |s|>\ell^{-1}t_0\})}. 
 \end{align*}
 Here $G_-$ is the radiation profile of $v$ in the negative time direction. By the symmetry of radiation profiles $G_\pm$, we may also substitute $G_-$ by $G_+$ in the inequality above. 
\end{remark}

The following proposition is the main result of this section. It separates each bubble one-by-one as long as the radiation is sufficiently weak in the corresponding light cone. 
\begin{proposition} \label{main tool} 
 Let $n$ be a positive integer and $E>E(W,0)$. Then (I) there exists a large constant $\ell_0 = \ell_0(n,E)$, a small constant $\delta_0 = \delta_0(n,E)$ and an absolute constant $c_2 \gg 1$, such that if $u$ is a radial solution to (CP1) defined for all $t>0$ with an energy $E(u) \in [E(W,0),E)$ and $t_0$ is a time satisfying $\mathcal{L}(\ell, \delta)$ for some constants $\ell > \ell_0$, $\delta<\delta_0$, then either of the following holds
  \begin{itemize}
  \item [(a)] There exists a sequence $\{(\lambda_j,\zeta_j)\}_{j=1,2,\cdots,J}$ with $0 \leq J\leq n-1$, $\zeta_j \in \{+1,-1\}$ and $0 < \lambda_J < \lambda_{J-1} < \cdots < \lambda_1<t_0$ such that  
  \begin{align*}
  \max\left\{\frac{\lambda_1}{t_0}, \frac{\lambda_2}{\lambda_1}, \cdots, \frac{\lambda_J}{\lambda_{J-1}} \right\}& \leq \kappa_{n,E} (\ell, \delta); \\
   \left\|\vec{u}(t_0)-\sum_{j=1}^J (\zeta_j W_{\lambda_j},0) - \vec{v}_{t_0,L} (t_0)\right\|_{\dot{H}^1\times L^2} 
   & \leq \varepsilon_{n,E} (\ell, \delta);\\
   \left\|\chi_{|x|>|t-t_0|} \left(u - \sum_{j=1}^J \zeta_j W_{\lambda_j} \right)\right\|_{Y(\Rm)}  & \leq \varepsilon_{n,E}(\ell,\delta). 
  \end{align*}
  \item[(b)] There exists a sequence $\{(\lambda_j,\zeta_j)\}_{j=1,2,\cdots,n}$ with $\zeta_j \in \{+1,-1\}$ and $0<\lambda_n < \lambda_{n-1}<\cdots < \lambda_{1}$ such that  
  \begin{align*}
  \max\left\{\frac{\lambda_1}{t_0}, \frac{\lambda_2}{\lambda_1}, \cdots, \frac{\lambda_n}{\lambda_{n-1}} \right\}& \leq \kappa_{n,E} (\ell, \delta); \\
   \left\|\vec{u}(t_0)-\sum_{j=1}^n (\zeta_j W_{\lambda_j},0) - \vec{v}_{t_0,L} (t_0)\right\|_{\mathcal{H}(c_2 \lambda_n)} &  \leq \varepsilon_{n,E} (\ell, \delta);\\
   \left\|\chi_{|x|>|t-t_0|+c_2 \lambda_n} \left(u - \sum_{j=1}^n \zeta_j W_{\lambda_j} \right)\right\|_{Y(\Rm)}  &\leq \varepsilon_{n,E} (\ell, \delta).
  \end{align*}
 \end{itemize}
 Here $v_{t_0,L}$ is the linear free wave given in Lemma \ref{vL lemma}. $\varepsilon_{n,E}(\ell, \delta)$ and  $\kappa_{n,E}(\ell, \delta)$ are positive functions of $\ell$ and $\delta$ satisfying (these functions depend solely on $n$ and $E$)
 \begin{align*}
  &\lim_{\ell \rightarrow +\infty, \delta \rightarrow 0^+} \varepsilon_{n,E} (\ell, \delta) = 0;& & \lim_{\ell \rightarrow +\infty, \delta \rightarrow 0^+} \kappa_{n,E} (\ell, \delta) = 0. 
 \end{align*}
 (II) Furthermore, given any positive constant $c \leq c_2$, there exist two constants $\ell_0(n,E,c)$, $\delta_0(n,E,c)$ and a function $\varepsilon_{n,E,c}(\ell, \delta)$ satisfying 
 \[
  \lim_{\ell \rightarrow +\infty, \delta \rightarrow 0^+} \varepsilon_{n,E,c} (\ell, \delta) = 0,
 \]
 such that if $(u,t_0)$ is a pair of solution/time discussed above in case (b) with $\ell > \ell_0(n,E,c)$ and $\delta<\delta_0(n,E,c)$, then 
 \begin{align*}
    \left\|\vec{u}(t_0)-\sum_{j=1}^n (\zeta_j W_{\lambda_j},0) - \vec{v}_{t_0,L}(t_0)\right\|_{\mathcal{H}(c \lambda_n)} + \left\|\chi_{|x|>|t-t_0|+c \lambda_n} \!\!\left(u - \sum_{j=1}^N \zeta_j W_{\lambda_j} \right)\right\|_{Y(\Rm)}  \leq \varepsilon_{n,E,c} (\ell, \delta).
  \end{align*}
  Here $(\zeta_j, \lambda_j)$ are still the parameters given in case (b) of part (I). 
\end{proposition} 

\begin{remark}
 A direct calculation of nonlinear estimate shows that if $\ell > \ell_0(n,E)$ is sufficiently large and $\delta < \delta_0(n,E)$ is sufficiently small, then a solution $u$ in case (a) satisfies 
\[
 \left\|\chi_{|x|>|t-t_0|} \left(F(u) - \sum_{j=1}^n F(\zeta_j W_{\lambda_j}) \right)\right\|_{L^1 L^2(\Rm \times \Rm^5)} \leq \varepsilon_{n,E}(\ell, \delta).
\]
 Similarly, if $c<c_2$, $\ell>\ell_0(n,E,c)$ and $\delta < \delta_0(n,E,c)$, then a solution in case (b) satisfies 
\[
 \left\|\chi_{|x|>|t-t_0|+c \lambda_n} \left(F(u) - \sum_{j=1}^n F(\zeta_j W_{\lambda_j}) \right)\right\|_{L^1 L^2(\Rm \times \Rm^5)} \leq \varepsilon_{n,E,c}(\ell,\delta). 
\]
Please note that the functions $\varepsilon_{n,E}(\ell,\delta)$ and $\varepsilon_{n,,E,c}(\ell, \delta)$ here may be different from those in the proposition but still satisfy the same limit conditions. For convenience in this section the notation $\varepsilon_{n,E}(\ell, \delta)$ represent a positive function of $n$, $E$, $\ell$ and $\delta$, which satisfies 
\[
 \lim_{\ell\rightarrow +\infty, \delta \rightarrow \delta^+} \varepsilon_{n,E}(\ell, \delta) = 0
\]
for any given positive integer $n$ and $E>E(W,0)$. It may represent different functions at different places. The notations $\varepsilon_{n,E,c}(\ell,\delta)$ can be understood in the same way. Similarly $\delta_0(n,E)$, $\delta_{0}(n,E,c)$ or similar notations represent positive constants depending on $n$ and $E$ (or $n$ , $E$ and $c$). Again they may represent different constants at different places. 
\end{remark} 

The rest of this section is devoted to the proof of Proposition \ref{main tool}. We will apply an induction in the positive integer $n$. More precisely we split the proposition into part I and II, as marked in the proposition. We prove via a bootstrap argument that
\begin{itemize}
 \item[(a)] Part I holds for $n=1$; 
 \item[(b)] Part II holds if Part I holds, for any given $n\geq 1$; 
 \item[(c)] Part I holds for $n+1$ as long as the whole proposition holds for $n$. 
\end{itemize}
The step (b) is simply an application of Lemma \ref{lemma connection 2}. More details can be found in \cite{dynamics3d}. Here we only give the details of step (a) and (c). The proof of these two steps is roughly the same and differs only in the beginning part. We first choose a few constants. Let $\rho$ be the small constant in classification of non-radiative solutions. Without loss of generality, we may assume $\rho < \beta/2$, where $\beta$ is the constant in Lemma \ref{lemma connection}, and let $c_2$ be the large radius such that the radiation residues $\tau_j$ of $W$ satisfy 
\begin{align*}
 &\tau_1(r) = 0;& &\tau_2(c_2) = \rho;& &0<\tau_2(r) < \rho, \quad r>c_2. 
\end{align*}
By Lemma \ref{small tau}, without loss of generality, we may further reduce the value of $\rho$ such that any radial non-radiative solution $v$ to (CP1) with radiation residues $\vec{\tau}$ satisfying
\begin{align*}
 &|\vec{\tau}(R)| = \rho,& &|\vec{\tau}(r)|<\rho, \quad r>R
\end{align*}
must satisfy 
\[
 \|\chi_R v\|_{Y(\Rm)} < \eta/3,
\]
where $\eta$ is the constant in Lemma \ref{lemma connection}. Next we introduce a few notations. Let $v_{t_0,L}$ be the linear free wave given in Lemma \ref{vL lemma}. Given $n\geq 0$ and $(\zeta_j,\lambda_j) \in \{+1,-1\}\times \Rm^+$ for $j=1,2,\cdots,n$, we define 
\begin{align*}
 S_n(x,t) & = v_{t_0,L} + \sum_{j=1}^n \zeta_j W_{\lambda_j}; & e_n(x,t) = \sum_{j=1}^n \zeta_j F(W_{\lambda_j}) - F(S_n).
\end{align*}
In particular, $S_0 = v_{t_0,L}$. The approximated solution $S_n$ solves the following wave equation 
\[
 (\partial_t^2 - \Delta) S_n = F(S_n) + e_n(x,t).
\]
We also define $w_n = u - S_n$, which solves the wave equation 
\[
 (\partial_t^2 - \Delta) w_n = F(u) - F(S_n) - e_n(x,t).
\]
We use the notation $\vec{\tau}_n(r)$ for the radiation residue of $\vec{w}_n(t_0)$ and $\tau_{j,n}(r)$ for its components($j=1,2$). For convenience we also define $\chi_{(r,t_1)}$ to be the characteristic function of the region
\[
 \{(x,t): |x|>r+|t-t_1|\};
\]
and $\chi_{(r_1,r_2,t_1)}$ to be the characteristic function of the region
\[
 \{(x,t): r_1+|t-t_1|<|x|<r_2+|t-t_1|\}. 
\]

\subsection{Beginning of step (a)} 

In this subsection we give the beginning part of Step (a), i.e. Part (I) of Proposition \ref{main tool} for $n=1$. We let $S_0 = v_{t_0,L}$ and $w_0=u-S_0$. Please note this case $e_0 (x,t) = -F(v_{t_0,L})$. Thus 
\begin{align*}
& \|\chi_{(0,t_0)} e_0 (x,t)\|_{L^1 L^2(\Rm \times \Rm^5)} \leq \varepsilon_E (\ell, \delta);& & \|\chi_{(0,t_0)} S_0\|_{Y(\Rm)} \leq \varepsilon_E (\ell,\delta).
\end{align*}
In addition, A combination of finite speed of propagation, small data theory and Lemma \ref{vL lemma} gives 
\begin{equation} \label{tail estimate w0}
 \|\vec{w}_0(t_0)\|_{\mathcal{H}(\ell^{-1} t_0 + t_0)} + \|\chi_{(\ell^{-1}t_0 + t_0,t_0)} (F(u)-F(S_0) - e_0)\|_{L^1 L^2 (\Rm \times \Rm^5)} \leq \varepsilon_E (\ell, \delta). 
\end{equation}
As a result, given any constant $c\in (0,1)$, we may apply Lemma \ref{lemma connection 2} repeatedly to deduce that 
\begin{equation} \label{middle estimate w0}
 \|\vec{w}_0(t_0)\|_{\mathcal{H}(c(\ell^{-1} t_0 + t_0))} \leq \varepsilon_{E,c}(\ell, \delta),
\end{equation}
as long as $\ell$ is sufficiently large and $\delta$ is sufficiently small. Now let us consider the global behaviour of the radiation residue $\vec{\tau}_0(r)$ of $\vec{w}_0(t_0)$. The radiation residue $\vec{\tau}_0(r)$ is very small for $r\geq \ell^{-1} t_0 + t_0$ by \eqref{tail estimate w0}. There are two cases for the behaviour of $\vec{\tau}(r)$ for smaller $r$'s: 

\paragraph{Case one} In this case we assume 
\[
 \sup_{r>0} |\vec{\tau}_0 (r)| \leq \rho < \beta/2. 
\]
An application of Lemma \ref{lemma connection} with $R_1=0$, $R_2= \ell^{-1} t_0 + t_0$ then yields
\[
 \|\vec{w}_0(t_0)\|_{\dot{H}^1\times L^2} + \|\chi_{0,t_0} w_n\|_{Y(\Rm)} \leq \varepsilon_E (\ell, \delta). 
\]
This verifies (a) with $J=0$. 

\paragraph{Case two} In this case we have 
\[
 \sup_{r>0} |\vec{\tau}_0(r)| > \rho. 
\]
Let $R$ be the radius satisfies 
\begin{align*}
 &|\vec{\tau}_0(R)| =\rho;& &|\vec{\tau}_0(r)|<\rho, \quad r>R. 
\end{align*}
By \eqref{middle estimate w0}, the radius $R$ satisfies $R/t_0 \rightarrow 0^+$ as $\ell \rightarrow +\infty$ and $\delta \rightarrow 0^+$, namely
\[
 \frac{R}{t_0} \leq \varepsilon_E (\ell, \delta). 
\]
We temporarily stop here and consider the beginning part of Step (c), and then deal with the remaining part of Step (a) and (c) in a unified way. 

\subsection{Beginning of step (c)}

Let us assume that the proposition holds for a positive integer $n$. We first choose a constant $c_1= c_1(n)$ such that 
\[
 \|\chi_{0,c_1} W\|_{Y(\Rm)} < \frac{\eta}{3n}. 
\]
It is not difficult to see that we only need to consider solutions $u$ and time $t_0$ satisfying case (b) of the proposition for the positive integer $n$. In fact, if $u$ satisfies case (a) for the positive integer $n$, then it also satisfies case (a) for the positive integer $n+1$, with the same choice of $J$ and $(\zeta_j, \lambda_j)$'s. By the induction hypothesis, a solution $u$ in case (b) for $n$ and associated solutions/functions $w_n$, $S_n$ and $e_n(x,t)$ defined at the beginning of this section satisfy 
\begin{align*}
 \left\|\vec{w}_n(t_0)\right\|_{\mathcal{H}(c_1 \lambda_n)} + \|\chi_{(c_1\lambda_n, t_0)}\left(F(u) - F(S_n)- e_n\right)\|_{L^1 L^2(\Rm \times \Rm^5)} & \leq \varepsilon_{n,E}(\ell, \delta); \\
 \|\chi_{(0,t_0)} e_n(x,t)\|_{L^1 L^2(\Rm \times \Rm^5)} & \leq \varepsilon_{n,E}(\ell, \delta); \\
 \left\|\chi_{c_1 \lambda_n} w_n\right\|_{Y(\Rm)} & \leq \varepsilon_{n,E} (\ell,\delta);\\
 \|\chi_{(0,c_1\lambda_n,t_0)} S_n\|_{Y(\Rm)} & \leq \eta/3 + \varepsilon_E (\ell,\delta). 
\end{align*}
In addition, given any positive constant $c<c_1$, the following estimate holds for sufficiently large $\ell$ and sufficient small $\delta$:
\begin{equation} \label{wn estimate clambdan}
 \left\|\vec{w}_n(t_0)\right\|_{\mathcal{H}(c \lambda_n)} + \|\chi_{(c\lambda_n,t_0)} w_n\|_{Y(\Rm)} \leq \varepsilon_{n,E,c}(\ell,\delta). 
\end{equation}
Again we consider the behaviour of the radiation residue $\vec{\tau}_n(r)$ of $\vec{w}_n(t_0)$. There are two cases: 

\paragraph{Case one} In this case we assume 
\[
 \sup_{r>0} |\vec{\tau}_n(t_0)| \leq \rho < \beta/2.
\]
An application of Lemma \ref{lemma connection} with $R_1=0$, $R_2= c_1 \lambda_n$ then gives 
\[
 \|\vec{w}_n(t_0)\|_{\dot{H}^1\times L^2} + \|\chi_{(0,t_0)} w_n\|_{Y(\Rm)} \leq \varepsilon_{n,E} (\ell, \delta). 
\]
This verifies (a) with $J=n$. 

\paragraph{Case two} In this case we have 
\[
 \sup_{r>0} |\vec{\tau}_n(t_0)| > \rho. 
\]
Again we let $R$ be the radius satisfying
\begin{align*}
 &|\vec{\tau}_n(R)| =\rho;& &|\vec{\tau}_n(r)|<\rho, \quad r>R. 
\end{align*}
By \eqref{wn estimate clambdan}, the radius satisfies 
\[
 \frac{R}{\lambda_n} \leq \varepsilon_{n,E} (\ell, \delta).
\]
Now we may complete Step (c) and Step (a) in a unified way. The latter case corresponds to $n=0$. Please note that we have excluded some solutions $u$, which satisfy case (a) and whose corresponding radiation residues $\vec{\tau}_{n}(r)$ are very small for all radii. 

\subsection{Decay of $\tau_1$ and completion of proof}

In this subsection we prove that the first radiation residue $\tau_{1,n}(R)$ must be very small as long as $\ell$ is very large and $\delta$ is very small, in the setting of the last two subsections, which finishes the proof of Proposition \ref{main tool}. More precisely we show that 
\begin{equation} \label{to prove tau1}
 |\tau_{1,n}(R)| \leq \varepsilon_{n,E}(\ell,\delta). 
\end{equation}
Here $\vec{\tau}_n(r) = (\tau_{1,n}(r),\tau_{2,n}(r))$ is the radiation residue of $\vec{w}_n(t_0)$ and $R$ is the radius determined in the last two subsections. We temporarily assume \eqref{to prove tau1} and explain why this implies part (b) holds. We let $\zeta_{n+1} = \tau_2(R)/|\tau_2(R)|$, $\lambda_{n+1} = R/c_2$ and define 
\[
 S_{n+1} = \sum_{j=1}^{n+1} \zeta_j W_{\lambda_j} + v_{t_0,L}, 
\]
which solves the approximated equation 
\begin{align*}
 &(\partial_t^2 u -\Delta u) S_{n+1} = F(S_{n+1}) + e_{n+1}, & &e_{n+1} = \sum_{j=1}^{n+1} F(\zeta_j W_{\lambda_j}) - F(S_{n+1}). 
\end{align*}
Let $w_{n+1}= u - S_{n+1}$ and $\vec{\tau}_{n+1}$ be the radiation residue of $\vec{w}_{n+1}(t_0)$. Combing the ratio inequality 
\[
 \max \left\{\frac{\lambda_1}{t_0}, \frac{\lambda_2}{\lambda_1}, \cdots, \frac{\lambda_{n+1}}{\lambda_n}\right\} \leq \kappa_{n,E}(\ell,\delta), 
\]
with the corresponding estimates for $\vec{w}_n$, we have
\begin{align*}
 \left\|\vec{w}_{n+1} (t_0)\right\|_{\mathcal{H}(c_1 \lambda_n)} + \|\chi_{(c_1\lambda_n, t_0)}\left(F(u) - F(S_{n+1})- e_{n+1}\right)\|_{L^1 L^2(\Rm \times \Rm^5)} & \leq \varepsilon_{n,E}(\ell, \delta); \\
 \|\chi_{(0,t_0)} e_{n+1}(x,t)\|_{L^1 L^2(\Rm \times \Rm^5)} & \leq \varepsilon_{n,E}(\ell, \delta); \\
 \sup_{R\leq r\leq c_1\lambda_n} |\vec{\tau}_{n+1}(r)| & \leq 2\rho < \beta; \\
 \|\chi_{(R,c_1\lambda_n,t_0)} S_{n+1}\| & \leq \frac{2}{3}\eta + \varepsilon_{n,E}(\ell, \delta). 
 \end{align*}
 Here we need to substitute the radius $c_1\lambda_n$ by $t_0+\ell^{-1}t_0$ if $n=0$. When $\ell > \ell(n,E)$ is sufficiently large and $\delta < \delta(n,E)$ is sufficiently small, we may apply Lemma \ref{lemma connection} with $R_1 = R = c_2\lambda_{n+1}$ and $R_2 = c_1\lambda_n$ (or $t_0+\ell^{-1} t_0$ if $n=0$) to conclude 
 \begin{align*}
  \left\|\vec{w}_{n+1} (t_0)\right\|_{\mathcal{H}(c_2 \lambda_{n+1})} + \|\chi_{(c_2\lambda_{n+1},t_0)} w_{n+1}\|_{Y(\Rm)} & \lesssim_1 \left|\vec{\tau}_{n}(R) - \zeta_{n+1}(0,\rho)\right| + \varepsilon_{n,E}(\ell,\delta)\\
  & \lesssim_1 \varepsilon_{n,E}(\ell,\delta). 
 \end{align*}
 Here we utilize \eqref{to prove tau1} and verifies (b) in the proposition. 
 
The remaining task is to prove \eqref{to prove tau1}. This is the most difficult part of the proof, which involves the time evolution and depends on the global behaviour of non-radiative solutions. In the 3-dimensional case, however, this part of argument is not necessary because all non-radiative solutions are ground states, up to the dilation/sign symmetry. By the classification given in the previous sections, it suffices to show that given any constant $T > 0$, the non-radiative class determined by the radiation residue $\vec{\tau}_n (R)$ does not match any non-radiative class of $\vec{U}(t)$ or $-\vec{U}(t)$ for $t\in [-T,T]$, as long as $\ell > \ell_1(n,E,T)$ is sufficiently large and $\delta< \delta_1(n,E,T)$ is sufficiently small. Here $U$ is the universal cover of non-stationary non-radiative solutions. We prove this by a contradiction argument. Let us assume that $\vec{\tau}_n(R)$ matches that of $\vec{U}(T_1)$, without loss of generality, where $T_1\in [-T,T]$. Now let us define 
\[
 S_n^+ = S_n + \frac{1}{(R/R_1)^{3/2}}U\left(\frac{x}{R/R_1}, \frac{t-t_0}{R/R_1}+T_1\right), \qquad |x|>|t-t_0|.
\]
Here $R_1$ is the scale of $U(T_1)$, i.e. the radius at which the radiation residue comes with a length $\rho$, which is uniformly bounded from both the above and the below for all $T_1 \in [-T,T]$. For convenience we let $\lambda = R/R_1$. Clearly we have 
\begin{align*}
 &\frac{\lambda}{\lambda_n} \leq \varepsilon_{n,E,T}(\ell, \delta),\; n\geq 1;& &\frac{\lambda}{t_0} \leq \varepsilon_{E,T}(\ell, \delta),\; n=0.
\end{align*}
The approximated solution $S_n^+$ solves the wave equation 
\[
 (\partial_t^2 - \Delta) S_n^+ = F(S_n^+) + e_n^+.
\]
Here the error term $e_n^+$ is defined by 
\begin{align*}
 &e_n^+ = F(S_n) + F\left(U_\lambda\right) - F(S_n^+) + e_n(x,t); & & U_\lambda = \frac{1}{\lambda^{3/2}} U\left(\frac{x}{\lambda}, \frac{t-t_0}{\lambda}+T_1\right). 
\end{align*}
We recall the Strichartz estimates and Corollary \ref{decay of initial data non-ra} to deduce that ($n\geq 1$)
\begin{align*}
 \|\chi_{(R,t_0)} e_n^+\|_{L^1 L^2} &\leq \left\|\chi_{(R, R^{\frac{1}{2}}\lambda_n^{\frac{1}{2}},t_0)} (e_n^+-e_n)\right\|_{L^1 L^2} + \left\|\chi_{(R^{\frac{1}{2}}\lambda_n^{\frac{1}{2}},t_0)}  (e_n^+ -e_n)\right\|_{L^1 L^2} +\|\chi_{(R,t_0)}e_n\|_{L^1 L^2} \\
 & \lesssim_1 \left\|\chi_{(R, R^{\frac{1}{2}}\lambda_n^{\frac{1}{2}},t_0)} S_n U_\lambda \right\|_{L^{\frac{7}{6}} L^{\frac{7}{3}}}\left(\left\|\chi_{(R, R^{\frac{1}{2}}\lambda_n^{\frac{1}{2}},t_0)} S_n \right\|_{Y(\Rm)}^{1/3} + \left\|\chi_{(R, R^{\frac{1}{2}}\lambda_n^{\frac{1}{2}},t_0)} U_\lambda \right\|_{Y(\Rm)}^{1/3}\right)\\
 & \quad + \left\|\chi_{(R^{\frac{1}{2}}\lambda_n^{\frac{1}{2}},t_0)} S_n U_\lambda \right\|_{L^{\frac{7}{6}} L^{\frac{7}{3}}}\left(\left\|\chi_{(R^{\frac{1}{2}}\lambda_n^{\frac{1}{2}},t_0)} S_n \right\|_{Y(\Rm)}^{1/3} + \left\|\chi_{(R^{\frac{1}{2}}\lambda_n^{\frac{1}{2}},t_0)} U_\lambda \right\|_{Y(\Rm)}^{1/3}\right)\\
 & \qquad + \varepsilon_{n,E}(\ell, \delta)\\
 & \lesssim_{n,E} \left\|\chi_{(R, R^{\frac{1}{2}}\lambda_n^{\frac{1}{2}},t_0)} S_n \right\|_{Y(\Rm)} + \left\|\chi_{(R^{\frac{1}{2}}\lambda_n^{\frac{1}{2}},t_0)} U_\lambda \right\|_{Y(\Rm)} +  \varepsilon_{n,E}(\ell, \delta)\\
 & \leq \varepsilon_{n,E}(\ell,\delta).
\end{align*}
If $n=0$, then we also have 
\begin{align*}
 \|\chi_{(R,t_0)} e_0^+\|_{L^1 L^2} & = \left\|\chi_{(R,t_0)} \left(F(U_\lambda) - F(S_0^+)\right)\right\|_{L^1 L^2} \\
 & \lesssim_1 \|\chi_{(R,t_0)} v_{t_0,L}\|_{Y(\Rm)} \left(\|\chi_{(R,t_0)} v_{t_0,L}\|_{Y(\Rm)}^{4/3} + \|\chi_{(R,t_0)} U_\lambda\|_{Y(\Rm)}^{4/3}\right)\\
 & \leq \varepsilon_E(\ell, \delta). 
\end{align*}
In addition, the way we choose $\rho$ and $c_1$ implies that 
\begin{align*}
 \|\chi_{(R,c_1\lambda_n,t_0)} S_n^+\|_{Y(\Rm)} \leq n\|\chi_{0,c_1} W\|_{Y(\Rm)} + \|\chi_{(R_1,T_1)} U\|_{Y(\Rm)} + \|\chi_0 v_{t_0,L}\|_{Y(\Rm)} \leq \frac{2}{3}
 \eta + \varepsilon_{E}(\ell, \delta),
\end{align*}
for $n\geq 1$; or 
\begin{align*}
 \|\chi_{(R,t_0+\ell^{-1} t_0,t_0)} S_0^+\|_{Y(\Rm)} \leq  \|\chi_{(R_1,T_1)} U\|_{Y(\Rm)} + \|\chi_0 v_{t_0,L}\|_{Y(\Rm)} \leq \frac{1}{3}
 \eta + \varepsilon_{E}(\ell, \delta),
\end{align*}
for $n=0$. In addition, if we define $w_n^+ = u - S_n^+ = w_n - U_\lambda$, then our choice of $T_1$ guarantees that the radiation residue $\vec{\tau}_n^+(r)$ of $\vec{w}_n^+$ satisfies 
\begin{align*}
 &\vec{\tau}_n^+(R) = 0;& &|\vec{\tau}_n^+(r)| \leq 2\rho < \beta, \quad r\in (R,c_1\lambda_n] \quad (\hbox{or}\; r\in (R,t_0+\ell^{-1} t_0] \; \hbox{if}\; n=0);
\end{align*}
as long as $\ell$ is sufficiently large and $\delta$ is sufficiently small. Finally by induction hypothesis and the decay of non-radiative solution $U_\lambda$, we also have 
\begin{align*}
 \|\chi_{(c_1\lambda_n,t_0)} w_n^+\|_{Y(\Rm)} + \|\chi_{(c_1\lambda_n,t_0)} (F(u)-F(S_n^+)-e_n^+)\|_{L^1 L^2} & \leq \varepsilon_{n,E}(\ell,\delta),\quad n\geq 1;\\
 \|\chi_{(t_0+\ell^{-1}t_0,t_0)} w_0^+\|_{Y(\Rm)} + \|\chi_{(t_0+\ell^{-1}t_0,t_0)} (F(u)-F(S_0^+)-e_0^+)\|_{L^1 L^2} & \leq \varepsilon_{E}(\ell,\delta). 
\end{align*}
Putting all these together, we are able to apply Lemma \ref{lemma connection} and deduce that 
\begin{equation} \label{single time t0 estimate}
 \sup_{t\in \Rm} \|\vec{w}_n^+(t)\|_{\mathcal{H}(R+|t-t_0|)} + \|\chi_{(R,t_0)} w_n^+\|_{Y(\Rm)} \leq \varepsilon_{n,E}(\ell,\delta). 
\end{equation}
It immediately follows that 
\begin{equation} \label{single time t0 estimate 2}
 \|\chi_{(R,t_0)}(F(u)-F(S_n^+)-e_n^+)\|_{L^1 L^2 (\Rm \times \Rm^5)} \leq \varepsilon_{n,E}(\ell, \delta). 
\end{equation}
Next we give an improved estimate for time near $t_0$. For each 
\[
 t_1 \in I \doteq [-\lambda(T+T_1)+t_0, \lambda(T-T_1)+t_0],
 \] 
 we define 
\[
 S_{n,t_1}^+ (x,t) = U_\lambda (x,t) + \sum_{j=1}^n \zeta_j W_{\lambda_j}(x) + v_{t_1,L}(x,t), \quad |x|>|t-t_1|;
\]
and $e_{n,t_1}^+$, $w_{n,t_1}^+$ accordingly. Please note that when $\ell$ is sufficiently large and $\delta$ is sufficiently small, each time $t_1\in I$ satisfies 
\[
 |t_1 -t_0| \leq 2T\lambda\lesssim_T R \ll t_0 \qquad \Rightarrow \qquad t_1 \in [t_0/2,3t_0/2].
\]
By comparing the radiation profiles, we immediately obtain $\|\vec{v}_{t_1,L}-\vec{v}_{t_0,L}\|_{\dot{H}^1\times L^2} \lesssim_1 \delta$. Combining this with \eqref{single time t0 estimate} and \eqref{single time t0 estimate 2}, we obtain 
\begin{align*}
 \sup_{t_1\in I} \left(\left\|\vec{w}_{n,t_1}^+(t_1)\right\|_{\mathcal{H}(|t_1-t_0|+R)} + \left(\|\chi_{(R+|t_1-t_0|,t_1)}(F(u)-F(S_{n,t_1}^+)-e_{n,t_1}^+)\right\|_{L^1 L^2}\right) \leq \varepsilon_{n,E}(\ell, \delta). 
\end{align*}
Now we recall the following property of odd non-radiative solutions given in the previous section 
\begin{align*}
 &\|U_t(\cdot,0)\|_{L^2(\{x: |x|<2\})} = +\infty;& &|U(x,t)| \leq |x|^{-3/2} \varphi_{\nu_0} (|t|/|x|) \lesssim_1 |t|/|x|^{5/2}, \; |x|>|t|. 
\end{align*}
By continuity, we may find a constant $\gamma = \gamma(E) \ll 1/10$, such that 
\begin{equation} \label{lower bound of Ut 20E}
 \int_{\gamma<|x|<2} |U_t(x,t)|^2 {\rm d} x > 20 E + 10^{15}, \qquad \forall t \in [-\gamma,\gamma].
\end{equation}
Since $\|\chi_{(\gamma^{2},t)}  U\|_{Y(\Rm)}$ and the scale of $\vec{U}(t)$ are uniformly bounded for all $t\in [-T,T]$, we may follow a similar argument to the one given above and deduce that the following estimates hold for sufficiently large $\ell$ and sufficiently small $\delta$:
\begin{align*}
 \sup_{t_1\in I} \left\|\chi_{(\gamma^2 \lambda, t_1)}S_{n,t_1}^+\right\|_{Y(\Rm)} & \lesssim_{n,E,T} 1;\\
 \sup_{t_1\in I} \left\|\chi_{(\gamma^2 \lambda, t_1)} e_{n,t_1}^+\right\|_{L^1 L^2(\Rm \times \Rm^5)} & \leq \varepsilon_{n,E,T}(\ell,\delta).
\end{align*}
In addition, we have 
\[
 \frac{R+|t_1-t_0|}{\gamma^2 \lambda} \leq \frac{R_1 \lambda + 2T \lambda}{\gamma^2 \lambda} = \frac{R_1 + 2T}{\gamma^2} \lesssim_{E,T} 1.
\]
All these estimates enable us to apply Lemma \ref{lemma connection 2} repeatedly and obtain 
\begin{equation} \label{estimate time t1 gamma2}
 \sup_{t_1\in I} \left(\left\|\vec{w}_{n,t_1}^+(t_1)\right\|_{\mathcal{H}(\gamma^2 \lambda)} + \left\|\chi_{(\gamma^2 \lambda, t_1)}w_{n,t_1}^+\right\|_{Y(\Rm)}\right) \leq \varepsilon_{n,E,T}(\ell,\delta), 
\end{equation}
as long as $\ell > \ell(n,E,T)$ is sufficiently large and $\delta < \delta(n,E,T)$ is sufficiently small. Now let us show that $\|\vec{u}(t)\|_{\dot{H}^1\times L^2}$ is NOT bounded near the time $t_1^\ast = -\lambda T_1 + t_0 \in I$, which gives a contradiction. We shall show this by the virial identity. Let us choose a smooth cut-off function $\varphi$ satisfying 
 \begin{align*}
  &\varphi(s) = \left\{\begin{array}{ll} 1, & s\leq 2\\ 0, & s\geq 3; \end{array}\right.& &|\varphi'(s)| \leq 2,\; s \in \Rm;
 \end{align*}
 and let $\phi(s) = \varphi^2(s)$. We then define 
 \[
  J(t) = \int_{\Rm^5} |u(x,t)|^2 \phi(|x|/\lambda) {\rm d} x, \qquad t\in I_1 \doteq [t_1^\ast-\lambda \gamma, t_1^\ast + \lambda \gamma]\subset I. 
 \]
 A direction calculation shows that 
 \begin{align*}
  J'(t) & = 2\int_{\Rm^5} u(x,t) u_t(x,t) \phi(|x|/\lambda) {\rm d} x; \\
  J''(t) & = 2\int_{\Rm^5} \left(|u_t(x,t)|^2 + u(x,t)u_{tt}(x,t)\right) \phi(|x|/\lambda) {\rm d} x \\
  & = \int_{\Rm^5} \left(2|u_t(x,t)|^2 - 2 |\nabla u(x,t)|^2 + 2 |u(x,t)|^{10/3}\right)\phi(|x|/\lambda) {\rm d} x \\
  & \qquad  - \int_{\Rm^5} 2 u \phi'(|x|/\lambda)\nabla u \cdot \frac{x}{\lambda |x|} {\rm d} x
 \end{align*}
 By our choice of the universal cover $U$, Lemma \ref{vL lemma}, the estimate \eqref{estimate time t1 gamma2} and the ratio inequality $\lambda/\lambda_n \leq \varepsilon_{n,E,T}(\ell,\delta)$, or $\lambda/t_0\leq \varepsilon_{E,T}(\ell, \delta)$ if $n=0$, we have 
 \begin{equation} \label{neck estimate of energy}
  \sup_{t\in I_1} \int_{2\lambda < |x| < 3\lambda} \left(|\nabla u(x,t)|^2 + |u_t(x,t)|^2 + \frac{|u(x,t)|^2}{|x|^2} + |u(x,t)|^{10/3} \right) {\rm d} x < \frac{1}{6}.
 \end{equation}
 as long as $\ell > \ell(n,E,T)$ is sufficiently large and $\delta < \delta(n,E,T)$ is sufficiently small. Similarly we have the following almost orthogonality properties in the exterior regions 
 \begin{align*}
  \sup_{t\in I_1} \left|\|\vec{u}(t)\|_{\mathcal{H}(3\lambda)}^2 - n \|W\|_{\dot{H}^1(\Rm^5)}^2 - \|\vec{U}_\lambda(t)\|_{\mathcal{H}(3\lambda)}^2 - \|\vec{v}_{t,L}\|_{\dot{H}^1\times L^2}^2\right| & \leq \varepsilon_{n,E,T}(\ell,\delta); \\
  \sup_{t\in I_1} \left|\int_{|x|>3\lambda} |u(x,t)|^{\frac{10}{3}}{\rm d}x - n \|W\|_{L^\frac{10}{3}(\Rm^5)}^{\frac{10}{3}} - \int_{|x|>3\lambda} |U_\lambda(x,t)|^{\frac{10}{3}} {\rm d} x \right| & \leq \varepsilon_{n,E,T}(\ell,\delta). 
 \end{align*} 
 It immediately follows from the almost orthogonality and \eqref{lower energy U} that the following inequality holds for sufficiently large $\ell$ and sufficiently small $\delta$:
 \[
  \int_{|x|>3\lambda} \left(\frac{1}{2}|\nabla u(x,t)|^2 + \frac{1}{2}|u_t(x,t)|^2 - \frac{3}{10}|u(x,t)|^\frac{10}{3}\right) {\rm d}x > 0, \qquad \forall t\in I_1. 
 \]
 A combination of this with \eqref{neck estimate of energy} and the energy conservation law gives 
 \begin{equation} \label{upper bound East}
  E_\ast (t) = \int_{\Rm^5} \left(\frac{1}{2}|\nabla u(x,t)|^2 + \frac{1}{2}|u_t(x,t)|^2 - \frac{3}{10}|u(x,t)|^\frac{10}{3}\right) \phi(|x|/\lambda) {\rm d}x < E + 1, \qquad t\in I_1. 
 \end{equation}
 Using this and \eqref{neck estimate of energy}, we may rewrite 
 \begin{align} \label{lower bound Jpp 1}
  J''(t) \geq \int_{\Rm^5} \left(4|u_t(x,t)|^2 + \frac{4}{5}|u(x,t)|^\frac{10}{3}\right) \phi\left(\frac{|x|}{\lambda}\right) {\rm d} x -4 E_\ast(t) - 1; 
 \end{align}
 and
 \begin{align} \label{lower bound Jpp 2}
  J''(t) \geq \int_{\Rm^5} \left(\frac{16}{3}|u_t(x,t)|^2 + \frac{4}{3}|\nabla u(x,t)|^2\right) \phi\left(\frac{|x|}{\lambda}\right) {\rm d} x - \frac{20}{3} E_\ast(t) - 1. 
 \end{align}
 By \eqref{lower bound of Ut 20E}, \eqref{estimate time t1 gamma2} and Lemma \ref{vL lemma}, we also have 
 \begin{equation} \label{low bound ut J}
  \inf_{t\in I_1} \int_{\gamma \lambda<|x|<2\lambda} |u_t(x,t)|^2 {\rm d} x \geq 20 E +10^{15}-1
 \end{equation}
 for sufficiently large $\ell$ and small $\delta$. Therefore by \eqref{lower bound Jpp 2} we have ($t\in I_1$)
 \[
  J''(t) \geq 5 \int_{\Rm^5} |u_t(x,t)|^2 \phi(|x|/\lambda) {\rm d} x> 100E + 10^{15}; \qquad \Longrightarrow \qquad J''(t) J(t) \geq \frac{5}{4} |J'(t)|^2. 
 \]
 Next we prove 
 \begin{lemma} \label{Jpp to J}
  Let $J(t)$, $\gamma$, $\lambda$ be defined as above and $M>0$, $\gamma_1 \in (0,\gamma \lambda/5]$ be constants. If $J''(t) \geq M$ holds for all $t\in [t_1^\ast-\gamma_1, t_1^\ast+\gamma_1]$, then $|J(t)|\geq \gamma_1^2 M/8$ for all $t \in [t_1^\ast-\gamma_1/2, t_1^\ast+\gamma_1/2]$. 
 \end{lemma}
 \begin{proof}
  We prove this lemma by contradiction. If there existed a time $t_2^\ast \in [t_1^\ast-\gamma_1/2, t_1^\ast+\gamma_1/2]$ such that $|J(t_2^\ast)| < \gamma_1^2 M/8$, then we would have 
  \begin{align*}
   J'(t_2^\ast + \gamma_1/2) & = J'(t_2^\ast) + \int_{t_2^\ast}^{t_2^\ast+\gamma_1/2} J''(t) {\rm d} t; \\
   J(t_2^\ast +\gamma_1/2) & = J(t_2^\ast) + (\gamma_1/2)J'(t_2^\ast) + \int_{t_2^\ast}^{t_2^\ast+\gamma_1/2} (t_2^\ast+\gamma_1/2-t)J''(t) {\rm d} t. 
  \end{align*}
  Here we assume $J'(t_2^\ast) \geq 0$, without loss of generality, otherwise we consider the negative time direction instead. The lower bound of $J''(t)$ and the upper bound of $J(t_2^\ast)$ implies that 
  \[
   |J(t_2^\ast)| < \int_{t_2^\ast}^{t_2^\ast+\gamma_1/2} (t_2^\ast+\gamma_1/2-t)J''(t) {\rm d} t. 
  \]
  This implies that $J'(t_2^\ast+\gamma_1/2) > \gamma_1^{-1} J(t_2^\ast+\gamma_1/2) > 0$. Thus $Q(t) = J'(t)/J(t)$ is well-defined for all $t\in [t_2^\ast+\gamma_1/2, t_1^\ast+\gamma \lambda]$ and satisfies $Q(t_2^\ast+\gamma_1/2) > \gamma_1^{-1}$. In addition, $Q(t)$ satisfies 
  \[
   Q'(t) = \frac{J''(t) J(t) - (J'(t))^2}{J(t)^2} \geq \frac{(J'(t))^2}{4 J(t)^2} = \frac{1}{4} Q^2(t)  > 0. 
  \]
  Thus 
  \[
   \gamma_1> \frac{1}{Q(t_2^\ast+\gamma_1/2)} - \frac{1}{Q(t)} \geq \frac{1}{4} (t - t_2^\ast - \gamma_1/2),\qquad \forall t\in [t_2^\ast +\gamma_1/2, t_1^\ast + \gamma \lambda]. 
  \]
  This is a contradiction when $t = t_1^\ast + \gamma \lambda \geq t_1^\ast + 5\gamma_1$. 
 \end{proof}
 
 An application of this lemma with $M = M_0 = 100 E + 10^{15}$ and $\gamma_1 = \gamma \lambda/5$ immediately gives 
 \begin{equation} \label{lower bound of J initial}
  J(t) \geq \frac{\gamma^2 \lambda^2 M_0}{200}, \qquad t \in [t_1^\ast-\gamma \lambda/10, t_1^\ast + \gamma \lambda/10]. 
 \end{equation}
 On the other hand, if $|x|>|t-t_1^\ast|$, then $U_\lambda(x,t)$ satisfies 
\begin{align*}
 \left|U_\lambda (x,t)\right| = \frac{1}{\lambda^{3/2}}\left|U\left(\frac{x}{\lambda}, \frac{t-t_1^\ast}{\lambda}\right) \right| \leq |x|^{-3/2} \varphi_{\nu_0} \left(\frac{|t-t_1^\ast|}{|x|}\right) \leq \frac{2|t-t_1^\ast|}{|x|^{5/2}}. 
\end{align*}
This immediately gives for $t \in [t_1^\ast-\gamma^2 \lambda, t_1^\ast + \gamma^2 \lambda]$ that 
\begin{equation*} 
 \int_{\gamma^2 \lambda < |x|< 3\lambda} |U_\lambda(x,t)|^2 {\rm d} x \leq 4\sigma_4 |t-t_1^\ast|^2 \ln \frac{3}{\gamma^2} \leq 4\sigma_4 \gamma^3 \lambda^2. 
\end{equation*}
We combine this with \eqref{estimate time t1 gamma2}, the scale separation, Lemma \ref{vL lemma} and the Hardy inequality to deduce for sufficiently large $\ell$ and sufficiently small $\delta$ that 
\begin{equation} \label{upper bound by U}
  \int_{\gamma^2 \lambda < |x|< 3\lambda} |u(x,t)|^2 {\rm d} x \leq 5\sigma_4 \gamma^3 \lambda^2, \qquad t\in [t_1^\ast-\gamma^2 \lambda, t_1^\ast + \gamma^2 \lambda].
\end{equation}
Combining this with \eqref{lower bound of J initial}, we obtain 
\[
 \int_{|x|<\gamma^2 \lambda} |u(x,t)|^2 {\rm d} x > \frac{\gamma^2 \lambda^2 M_0}{250}, \qquad t\in [t_1^\ast-\gamma^2 \lambda, t_1^\ast + \gamma^2 \lambda]. 
\]
Thus 
\[
 \int_{|x|<\gamma^2 \lambda} |u(x,t)|^{\frac{10}{3}} {\rm d} x \geq |B(0,\gamma^2 \lambda)|^{-2/3}\left(\int_{|x|<\gamma^2 \lambda} |u(x,t)|^2 {\rm d} x \right)^\frac{5}{3} \geq \gamma^{-\frac{10}{3}} M_0^{4/3}. 
\]
Next we prove for $k\in \{0\}\cup \mathbb{N}$ that 
\begin{equation} \label{to prove ind} 
 \int_{|x|<\gamma^2 \lambda} |u(x,t)|^{\frac{10}{3}} {\rm d} x \geq \gamma^{-2 (\frac{5}{3})^{k+1}} M_0^{(\frac{4}{3})^{k+1}}, \qquad  t\in [t_1^\ast-2^{-k}\gamma^2 \lambda, t_1^\ast + 2^{-k}\gamma^2 \lambda].
\end{equation}
This follows an induction. The case $k=0$ has been proved above. Now we assume that \eqref{to prove ind} holds for $k$. We combine it with \eqref{upper bound East}, \eqref{lower bound Jpp 1}, \eqref{low bound ut J} and obtain a stronger lower bound estimate of $J''(t)$
\[
 J''(t) \geq \frac{4}{5}\gamma^{-2 (\frac{5}{3})^{k+1}} M_0^{(\frac{4}{3})^{k+1}}, \qquad  t\in [t_1^\ast-2^{-k} \gamma^2 \lambda, t_1^\ast + 2^{-k}\gamma^2 \lambda]. 
\]
We apply Lemma \ref{Jpp to J} and obtain 
\[
  J(t) \geq \frac{2^{-2k}}{10} \gamma^{4-2 (\frac{5}{3})^{k+1}} \lambda^2 M_0^{(\frac{4}{3})^{k+1}}, \qquad  t\in [t_1^\ast-2^{-k-1} \gamma^2 \lambda, t_1^\ast + 2^{-k-1}\gamma^2 \lambda]. 
\]
We compare it with \eqref{upper bound by U} and deduce 
\[ 
 \int_{|x|<\gamma^2 \lambda} |u(x,t)|^2 {\rm d} x \geq \frac{2^{-2k}}{20} \gamma^{4-2 (\frac{5}{3})^{k+1}} \lambda^2 M_0^{(\frac{4}{3})^{k+1}}, \qquad  t\in [t_1^\ast-2^{-k-1} \gamma^2 \lambda, t_1^\ast + 2^{-k-1}\gamma^2 \lambda]. 
\]
Again for these time $t$'s we have 
\begin{align*}
  \int_{|x|<\gamma^2 \lambda} |u(x,t)|^{\frac{10}{3}} {\rm d} x &\geq |B(0,\gamma^2 \lambda)|^{-2/3}\left(\int_{|x|<\gamma^2 \lambda} |u(x,t)|^2 {\rm d} x \right)^\frac{5}{3} \\
  & \geq \frac{2^{-\frac{10}{3}k}}{20^{5/3}|B(0,1)|^{2/3}}\gamma^{-2(\frac{5}{3})^{k+2}} M_0^{\frac{5}{3}(\frac{4}{3})^{k+1}}\\
  & \geq \gamma^{-2(\frac{5}{3})^{k+2}} M_0^{(\frac{4}{3})^{k+2}}.
\end{align*}
This verifies \eqref{to prove ind}, which gives a contradiction, because the $\dot{H}^1$ norm, thus $L^{10/3}$ norm of $u$ must be uniformly bounded in the compact time interval $I_1$. 
\begin{remark}
 Given a positive integer $n$, the values of parameters $J$ and $(\zeta_1,\lambda_1), \cdots, (\zeta_J, \lambda_J)$ are uniquely determined by a pair $(u,t_0)$ if $t_0$ satisfies $\mathcal{L}(\ell,\delta)$ for sufficiently large $\ell$ and sufficiently small $\delta$, via the procedure given above. Please note that a small perturbation of $\lambda_j$'s may still satisfy the conditions given in Proposition \ref{main tool}. 
\end{remark}

\section{Proof of main theorem}

In this section we prove Theorem \ref{main thm}. We first give a quantitative result of soliton resolution for a single time

\begin{lemma} \label{pre soliton resolution} 
 Let $E>E(W,0)$, $\kappa$ and $\varepsilon$ be positive constants. Then exists two small positive constants $\ell_1 = \ell_1(E,\kappa,\varepsilon) > E/\varepsilon^2$ and $\delta_1 = \delta_1 (E, \kappa, \varepsilon) \leq \varepsilon$ such that if  
 \begin{itemize}
  \item $u$ is a solution to (CP1) defined for all $t\geq 0$ with an energy $E(u)\in [E(W,0),E)$;
 \item $t_0$ is a time satisfying the condition $\mathcal{L}(\ell_1,\delta_1)$;
 \end{itemize}
 then there exists a sequence $(\zeta_j, \lambda_j)\in \{+1,-1\}\times \Rm^+$ for $j=1,2,\cdots,J$ with $J\geq 0$ and 
 \begin{align*}
  \max\left\{\frac{\lambda_1}{t_0}, \frac{\lambda_2}{\lambda_1}, \cdots, \frac{\lambda_{J}}{\lambda_{J-1}}\right\} < \kappa,
 \end{align*}
 such that ($v_{t_0,L}$ is defined in Lemma \ref{vL lemma})
 \begin{align*}
   \left\|\vec{u}(\cdot, t_0)-\sum_{j=1}^J \zeta_j (W_{\lambda_j},0) - \vec{v}_{t_0,L} (\cdot,t_0)\right\|_{\dot{H}^1\times L^2} & \leq \varepsilon; \\
   \left|\|\vec{u}(t_0)\|_{\dot{H}^1\times L^2}^2 - J \|W\|_{\dot{H}^1}^2 - 2\sigma_4 \|G_+\|_{L^2([-t_0,+\infty))}^2 \right| & \leq \varepsilon^2; \\
    \left|E(u) - J E(W,0) - \sigma_4 \|G_+\|_{L^2([-t_0,+\infty))}^2\right| &\leq \varepsilon^2. 
 \end{align*}
 \end{lemma}
 \begin{proof}
  Without loss of generality we assume that $\kappa, \varepsilon \ll 1$. We start by choosing a positive integer $n = n(E)$ such that 
  \[
   (n-1) \|W\|_{\dot{H}^1(\Rm^5)}^2 > 20 E + K_0 + 1. 
  \]
 Here $K_0 = K_0(3)$ is the constant in Lemma \ref{local boundedness} with $L=3$. We then apply Proposition \ref{main tool} and obtain that if $(u,t)$ satisfies $\mathcal{L}(\ell,\delta)$ for a large constant $\ell \geq \ell_2(E)$ and $\delta \leq \delta_2(E)$, then either of the following holds 
 \begin{itemize} 
  \item There exists a sequence $(\zeta_j, \lambda_j)\in \{+1,-1\}\times \Rm^+$ for $j=1,2,\cdots,J$ with $0\leq J<n$ and 
 \begin{align*}
  \max\left\{\frac{\lambda_1}{t}, \frac{\lambda_2}{\lambda_1}, \cdots, \frac{\lambda_{J}}{\lambda_{J-1}}\right\} < \kappa_E (\ell,\delta),
 \end{align*}
 such that
 \begin{align*}
   \left\|\vec{u}(\cdot, t)-\sum_{j=1}^J \zeta_j (W_{\lambda_j},0) - \vec{v}_{t,L} (\cdot,t)\right\|_{\dot{H}^1\times L^2} & \leq \varepsilon_E (\ell,\delta); \\
   \left|\|\vec{u}(t)\|_{\dot{H}^1\times L^2}^2 - J \|W\|_{\dot{H}^1}^2 - 2\sigma_4 \|G_+\|_{L^2([-t,+\infty))}^2 \right| & \leq \varepsilon_E (\ell,\delta); \\
    \left|E(u) - J E(W,0) - \sigma_4 \|G_+\|_{L^2([-t,+\infty))}^2\right| &\leq \varepsilon_E(\ell,\delta). 
 \end{align*}
 \item There exists a sequence $(\zeta_j, \lambda_j)\in \{+1,-1\}\times \Rm^+$ for $j=1,2,\cdots,n$ and 
 \begin{align*}
  \max\left\{\frac{\lambda_1}{t}, \frac{\lambda_2}{\lambda_1}, \cdots, \frac{\lambda_{n}}{\lambda_{n-1}}\right\} < \kappa_E (\ell,\delta),
 \end{align*}
 such that
 \begin{align*}
   \left\|\vec{u}(\cdot, t)-\sum_{j=1}^n \zeta_j (W_{\lambda_j},0) - \vec{v}_{t,L} (\cdot,t)\right\|_{\mathcal{H}(c_2 \lambda_n)} & \leq \varepsilon_E(\ell,\delta); \\
   \left|\|\vec{u}(t)\|_{\mathcal{H}(c_2 \lambda_n)}^2 - (n-1) \|W\|_{\dot{H}^1}^2 - \|(W,0)\|_{\mathcal{H}(c_2)}^2 - 2\sigma_4 \|G_+\|_{L^2([-t,+\infty))}^2 \right| & \leq \varepsilon_E (\ell,\delta).
 \end{align*}
 \end{itemize}
 Here the estimates on the energy norm or energy of $u$ is a direct consequence of the soliton approximated resolution and the almost orthogonality, as well as the energy estimate on $v_{t,L}$ given in Lemma \ref{vL lemma}. The calculation is similar to the 3-dimensional case. More details can be found in \cite{dynamics3d}.  Please note that 
 \[
  \|v_{t,L}\|_{\dot{H}^1\times L^2} = 2\sigma_4 \left(\|G_+\|_{L^2([-t,+\infty))}^2 + \|G_-\|_{L^2([t,+\infty))}^2\right) = 2\sigma_4 \|G_+\|_{L^2([-t,+\infty))}^2 + O(\delta^2). 
 \]
 It follows that if $(u,t)$ satisfies $\mathcal{L}(\ell_3,\delta_3)$ for a sufficiently large constant $\ell_3 = \ell_3(E,\kappa,\varepsilon)$ and a sufficiently small $\delta_3 = \delta_3(E,\kappa,\varepsilon) < \delta_0$, then either of the following holds 
 \begin{itemize}
  \item[(a)] The conclusion of Lemma \ref{pre soliton resolution} holds for the time $t$; In this case a comparison of the $\dot{H}^1\times L^2$ norm and the energy estimates gives 
  \[
   \|\vec{u}(t)\|_{\dot{H}^1\times L^2}^2 \leq 5 E(u) + 6\varepsilon^2. 
  \]
  \item[(b)] The soliton resolution holds in the exterior region $\{x: |x|>c_2\lambda_n\}$. The norm of $\vec{u}(t)$ satisfies 
  \[
   \|\vec{u}(t)\|_{\dot{H}^1\times L^2}^2 > (n-1) \|W\|_{\dot{H}^1}^2. 
  \]
 \end{itemize}
 To make sure that case (a) happens at time $t_0$, we choose a large constant $\ell_1 = \ell_1(E,\kappa,\varepsilon)$ and a smaller constant $\delta_1 = \delta_1(E,\kappa,\varepsilon) < \min\{\varepsilon_0, \varepsilon\}$, where $\varepsilon_0$ is the small constant in Lemma \ref{local boundedness}, such that if $(u,t_0)$ satisfies $\mathcal{L}(\ell_1,\delta_1)$, then $(u,t)$ satisfies $\mathcal{L}(\ell_3, \delta_3)$ for all $t\in [t_0/2,3t_0/2]$. Without loss of generality we may choose $\ell_1>E/\varepsilon^2$. By the discussion above, either case (a) or (b) holds for each given time $t\in [t_0/2,3t_0/2]$. According to Lemma \ref{local boundedness}, case (a) holds for at least one time $t \in [t_0/2,3t_0/2]$, since 
 \[
  (n-1) \|W\|_{\dot{H}^1}^2 > 20 E + K_0 + 1 > 20 E(u) + K_0.
 \]
 By the continuity of $\|\vec{u}(t)\|_{\dot{H}^1\times L^2}$ and the inequality 
 \[
  5 E(u) + 6\varepsilon^2 < (n-1) \|W\|_{\dot{H}^1}^2,
 \]
 Case (a) happens for all time $t\in [t_0/2,3t_0/2]$, which finishes the proof.
 \end{proof}
  
 \paragraph{Proof of Theorem \ref{main thm}}  
 Finally we are able to prove Theorem \ref{main thm}. We let $\ell_\ast = \ell_1 (E_0, \kappa, \varepsilon)$ and $\delta_\ast = \delta_1 (E_0, \kappa, \varepsilon) \leq \varepsilon$ be the constants given by Lemma \ref{pre soliton resolution}. Next we choose a new parameter $\tilde{\varepsilon} = \min\{\delta_\ast/4, \varepsilon_0\}$, where $\varepsilon_0$ is the constant given in Lemma \ref{local boundedness}. and let $\ell = \ell_1(E_0,\kappa,\tilde{\varepsilon})$ and $\delta = \delta_1(E_0, \kappa, \tilde{\varepsilon})\leq \tilde{\varepsilon}$  be the corresponding constants given by Lemma \ref{pre soliton resolution}. Without loss of generality we assume $\ell >\ell_\ast$, otherwise we may slightly enlarge the value of $\ell$. It is not difficult to see that $\delta_\ast$, $\delta$, $\ell$ depend on $E_0$, $\varepsilon$ and $\kappa$ only. Now let $u$ be a solution as in the main theorem and $R$ be a large radius such that 
\[
 \|\vec{u}(0)\|_{\mathcal{H}(R)} < \delta.
\]
Since $\delta$ is a small constant, by small data theory the exterior solution with initial data $\vec{u}(0)$ exists in the region $\Omega_R$. By finite speed of propagation this exterior solution coincides with $u$ in the overlapping region of their domains. As a result, we may extend the domain of the solution $u$ to $\Omega_R \cup (\Rm^5 \times \Rm^+)$. Let $G_- \in L^2([R,+\infty))$ be its nonlinear radiation profile in the negative direction, as given in Lemma \ref{scatter profile of nonlinear solution}; and $G_+\in L^2(\Rm)$ be the nonlinear radiation profile in the positive time direction, whose existence follows from Lemma \ref{radiation part}. By small data theory and finite speed of propagation we have 
\[
 \sup_{t\in \Rm} \|\vec{u}(t)\|_{\mathcal{H}(R+|t|)} \leq 2\delta. 
\]
Thus 
\[
 2\sigma_4 \|G_\pm\|_{L^2([R,+\infty))}^2 =   \lim_{t\rightarrow \pm \infty} \int_{|x|>|t|+R} |\nabla_{t,x} u(x,t)|^2 {\rm d} x \leq 4\delta^2 \quad \Rightarrow \quad \|G_\pm\|_{L^2([R,+\infty))} < \delta/2. 
\]
\paragraph{Radiation strength} We start by considering the local radiation strength function defined by 
\[
  \varphi_\ell (t) \doteq \|G_+\|_{L^2([-\ell t, -\ell^{-1} t])}, \qquad t \geq \ell R.
\]
We define the set of time with weak local radiation strength:
\[
 Q = \{t>\ell R: \varphi_\ell (t) < \delta_\ast\}. 
\]
By the continuity of $\varphi_\ell$ and the finiteness $\|G_+\|_{L^2(\Rm)} < +\infty$, the set $Q$ is an open set containing a neighbourhood of $+\infty$. For each $t \in \bar{Q}$, i.e. the closure of $Q$, the pair $(u,t)$ satisfies the condition $\mathcal{L}(\ell_\ast, \delta_\ast)$. As a result, we may apply Proposition \ref{pre soliton resolution} and obtain that there exists a sequence $\{(\zeta_j(t), \lambda_j(t))\}_{j=1,2,\cdots,J(t)}$ satisfying
 \begin{align*}
    &\max\left\{\frac{\lambda_1(t)}{t}, \frac{\lambda_2(t)}{\lambda_1(t)}, \cdots, \frac{\lambda_{J(t)}(t)}{\lambda_{J(t)-1}(t)}\right\} < \kappa;& & \zeta_j(t) \in \{+1,-1\};  
 \end{align*}
   such that 
   \begin{equation} \label{soliton resolution 12}
     \left\|\vec{u}(\cdot, t)-\sum_{j=1}^{J(t)} \zeta_j(t) (W_{\lambda_j(t)},0) - \vec{v}_{t,L} (\cdot,t)\right\|_{\dot{H}^1\times L^2} \leq \varepsilon.
   \end{equation}
   In addition, the energy $E$ of $u$ satisfies 
   \begin{align} \label{energy estimate 12}
    \left|E - J(t) E(W,0) - \sigma_4 \|G_{+}\|_{L^2([-t,+\infty))}^2\right| \leq \varepsilon^2.
   \end{align}
 This energy estimate implies that $J(t)$ is a non-increasing function of $t\in \bar{Q}$ and 
 \[
  J(t) \leq \left\lfloor \frac{\varepsilon^2 + E}{E(W,0)} \right\rfloor. 
 \]
 \paragraph{Determination of stable time periods} Let $J_1 > J_2 > \cdots > J_m$ be all possible values of $J(t)$ for $t\in Q$. We may split $Q$ into a few parts 
 \begin{align*}
  & Q = \bigcup_{k=1}^m Q_k; & & Q_k = \{t\in Q: J(t) = J_k\}. 
 \end{align*}
 By the non-increasing property of $J(t)$, the inequality $t_1<t_2$ holds if $t_1 \in Q_{k_1}$, $t_2 \in Q_{k_2}$ and $k_1<k_2$. By continuity and \eqref{energy estimate 12}, $Q_k$ are all nonempty open sets. Thus each $Q_k$ is a union of disjoint open intervals, each of which is exactly a connected component of $Q_k$. We write 
 \[
  Q_k = \bigcup_{i\geq 1} I_{k,i}.
 \]
 Next we define the set of time with very weak local radiation
 \[
  P = \{t > \ell R: \varphi_\ell(t) \leq \delta\} \subset Q. 
 \]
  We claim that given $k \in \{1,2,\cdots, m\}$, there is at most one open interval $I_{k,i}$, such that $I_{k,i}\cap P \neq \varnothing$. Indeed, if $t \in P \cap Q_k$, then $(u,t)$ satisfies the condition $\mathcal{L}(\ell, \delta)$. Thus we may apply Proposition \ref{pre soliton resolution} and obtain 
 \[
  \left|E - J_k E(W,0) - \sigma_4 \|G_{+}\|_{L^2([-t,+\infty))}^2\right| \leq \tilde{\varepsilon}^2 \leq \delta_\ast^2 /16. 
 \]
 Thus if $t_1, t_2\in Q_k \cap P$ and $t_1<t_2$, then 
 \[
  \sigma_4 \|G_+\|_{L^2([-t_2,-t_1])}^2 \leq \delta_\ast^2/8\quad \Rightarrow \quad \|G_+\|_{L^2([-t_2,-t_1])} < \delta_\ast/8. 
 \] 
 It follows that 
 \[
  \varphi_\ell (t) \leq \varphi_\ell (t_1) + \|G_+\|_{L^2([-t_2,-t_1])} + \varphi_\ell(t_2) < 2\delta + \delta_\ast/8 < \delta_\ast, \qquad t\in [t_1,t_2]. 
 \]
 Thus $[t_1,t_2] \subset Q$. By the non-increasing property of $J(t)$, we have $[t_1,t_2] \subset Q_k$. This means that all times in $Q_k \cap P$, if there are any, are all contained in the same connected component of $Q_k$, which verifies our claim. Now we pick up an open interval $I_{k,i}$ for each $k$ and choose the corresponding stable interval $[a_k,b_k]$ to be its closure. There are two cases:
 \begin{itemize}
  \item If there exists an open interval $I_{k,i} = (a_k,b_k)$ such that $I_{k,i}\cap P \neq \varnothing$, then we choose $[a_k,b_k]$ to be the $k$'s stable time period.
  \item If such open interval does not exist, i.e. $Q_k \cap P=\varnothing$, then we pick up an arbitrary interval $I_{k,i} = (a_k,b_k)$ and choose $[a_k,b_k]$ to be the $k$'s stable time period.
 \end{itemize}
 Please note that $G_+ \in L^2(\Rm)$ implies that the last stable time period must be $[a_m, +\infty)$, namely $b_m = +\infty$. 
 \paragraph{Soliton resolution in stable periods} Now we are able to verify the soliton resolution in each stable time period. Since $[a_k,b_k]\subset \bar{Q}$, the approximated soliton resolution \eqref{soliton resolution 12} holds for each $t\in [a_k,b_k]$. Here (and in the argument below) we need to substitute $[a_k,b_k]$ by $[a_m, +\infty)$ for the last stable period. A combination of \eqref{energy estimate 12} and the continuity of $\|G_+\|_{L^2([-t,+\infty))}$ shows that $J(a_k) = J(b_k) = J_k$ holds for the endpoints as well. Therefore $J(t) = J_k$ is a constant for all $t\in [a_k,b_k]$. This also gives the estimate 
 \[
  \sigma_4 \|G_+\|_{L^2((-b_k,-a_k])}^2 \leq 2\varepsilon^2. 
 \]
 Next we may use the continuity of $\vec{u}(t)$ and $\vec{v}_{t,L}(t)$ to deduce that $\zeta_{j}(t)$ never changes in the time interval $[a_k,b_k]$ for each $1\leq j \leq J_k$.  Thus for each $k$ and $1\leq j \leq J_k$, we may choose 
 \begin{align}
   \zeta_{k,j} = \zeta_{j}(t), \quad t\in [a_k,b_k].
 \end{align}
We still need to substitute $\vec{v}_{t,L}$ by a linear free wave independent of $t$ for each stable time period. We let $v_{k,L}$ be the linear free wave with the following radiation profile in the positive time direction:
 \[
  G_{k,+} (s) = \left\{\begin{array}{ll} G_+(s), & s> -b_k; \\ 0, & s < -b_k. \end{array}\right.
 \]
 By comparing the radiation profiles we have ($t \in [a_k,b_k]$) 
 \begin{align*}
  \|\vec{v}_{t,L}(\cdot,t)-\vec{v}_{k,L}(\cdot,t)\|_{\dot{H}^1 \times L^2}^2 & = 2\sigma_4 \int_{-\infty}^{-b_k} |G_-(-s)|^2 {\rm d} s + 2\sigma_4 \int_{-b_k}^{-t} |G_-(-s)-G_+(s)|^2 {\rm d} s \\
  &\lesssim_1 \|G_+\|_{L^2((-b_k,-t])}^2 + \|G_-\|_{L^2([t,+\infty))}^2 \lesssim_1 \varepsilon^2 + \delta^2 \lesssim_1 \varepsilon^2. 
 \end{align*}
 Thus we have 
  \begin{equation} \label{soliton resolution 15}
     \left\|\vec{u}(\cdot, t)-\sum_{j=1}^{J_k} \zeta_{k,j} (W_{\lambda_{k,j}(t)},0) - \vec{v}_{k,L} (\cdot,t)\right\|_{\dot{H}^1\times L^2} \lesssim_1 \varepsilon, \qquad t\in [a_k,b_k]. 
   \end{equation}
 In addition, we may recall Remark \ref{small Strichartz norm}, as well as the facts $2\sigma_4 \|G_+\|_{L^2}^2 \leq 20 E_0 + K_0(3) \lesssim_1 E_0$ and $\ell > E_0/\varepsilon^2$, to deduce 
  \begin{align*}
  \|\chi_{|x|>|t-a_k|} v_{k,L}\|_{Y(\Rm)} + \|v_{k,L}\|_{Y([a_k,+\infty))} & \lesssim_1 \ell^{-1/2} \|G_+\|_{L^2(\{s: |s|<\ell^{-1}a_k\})} + \|G_+\|_{L^2(-b_k,-a_k)}\\
  & \quad +\|G_+\|_{L^2(-a_k,-\ell^{-1}a_k)} + \|G_+\|_{L^2([\ell^{-1} a_k,+\infty))} \\
  & \lesssim_1 (E_0 /\varepsilon^2)^{-1/2} E_0^{1/2} + \varepsilon + \varphi_\ell(a_k) + \delta \lesssim_1 \varepsilon.
 \end{align*}
 Furthermore, the basic theory of radiation fields gives that 
 \begin{align*}
  \lim_{t\rightarrow +\infty} \int_{|x|<t-\ell^{-1} a_k} |\nabla_{t,x} v_{k,L}(x,t)|^2 {\rm d} x = 2\sigma_4 \|G_{k,+}\|_{L^2(-\infty, -\ell^{-1} a_k)}^2 \lesssim_1 \varepsilon^2.
 \end{align*}
 The finite speed of energy propagation then gives 
 \[
  \int_{|x|<t-\ell^{-1} a_k} |\nabla_{t,x} v_{k,L}(x,t)|^2 {\rm d} x \lesssim_1 \varepsilon^2, \qquad t> \ell^{-1} a_k. 
 \]
 \paragraph{Property of collision periods} Now let us consider the collision time periods $[b_k, a_{k+1}]$. By the choice of $a_k$, $b_k$ and the continuity of $\varphi_\ell$, we must have 
 \[
  \varphi_\ell (b_k) = \varphi_\ell (a_{k+1}) = \delta_\ast, \qquad k =1, 2, \cdots, m-1. 
 \]
 The way we choose $[a_k,b_k]$ guarantees that $P \cap [b_k, a_{k+1}] = \varnothing$. Therefore we have 
 \[
  \varphi_\ell (t) > \delta, \quad t\in [b_k, a_{k+1}]. 
 \]
 Now let us give an upper bound of the ratio $a_{k+1}/b_k$. First of all, \eqref{energy estimate 12} gives 
 \begin{align*}
   & p_k E(W,0) - 2\varepsilon^2 \leq \sigma_4 \int_{-a_{k+1}}^{-b_k} |G_+(s)|^2 {\rm d} s \leq p_k E(W,0) + 2\varepsilon^2; & & p_k = J_k - J_{k+1} \in \mathbb{N}. 
 \end{align*}
 We may combine this upper bound with the lower bound $\varphi_\ell(t) > \delta$ to deduce 
 \[
  \delta^2 \left\lfloor \log_{\ell^2} \frac{a_{k+1}}{b_k} \right\rfloor \leq \int_{-a_{k+1}}^{-b_k} |G_+(s)|^2 {\rm d} s \leq \frac{p_k E(W,0) + 2\varepsilon^2}{\sigma_4}.
 \]
 As a result, there exists a large constant $L = L(E_0, \varepsilon, \kappa)$, such that $a_{k+1}/b_k \leq L$. Similarly we may give the upper bound of the ratio $a_{1}/R$. In fact, if $a_1 > \ell R$, then the way we choose $[a_k,b_k]$ guarantees $P \cap (\ell R, a_1] = \varnothing$. This implies 
 \[
  \varphi_\ell (t) \geq \delta, \qquad t\in [\ell R, a_1]. 
 \]
 We may combine this with \eqref{energy estimate 12} to deduce 
 \[
   \delta^2 \left\lfloor \log_{\ell^2} \frac{a_1}{R} \right\rfloor \leq \int_{-a_1}^{-R} |G_+(s)|^2 {\rm d} s \leq \int_{-a_1}^{\infty} |G_+(s)|^2 {\rm d} s \leq \frac{E-J_1 E(W,0) + \varepsilon^2}{\sigma_4}
 \]
 This gives 
 \[
  a_1/R \leq L = L(E_0, \varepsilon, \kappa). 
 \]
 \paragraph{Completion of proof} Finally we combine the properties of stable/collision/preparation periods given above and complete the proof, except that the upper bounds we obtain are $C \varepsilon$ instead of $\varepsilon$, where $C$ is an absolute constant(or a constant multiple of $\varepsilon^2$ instead of $\varepsilon^2$). A substitution of $\varepsilon$ by $C^{-1} \varepsilon$ finishes the proof.

\section*{Acknowledgement}
The author is financially supported by National Natural Science Foundation of China Project 12471230.


\begin{thebibliography}{99}
% \bibitem{wavedr} J-P. Anker, P. Martinot, E. Pedon, and A. G. Setti. {``The shifted wave equation on Damek--Ricci spaces and on homogeneous trees''} \textit{Trends in Harmonic Analysis} (2013): 1-25.
% \bibitem{wavehyper} J-P. Anker, V. Pierfelice, and M. Vallarino. {``The wave equation on hyperbolic spaces''} \textit{Journal of Differential Equations} 252(2012): 5613-5661.
%\bibitem{interpolationspace} J. Bergh and J. L\"{o}fstr\"{o}m. ``Interpolation Space.'' Springer-Verlag, Berlin Heidelberg New York 1976. 
%\bibitem{stablemanit2} S. Burzio and J. Krieger. {``Type II blow up solutions with optimal stability properties for the critical focusing nonlinear wave equation on $\Rm^{3+1}$.''} \textit{Memoirs of the American Mathematical Society} 278(2022), no. 1369. 
%\bibitem{fchain} M. Christ and M. Weinstein {``Dispersion of small amplitude solutions of the generalized Korteweg-de Vries equation''} \textit{Journal of Functional Analysis} 100(1991): 87-109.
% \bibitem{wkghyper} J-P. Anker, and V. Pierfelice. {``Wave and Klein-Gordon equations on hyperbolic spaces''} (2011), preprint arXiv: 1104.0177v2.
 %\bibitem{ab1} A. Bulut. {``Global well-posedness and scattering for the defocusing energy-supercritical cubic nonlinear wave equation''} \textit{Journal of Functional Analysis} 263(2012): 1609-1660.
% \bibitem{ksph1} M. G. Cowling. {``The Kunze-Stein phenomenon''} \textit{Annals of Mathematics} 107(1978): 209-234.
%\bibitem{bahouri} H. Bahouri, and P. G\'{e}rard. {``High frequency approximation of solutions to critical nonlinear equations.''} \textit{American Journal of Mathematics} 121(1999): 131-175.
%\bibitem{ascattering} H. Bahouri and J. Shatach. {``Decay estimates for the critical semilinear wave equation.''} \textit{Annales de l'Institut Henri Poincar\'{e}. Analyse Non Lin\'{e}aire} 15(1998): 783-789. 
\bibitem{classNR} C. Collot, T. Duyckaerts, C. Kenig and F. Merle. {``On classification of non-radiative solutions for various energy-critical wave equations.''} \textit{Advances in Mathematics} 434(2023), Paper No. 109337: 91pp. 
\bibitem{soliton6d} C. Collot, T. Duyckaerts, C. Kenig and F. Merle. {``Soliton resolution for the radial quadratic wave equation in space dimension 6.''} \textit{Vietnam Journal of Mathematics} 52(2024), no. 3: 735-773.
%\bibitem{interaction} J. Colliander, M. Keel, G. Staffilani, H. Takaoka, and T. Tao. {``Global existence and scattering for rough solutions to a nonlinear Schr\"{o}dinger equation in $\Rm^3$.''} \textit{Communications on Pure and Applied Mathematics} 57(2004): 987-1014.
 %\bibitem{ctao} J. Colliander, M. Keel, G. Staffilani, H. Takaoka, and T. Tao. {``Global well-posedness and scattering in the energy space for the critical nonlinear nonlinear Schr\"{o}dinger equation in $\Rm^3$.''} \textit{Annals of Mathematics} 167(2007): 767-865.
 \bibitem{newradiation} R. C\^{o}te, and C. Laurent. {``Concentration close to the cone for linear waves.''} \textit{Revista Matem\'{a}tica Iberoamericana} 40(2024): 201-250.
 \bibitem{4dprofile} R. C\^{o}te, C. E. Kenig, A. Lawrie and W. Schlag. {``Profiles for the radial focusing 4d energy-critical wave equation.''} \textit{Communications in Mathematical Physics} 357(2018), no. 3: 943-1008.
% \bibitem{channeleven} R. C\^{o}te, C.E. Kenig and W. Schlag. {``Energy partition for linear radial wave equation.''} \textit{Mathematische Annalen} 358, 3-4(2014): 573-607.
 %\bibitem {truncation} M. Christ and A. Kiselev. {``Maximal functions associated to filtrations''} \textit{Journal of Functional Analysis} 179(2001): 409-425.
 %\bibitem{smallgs2} P. \'{D}Ancona, V. Georgiev, and H. Kubo. {``Weighted decay estimates for the wave equation''} \textit{Journal of Differential Equations} 177(2001): 146-208.
% \bibitem{claim1} B. Dodson. {``Global well-posedness and scattering for the radial, defocusing, cubic nonlinear wave equation.''} \textit{arXiv Preprint} 1809.08284.
 %\bibitem{claim2} B. Dodson. {``Global well-posedness for the radial, defocusing nonlinear wave equation for $3<p<5$.''} \textit{arXiv Preprint} 1810.02879.
% \bibitem{cubic3dwave} B. Dodson and A. Lawrie. {``Scattering for the radial 3d cubic wave equation.''} \textit{Analysis and PDE}, 8(2015): 467-497.
 %\bibitem{nonradial3p5} B. Dodson, A. Lawrie, D. Mendelson, J. Murphy {``Scattering for defocusing energy subcritical nonlinear wave equations''}, Analysis and PDE 13(2020): 1995-2090.
 \bibitem{stablet1} R. Donninger. {``Strichartz estimates in similarity coordinates and stable blowup for the critical wave equation.''} \textit{Duke Mathematical Journal} 166(2017), no. 9: 1627-1683.
 \bibitem{moreexamples} R. Donninger, M. Huang, J. Krieger and W. Schlag. {``Exotic blowup solutions for the $u^5$ focusing wave equation in $\Rm^3$.''} \textit{Michigan Mathematical Journal} 63(2014), no. 3: 451-501. 
 %\bibitem{nonscaglobal1} R. Donninger and J. Krieger. {``Nonscattering solutions and blowup at infinity for the critical wave equation.''} \textit{Mathematische Annalen} 357(2013), no. 1: 89-163.  
 \bibitem{djknonradial} T. Duyckaerts, H. Jia and C.E.Kenig {``Soliton resolution along a sequence of times for the focusing energy critical wave equation''}, \textit{Geometric and Functional Analysis} 27(2017): 798-862.
 \bibitem{soliton4d} T. Duyckaerts, C.E. Kenig, Y. Martel and F. Merle. {``Soliton resolution for critical co-rotational wave maps and radial cubic wave equation.''} \textit{Communications in Mathematical Physics} 391(2022), no. 2: 779-871. 
 \bibitem{tkm1} T. Duyckaerts, C.E. Kenig, and F. Merle. {``Universality of blow-up profile for small radial type II blow-up solutions of the energy-critical wave equation.''} \textit{The Journal of the European Mathematical Society} 13, Issue 3(2011): 533-599.
%\bibitem{dkmnonradial} T. Duyckaerts, C. E. Kenig, and F. Merle. {``Universality of blow-up profile for small type II blow-up solutions of the energy-critical wave equation: the nonradial case''} \textit{The Journal of the European Mathematical Society} 14, Issue 5(2012): 1389-1454.
 \bibitem{dkmradial} T. Duyckaerts, C.E. Kenig and F. Merle. {``Profiles of bounded radial solutions of the focusing, energy-critical wave equation''}, \textit{Geometric and Functional Analysis} 22(2012): 639-698.
  \bibitem{se} T. Duyckaerts, C.E. Kenig, and F. Merle. {``Classification of radial solutions of the focusing, energy-critical wave equation.''} \textit{Cambridge Journal of Mathematics} 1(2013): 75-144.
% \bibitem{dkm2} T. Duyckaerts, C.E. Kenig, and F. Merle. {``Scattering for radial, bounded solutions of focusing supercritical wave equations.''} \textit{International Mathematics Research Notices} 2014:  224-258.
 \bibitem{dkm3} T. Duyckaerts, C.E. Kenig, and F. Merle. {``Scattering profile for global solutions of the energy-critical wave equation.''} \textit{Journal of European Mathematical Society} 21 (2019): 2117-2162.
 \bibitem{oddtool} T. Duyckaerts, C. E. Kenig, and F. Merle. {``Decay estimates for nonradiative solutions of the energy-critical focusing wave equation.''} \textit{Journal of Geometric Analysis} 31(2021), no. 3: 7036-7074.
% \bibitem{exteriorW} T. Duyckaerts, C. E. Kenig, and F. Merle. {``Exterior energy bounds for the critical wave equation close to the ground state.''} \textit{arXiv preprint} 1912.07658.
 \bibitem{oddhigh} T. Duyckaerts, C. E. Kenig, and F. Merle. {``Soliton resolution for the critical wave equation with radial data in odd space dimensions.''}  \textit{Acta Mathematica} 230(2023), no. 1: 1-92.
 \bibitem{threshold} T. Duyckaerts and F. Merle. {``Dynamic on threshold solutions for energy-critical wave equation.''} \textit{International Mathematics Research Papers} 2007, no.4: Article ID. rpn002. 
 %\bibitem{sob} M. Cowling, S. Giulini, S. Meda {``$L^p-L^q$ estimates for functions of the Laplace-Beltrami operator on noncompact symmetric spaces I''} \textit{Duke Mathematical Journal} 72(1993): 109-150.
% \bibitem{pdeevans} L. C. Evans {``Partial Differential Equations, Second Edition.''} \textit{Graduate Studies in Mathematics} 19(2010), AMS, Providence.
% \bibitem{wavehyper97} J. Fontaine. {``A semilinear wave equation on hyperbolic spaces''} \textit{Communications in Partial Differential Equations} 22(1997): 633-659. 
 %\bibitem{superhyp} A. French. {``Existence and Scattering for Solutions to Semilinear Wave Equations on High Dimensional Hyperbolic Space''}(2014), Preprint Arxiv: 1407.2695v1.
\bibitem{radiation1} F. G. Friedlander. {``On the radiation field of pulse solutions of the wave equation.''}  \textit{Proceeding of the Royal Society Series A} 269 (1962): 53-65.
%\bibitem{inverseradiation} F. G. Friedlander. {``An inverse problem for radiation fields.''} \textit{Proceeding of the London Mathematical Society} 27, no 3(1973): 551-576.
\bibitem{radiation2} F. G. Friedlander. {``Radiation fields and hyperbolic scattering theory.''} \textit{Mathematical Proceedings of Cambridge Philosophical  Society} 88(1980): 483-515.
%\bibitem{fourierappli} G. B. Folland. {``Fourier analysis and its applications.''} \textit{The Wadsworth and Brooks/Cole mathematics series}, 1992, Pacific Grove, California. 
 %\bibitem{smallgs3} V. Georgiev. \textit{Semilinear hyperbolic equations, MSJ Memoirs 7}, Tokyo: Mathematical Society of Japan, 2000.
% \bibitem{gwpwrn} V. Georgiev, H. Lindblad, and C. Sogge {``Weighted Strichartz estimates and global existence for semilinear wave equations''}, \textit{American Journal of Mathematics} 119(1997): 1291-1319.
%\bibitem{locad1} J. Ginibre, A. Soffer and G. Velo. {``The global Cauchy problem for the critical nonlinear wave equation''} \textit{Journal of Functional Analysis} 110(1992): 96-130.
% \bibitem{globalwell} J. Ginibre, and G. Velo. {``The global Cauchy problem for the nonlinear Klein-Gordon equation.''} \textit{Mathematische Zeitschrift} 189, No 4(1982): 487-505
% \bibitem{conformal2} J. Ginibre, and G. Velo. {``Conformal invariance and time decay for nonlinear wave equations.''} \textit{Annales de l'institut Henri Poincar\'{e} (A) Physique th\'{e}orique} 47(1987): 221-276.
\bibitem{strichartz} J. Ginibre, and G. Velo. {``Generalized Strichartz inequality for the wave equation.''} \textit{Journal of Functional Analysis} 133(1995): 50-68.
 %\bibitem{blowup2} R. T. Glassey. {``Finite-time blow-up for solutions of nonlinear wave equations''} \textit{Mathematische Zeitschrift} 177(1981): 323-340.
 %\bibitem{subconformal2d} R. Glassey and H. Pecher. {``Time decay for nonlinear wave equations in two space dimensions.''} \textit{Manuscripta Mathematica} 38(1982): 387-400.
\bibitem{mg1} M. Grillakis. {``Regularity and asymptotic behaviour of the wave equation with critical nonlinearity.''} \textit{Annals of Mathematics} 132(1990): 485-509.
 %\bibitem{mg2} M. Grillakis. {``Regularity for the wave equation with a critical nonlinearity.''} \textit{Communications on Pure and Applied Mathematics} 45(1992): 749-774.
 \bibitem{4dtypeII} M. Hillairet and P. Rapha\"{e}l. {``Smooth type II blow-up solutions to the four-dimensional energy-critical wave equation.''} \textit{Analysis} \& \textit{PDE} 5(2012), no. 4: 777-829.
 %\bibitem{fourier1} S. Helgason. {``Radon-Fourier transform on symmetric spaces and related group representations''}, \textit{Bulletin of the American Mathematical Society} 71(1965): 757-763.
 %\bibitem{fourier2} S. Helgason. \textit{Differential geometry, Lie groups, and symmetric spaces, Graduate Studies in Mathematics 34}, Providence: American Mathematical Society, 2001.
 %\bibitem{radon1} S. Helgason. {``The Radon transform on Euclidean spaces, compact two-point homogeneous spaces and Grassmann Manifolds.''} \textit{Acta Mathematica} 113(1965): 153-180.
% \bibitem{conformal} K. Hidano. {``Conformal conservation law, time decay and scattering for nonlinear wave equation''} \textit{Journal D'analysis Math\'{e}matique} 91(2003): 269-295.
 %\bibitem{SThyper} A. D. Ionescu. {``Fourier integral operators on noncompact symmetric spaces of real rank one''}, \textit{Journal of Functional Analysis} 174 (2000): 274-300.
% \bibitem{ksph2} A. D. Ionescu. {``An endpoint estimate for the Kunze-Stein phenomenon and related maximal operators''} \textit{Annals of Mathematics} 152, No.2(2000): 259-275.
 %\bibitem{hypersdg} A. D. Ionescu, and G. Staffilani. {``Semilinear Schr\"{o}dinger flows on hyperbolic spaces: scattering in $\Hm^1$''} \textit{Mathematische Annalen} 345, No.1(2012): 133-158.
  \bibitem{5dtypeII} J. Jendrej. {``Construction of type II blow-up solutions for the energy-critical wave equation in dimension 5.''} \textit{Journal of Functional Analysis} 272(2017), no. 3: 866-917. 
  %\bibitem{twobubble6d} J. Jendrej. {``Construction of two-bubble solutions for energy-critical wave equations.''} \textit{American Journal of Mathemtics} 141(2019), no.1: 55-118. 
 \bibitem{anothersoliton} J. Jendrej and A. Lawrie. {``Soliton resolution for the energy-critical nonlinear wave equation in the radial case.''} \textit{Annal of PDE} 9(2023), no. 2: Paper No. 18. 
  %\bibitem{blowup3} F. John. {``Blow-up of solutions of nonlinear wave equations in three space dimensions''} \textit{Manuscripta Mathematica} 28(1979): 235-268.
 %\bibitem{continuousL} L. Kapitanski. {``Global and unique weak solutions of nonlinear wave equations.''} \textit{Mathematical Research Letters} 1(1994): 211-223. 
\bibitem{loc1} L. Kapitanski. {``Weak and yet weaker solutions of semilinear wave equations''} \textit{Communications in Partial Differential Equations} 19(1994): 1629-1676.
%\bibitem{katayamaradiation} S. Katayama. {``Asymptotic behavior for systems of nonlinear wave equations with multiple propagation speeds in three space dimensions.''} \textit{Journal of Differential Equations} 255(2013): 120-150.
%\bibitem{kato} T. Kato. {``On nonlinear schr\"{o}dinger equations, II. $\dot{H}^s$-solutions and unconditional well-posedness.''} \textit{Journal D'analyse Math\'{e}matique} 67(1995).
 %\bibitem{endpointStrichartz} M. Keel, and T. Tao. {``Endpoint Strichartz estimates''} \textit{American Journal of Mathematics} 120 (1998): 955-980.
 \bibitem{channel5d} C. E. Kenig, A. Lawrie, B. Liu and W. Schlag. {``Relaxation of wave maps exterior to a ball to harmonic maps for all data''} \textit{Geometric and Functional Analysis} 24(2014): 610-647.
 \bibitem{randomf} C. E. Kenig and D. Mendelson. {``The focusing energy-critical nonlinear wave equation with random initial data.''} \textit{International Mathematical Research Notices} 2021, no. 19: 14508-14615. 
% \bibitem{channel} C. E. Kenig, A. Lawrie, B. Liu and W. Schlag. {``Channels of energy for the linear radial wave equation.''} \textit{Advances in Mathematics}  285(2015): 877-936.
\bibitem{kenig} C. E. Kenig and F. Merle. {``Global Well-posedness, scattering and blow-up for the energy critical focusing non-linear wave equation.''} \textit{Acta Mathematica} 201(2008): 147-212.
%\bibitem{kenig1} C. E. Kenig, and F. Merle. {``Global well-posedness, scattering and blow-up for the energy critical, focusing, non-linear Schr\"{o}dinger equation in the radial case.''} \textit{Inventiones Mathematicae} 166(2006): 645-675.
% \bibitem{km} C. E. Kenig, and F. Merle. {``Nondispersive radial solutions to energy supercritical non-linear wave equations, with applications.''} \textit{American Journal of Mathematics} 133, No 4(2011): 1029-1065.
 %\bibitem{kenig2} C. E. Kenig, and F.Merle. {``Scattering for $\dot{H}^{1/2}$ bounded solutions to the cubic, defocusing NLS in 3 dimensions.''} \textit{Transactions of the American Mathematical Society} 362(2010): 1937-1962.
% \bibitem{fchain2} C.E.Kenig, G. Ponce and L.Vega. {``Well-posedness and scattering results for the generalized Korteweg-de Vries equation via the contraction principle''}, \textit{Communications on Pure and Applied Mathematics} 46(1993): 527-620.
 %\bibitem{kv1} R. Killip, B. Stovall and M. Visan. {``Blowup Behavior for the Nonlinear Klein-Gordan Equation.''}, preprint arXiv: 1203.4886v1.
% \bibitem{kv2} R. Killip, and M. Visan. {``The defocusing energy-supercritical nonlinear wave equation in three space dimensions''} \textit{Transactions of the American Mathematical Society}, 363(2011): 3893-3934.
% \bibitem{kv3} R. Killip, and M. Visan. {``The radial defocusing energy-supercritical nonlinear wave equation in all space dimensions''} \textit{Proceedings of the American Mathematical Society}, 139(2011): 1805-1817.
 %\bibitem{kv} R. Killip, and M. Visan {``The focusing energy-critical nonlinear Schr\"{o}dinger equation in dimensions five and higher.''} \textit{American Journal of Mathematics} 132(2010): 361-424.
 %\bibitem{tao} R. Killip, T. Tao, and M. Visan. {``The cubic nonlinear Schr\"{o}dinger equation in two dimensions with redial data.''} \textit{Journal of the European Mathematical Society} 11, Issue 6(2009): 1203-1258.
% \bibitem{smallgs1} S. Klainerman, and G. Ponce. {``Global, small amplitude solutions to nonlinear evolution equations''} \textit{Communications on Pure and Applied Mathematics} 36(1983): 133-141.
%\bibitem{stablem2} J. Krieger. {``On stability of type II blow-up for the critical nonlinear wave equation on $\Rm^{3+1}$.''} \textit{Memoirs of the American Mathematical Society} 267(2020), no. 1301. 
%\bibitem{instable} J. Krieger and J. Nahas. {``Instability of type II blow up for the quintic nonlinear wave equation on $\Rm^{3+1}$.''} \textit{Bulletin De La Societe Mathematique De France} 143(2015), no. 2: 339-355.
%\bibitem{aboveground1} J. Krieger, K. Nakanishi and W. Schlag. {``Global dynamics away from the ground state for the energy-critical nonlinear wave equation.''} \textit{American Journal of Mathematics} 135(2013), no. 4: 935-965.
%\bibitem{aboveground2} J. Krieger, K. Nakanishi and W. Schlag. {``Center-stable manifold of the ground state in the energy space for the critical wave equation.''} \textit{Mathematische Annalen} 361(2015), no. 1-2: 1-50.
 \bibitem{slowblowup1} J. Krieger, W. Schlag and D. Tataru. {``Slow blow-up solutions for the $\dot{H}^1(\Rm^3)$ critical focusing semilinear wave equation.''} \textit{Duke Mathematical Journal} 147(2009), no. 1: 1-53. 
 \bibitem{slowblowup2} J. Krieger and W. Schlag. {``Full range of blow up exponents for the quintic wave equation in three dimensions.''} \textit{Journal de Math\'{e}matiques Pures et Appliqu\'{e}es} 101(2014), issue 6: 873-900. 
 %\bibitem{typeIblowup} J. Krieger and W. Wong. {``On type I blow-up formation for the critical NLW.''} \textit{Communications in Partial Differential Equations} 39(2014), no. 9: 1718-1728. 
%\bibitem{lax} P. D. Lax. {``Functional analysis.''} \textit{Pure and Applied Mathematics (New York)}, Wiley-Interscience, New York, 2002.
 \bibitem{negativeenergy} H. Levine. {``Instability and nonexistence of global solutions to nonlinear wave equations of the form $\mathbf{P} u_{tt} = - \mathbf{A} u + F(u)$."} \textit{Transactions of the American Mathematical Society} 192(1974): 1-21.
 \bibitem{shenradiation} L. Li, R. Shen and L. Wei.  {``Explicit formula of radiation fields of free waves with applications on channel of energy''}, \textit{Analysis} \& \textit{PDE}  17(2024), no. 2: 723-748.
% \bibitem{3Dnonradiativedecay} L. Li, R. Shen and C. Wang. {``An inequality regarding non-radiative linear waves via a geometric method''}, \textit{arXiv preprint} 2201.02284.
% \bibitem{nonradialCE} L. Li, R. Shen, C. Wang and L. Wei. {``Asymptotic behaviour of non-radiative solution to the wave equations.''} \textit{arXiv} 2201.02286.
 %\bibitem{morawetzsch} J. Lin, and W. Strauss. {``Decay and scattering of solutions of a nonlinear Schr\"{o}dinger equation.''} \textit{Journal of Functional Analysis} 30(1978): 245-263.
 \bibitem{ls} H. Lindblad, and C. Sogge. {``On existence and scattering with minimal regularity for semi-linear wave equations''} \textit{Journal of Functional Analysis} 130(1995): 357-426.
% \bibitem{radon2} D. Ludwig. {``The Radon transform on Euclidean space.''} \textit{Communications on Pure and Applied Mathematics} 19, no. 1(1966): 49-81.
% \bibitem{morawetz} C. S. Morawetz {``Time decay for the nonlinear Klein-Gordon equations.''} \textit{Proceedings of the Royal Society. London. Series A} 306(1968): 291-296.
 \bibitem{msoliton5d} Y. Martel and F. Merle. {``Construction of multi-solitons for the energy-critical wave equation in dimension 5.''} \textit{Archive for Rational Mechanics and Analysis} 222(2016), no. 3: 1113-1160.
 \bibitem{msoliton5ds} Y. Martel and F. Merle. {``Inelasticity of soliton collisions for the 5D energy critical wave equation.''} \textit{Inventiones Mathematicae} 214(2018), no. 3: 1267-1363. 
 \bibitem{enscatter1} K. Nakanishi. {``Unique global existence and asymptotic behaviour of solutions for wave equations with non-coercive critical nonlinearity.''} \textit{Communications in Partial Differential Equations} 24(1999): 185-221.
 \bibitem{enscatter2} K. Nakanishi. {``Scattering theory for nonlinear Klein-Gordon equations with Sobolev critical power.''} \textit{International Mathematics Research Notices} 1999, no.1: 31-60.
 %\bibitem{morawetzNLKG} K. Nakanishi. {``Energy scattering for nonlinear Klein-Gordon and Schr\"{o}dinger equations in spatial dimensions 1 and 2.''} \textit{Journal of Functional Analysis} 169(1999): 201-225.
% \bibitem{Pecher} H. Pecher. {``Nonlinear small data scattering for the wave and Klein-Gordon equation.''} \textit{Mathematische Zeitschrift} 185(1984): 261-270.
% \bibitem{benoit} B. Perthame, and L. Vega. {``Morrey-Campanato estimates for Helmholtz equations.''} \textit{Journal of Functional Analysis} 164(1999): 340-355.
% \bibitem{local1} H. Pecher. {``Nonlinear small data scattering for the wave and Klein-Gordon equation''} \textit{Mathematische Zeitschrift} 185(1984): 261-270.
 %\bibitem{sub45} C. Rodriguez. {``Scattering for radial energy-subcritical wave equations in dimensions 4 and 5.''} \textit{Communications in Partial Differential Equations} 42(2017): 852-894. 
 %\bibitem{shen1} R. Shen. {``Global well-posedness and scattering of defocusing energy subcritical nonlinear wave equation in dimension 3 with radial data.''} (2011): preprint arXiv: 1111.1234.
 \bibitem{ss1} J. Shatah and M. Struwe. {``Regularity results for nonlinear wave equations''} \textit{Annals of Mathematics} 138(1993): 503-518.
 \bibitem{ss2} J. Shatah and M. Struwe. {``Well-posedness in the energy space for semilinear wave equations with critical growth''} \textit{International Mathematics Research Notices} 7(1994): 303-309.
 \bibitem{shen2} R. Shen. {``On the energy subcritical, nonlinear wave equation in $\Rm^3$ with radial data''}  \textit{Analysis and PDE} 6(2013): 1929-1987.
% \bibitem{compactness35} R. Shen. {``A semi-linear energy critical wave equation with an application.''} \textit{Journal of Differential Equations} 261(2016): 6437-6484. 
%\bibitem{subhyper} R. Shen and G. Staffilani. {``A Semi-linear Shifted Wave Equation on the Hyperbolic Spaces with Application on a Quintic Wave Equation on $\Rm^2$''}, \textit{Transactions of the American Mathematical Society} 368(2016): 2809-2864.
 %\bibitem{energyhyper} R. Shen. {``On the energy-critical semi-linear shifted wave equation on the hyperbolic spaces''}, to appear in \textit{Differential of Integral Equations}, preprint arXiv: 1408.0331.
% \bibitem{shen3} R. Shen. {``Scattering of solutions to the defocusing energy subcritical semi-linear wave equation in 3D''} \textit{Communications in Partial Differential Equations} 42(2017): 495-518.
 %\bibitem{pushmorawetz} R. Shen. {``Morawetz Estimates Method for Scattering of Radial Energy Sub-critical Wave Equation''} \textit{arXiv Preprint} 1808.06763. 
% \bibitem{shenenergy} R. Shen. {``Energy distribution of radial solutions to energy subcritical wave equation with an application on scattering theory''} to appear in \textit{Transactions of the American Mathematical Society}.
% \bibitem{sheninward} R. Shen. {``Scattering of solutions to NLW by inward energy decay''} \textit{arXiv Preprint} 1909.01881.
% \bibitem{shen3dnonradial} R. Shen {``Inward/outward Energy Theory of Non-radial Solutions to 3D Semi-linear Wave Equation''} \textit{Advances in Mathematics} 374(2020): 107384.
% \bibitem{shenhd} R. Shen {``Inward/outward Energy Theory of Wave Equation in Higher Dimensions''} \textit{arXiv Preprint} 1912.02428.
%\bibitem{ecarbitrary} R. Shen. {``The radiation theory of radial solutions to 3D energy critical wave equations.''} \textit{arXiv preprint} 2212.03405. 
\bibitem{dynamics3d} R. Shen. {``Dynamics of 3D focusing, energy-critical wave equation with radial data.''} \textit{arXiv preprint} 2507.00391.
% \bibitem{taomeasure} T. Tao. {``An introduction to measure theory''}, \textit{Graduate Studies in Mathematics} 126(2011), AMS, Providence.
% \bibitem{counter} T. Sideris. {``Nonexistence of global solutions to semilinear wave equations in high dimensions''} \textit{Journal of Differential Equations} 52(1984): 378-406. 
%\bibitem{fchain3} G. Staffilani. {``On the generalized Korteweg-de Vries-type equations''}, \textit{Differential and Integral Equations} 10(1997): 777-796.
%\bibitem{fchain4} M.Taylor. {``Tools for PDE. Pseudo differential operators, paradifferential operators and layer potentials''}, \textit{Mathematical Surveys and Monographs} 81(2000), AMS, Providence.
 %\bibitem{staten} R. J. Stanton and P. A. Tomas. {``Expansions for spherical functions on noncompact symmetric spaces''} \textit{Acta Mathematica} 140(1978): 251-271.
 %\bibitem{strcon} W. Strauss. \textit{Nonlinear wave equations, CBMS Regional Conference Series in Mathematics, Number 73}, Providence: American Mathematical Society, 1989.
% \bibitem{struwe} M. Struwe. {``Globally regular solutions to the $u^5$ Klein-Gordon equation.''} \textit{Annali della Scuola Normale Superiore di Pisa - Classe di Scienze} 15(1988): 495-513.
 %\bibitem{pertao} T. Tao and M. Visan. {``Stability of energy-critical nonlinear Schr\"{o}dinger equations in high dimensions''} \textit{Electronic Journal of Differential Equations} 118(2005), 28.
 %\bibitem{tataru} D. Tataru. {``Strichartz estimates in the hyperbolic space and global existence for the similinear wave equation''} \textit{Transactions of the American Mathematical Society} 353(2000): 795-807.
 %\bibitem{ktsutaya} K. Tsutaya. {``Scattering theory for semilinear wave equations with small data in two space dimensions''} \textit{Transactions of the American Mathematical Society} 342, No 2(1994): 595-618.
% \bibitem{yang1} S. Yang {``Global behaviors of defocusing semilimear wave equations''} \textit{arXiv Preprint} 1908.00606.
% \bibitem{yang2d} D. Wei and S. Yang. {``On the global behaviours for defocusing semilinear wave equation in $\Rm^{1+2}$.''} \textit{arXiv Preprint} 2003.02399.
 \bibitem{msoliton5dy} X. Yuan. {``On multi-solitons for the energy-critical wave equation in dimension 5.''} \textit{Nonlinearity} 32(2019), no. 12: 5017-5048.
\end{thebibliography}
\end{document}